\documentclass[12pt]{elsarticle}

\usepackage{graphicx}
\usepackage{amssymb}
\usepackage{amsmath}
\usepackage{amsthm}

\usepackage{subfig}
\usepackage{booktabs}
\usepackage{verbatim}

\usepackage[pagewise]{lineno}

\makeatletter

\def\ps@pprintTitle{%
 \let\@oddhead\@empty
 \let\@evenhead\@empty
 \def\@oddfoot{}%
 \let\@evenfoot\@oddfoot}
 
\makeatother

\newtheorem{theorem}{Theorem}

\newtheorem{remark}{Remark}

\newtheorem{algorithm}{Algorithm}

\def\w{\omega}
\def\k{\kappa}
\def\b{\beta}
\def\r{\rho}

\def\a{\alpha}
\def\b{\beta}
\def\g{\gamma}

\def\D{\partial}
\def\u{\mathbf{u}}

\def\vv{\mathbf{v}}
\def\v{\mathbf{v}}

\def\.{\cdot}

\def\x{\mathbf{x}}
\def\y{\mathbf{y}}

\newcommand{\ub}{{\bf u}}
\newcommand{\vb}{{\bf v}}

\newcommand{\R}{\mathbb{R}}
\newcommand{\N}{\mathbb{N}}
\newcommand{\A}{\mathcal{A}}

\newcommand{\norm}[2]{\mathopen\vert\hskip-1pt\vert#1 \mathclose\vert\hskip-1pt\vert_{#2}}

\newcommand{\E}{\mathbb{E}}
\newcommand{\Var}{\mathrm{Var}}

\newcommand{\mbf}[1]{\mathbf{#1}}
\newcommand{\mbb}[1]{\mathbb{#1}}
\newcommand{\mc}[1]{\mathcal{#1}}

\newcommand{\supp}{\mathop{\mathrm{supp}}}

\newcommand{\sgn}{\mathop{\mathrm{sgn}}}

\newcommand{\hv}{\hat{\vv}}

\newcommand{\hw}{\hat{\mbf{w}}}

\newcommand{\hby}{\hat{\y}}
\newcommand{\hz}{\hat{\mbf{z}}}

\DeclareMathOperator*{\essinf}{ess\,inf}

\DeclareMathOperator*{\spa}{span}

\journal{Computer Methods in Applied Mechanics and Engineering}

\begin{document} 
\begin{frontmatter}



\title{Multiparametric shell eigenvalue problems}


\author[hh]{Harri Hakula\fnref{support}}
\fntext[support]{H. Hakula: Supported by (FP7/2007--2013) ERC grant agreement no 339380.}

\ead{Harri.Hakula@aalto.fi}

\author[hh]{Mikael Laaksonen\fnref{support2}}
\fntext[support2]{M. Laaksonen: Supported by the Magnus Ehnrooth Foundation.}
\ead{Mikael.J.Laaksonen@aalto.fi}

\address[hh]{Aalto University\\
Department of Mathematics and Systems Analysis\\
P.O. Box 11100\\
FI--00076 Aalto, Finland
}

\begin{abstract}
The eigenproblem for thin shells of revolution under uncertainty in material
parameters is discussed. 
Here the focus is on the smallest eigenpairs.
Shells of revolution have natural eigenclusters due to symmetries,
moreover, the eigenpairs depend on a deterministic parameter, the dimensionless thickness.
The stochastic subspace iteration algorithms presented here are capable of resolving the smallest eigenclusters.
In the case of random material parameters, it is possible that the eigenmodes cross 
in the stochastic parameter space. This interesting phenomenon is demonstrated via
numerical experiments. Finally, the effect of the chosen material model on the
asymptotics in relation to the deterministic parameter is shown to be negligible.
\end{abstract}


\end{frontmatter}


\section{Introduction}
\label{sec:introduction}

In many engineering problems there are uncertainties concerning
material models and domains. Stochastic finite element methods (SFEM)
have received much attention over the past decade. This, however,
has not extended to eigenvalue problems despite their importance in
many applications, including the dynamic response of structures which
is the focus of this work. In recent paper Sousedik and Elman \cite{Sousedik:2016it} cover most of the rather
limited literature available with the notable exception of
work by Andreev and Schwab \cite{andreevschwab12}, which is to our knowledge the only
mathematically rigorous analysis of the collocation approach
to stochastic eigenproblems. 

Naturally the focus has been on second order problems where
properties such as simple smallest modes can be taken
advantage of. For instance, it is well-known that the first mode
of the Laplacian does not change its sign.
Here the eigenvalue problem, free vibration of thin shells of revolution,
is of fourth order where such properties are not guaranteed.
Indeed, one feature of such eigenvalue problems that complicates the analysis is
the inevitability of repeated or tightly-clustered eigenvalues, which arise 
naturally with symmetries, and will be central in our numerical experiments. 
Even in deterministic setting when such eigenvalues are to be
approximated in practice, it is futile to try to determine whether
computed eigenvalue approximations that are very close to each other are all
approximating the same (repeated) eigenvalue, or approximating eigenvalues that
just happen to be very close to each other. Instead, one should consider
such clusters as a collective whole via subspaces.

The main result is that the stochastic subspace iteration is capable of
resolving the eigenclusters. Moreover, we demonstrate that eigenvalue
crossings occur in the parameter space which leads to the concept of
the \textit{effective smallest mode}. For simplicity consider the case
where the smallest eigenvalue is of higher multiplicity. Over the parameter
space different eigenmodes may be associated with the lowest eigenvalue.
In practice this means that over a set of manufactured specimens it is possible
to measure or observe different lowest modes due to manufacturing or material
imperfections.

Here the shells of revolution are used as the model problem since the
eigenmodes have special properties that lend themselves well to our study.
The geometry of the shell of revolution is defined by an axial profile
function which is rotated about the axis of revolution. In dimensionally 
reduced setting the midsurface of the shell $D$ has a natural parametrization
in the axial and angular directions, $D = [x_0,x_1]\times [0,2\pi]$ (peridodic).
If the material properties in the angular direction are constant, the eigenmodes
will have integer valued wavenumbers in that direction and using a suitable
ansatz the spectrum can be computed over a set of one-dimensional problems.
(For an illustration see the Figure~\ref{fig:d3d} below.)
Therefore, introducing uncertainty in the 
material properties, for instance, Young's modulus, with randomness
in the axial direction only is a special case of interest.
It should be noted that
even though the shell geometry is periodic, it is perfectly feasible to consider
material uncertainties that are not periodic. This would correspond to
a situation where a cylinder, say, is formed by rolling a cut sheet of material.

Another aspect of the chosen experimental setting is that the numerical \textit{locking}
can be controlled. In \cite{bhp} it is shown that the eigenvalue problem is subject
to locking due to angular oscillations. The 1D formulation does not include this
form of locking since the integration in the angular direction is exact -- assuming
that the material parameter varies in the axial direction only. Thus, we can calibrate 
the 2D discretization to agree with the 1D results. Of course, the rate of convergence 
cannot be optimal in 2D, but the results can be sufficiently accurate.

Another salient feature of shell problems is the role of dimensionless thickness,
which in the context of this paper can be treated as a \emph{deterministic parameter}.
It is natural to consider the asymptotic behaviour of the quantities of interest
of the stochastic eigenproblems as functions of thickness. Through carefully designed
numerical experiments we show that the material imperfections considered here do
not affect the known asymptotics for the first eigenpair, and the standard
deviation of the smallest eigenvalue decreases linearly as the thickness tends to zero.

The rest of the paper is organized as follows:
First in Section~\ref{sec:natureofcrossings} the concept of an eigenvalue crossing is
illustrated in the context of 2D Laplacian.
The shell eigenproblem  and its
stochastic extension are defined in Sections \ref{sec:shell_eigenproblem} and \ref{sec:stochastic_eigenproblem}, respectively.
The algorithms necessary (collocation and Galerkin) for the solution of the problems 
are given in Section \ref{sec:strategies}.
The numerical experiments with related analysis of the results are discussed in Section \ref{sec:numerical_experiments}.
Finally, in Section \ref{sec:conclusions} the conclusions and directions for future
research are considered.
The shell models used in numerical experiments are outlined in \ref{sec:shell_models}.

\section{On the Nature of Eigenvalue Crossings}
\label{sec:natureofcrossings}
The concept of eigenvalue crossings inevitably arises when considering parameter dependent eigenvalue problems. For instance, the second smallest eigenvalue of the Laplace operator on a symmetric domain is of multiplicity two. Thus, as in the example considered in \cite{hakulalaaksonen17} for instance, the second and third eigenvalues of a parameter dependent extension of the problem typically cross within the parameter space.

As an illustrative example we consider the following eigenvalue problem: find $\lambda \in \R$ and $u \in H^1_0(D)$ such that
\[
c^2 \Delta u = \lambda u,
\]
where $D$ is the unit circle in $\R^2$ and $c > 0$. Instead of a constant value for $c$ we consider a random field $c: \Gamma \to L^{\infty}(D)$ given by
\[
c^2(\xi) = 1 + a_1 \xi_1 + a_2 \xi_2, \quad \xi = (\xi_1, \xi_2) \in \Gamma := [-1,1]^2,
\]
where $a_1(r,\varphi) = (\cos(\pi r) + 1)/3$ and $a_2(r,\varphi) = \sin(2 \varphi)(1 - \cos(\pi r))/3$ in polar coordinates. The eigenpairs $(\lambda, u)$ now become functions of $\xi \in \Gamma$. We define an ordering of the eigenvalues so that
\begin{equation}
\label{exenum}
\lambda^{(1)}(\xi) \le \lambda^{(2)}(\xi) \le \ldots \quad \forall \xi \in \Gamma
\end{equation}
and assume the associated eigenfunctions to be normalized in $L^2(D)$ for every $\xi \in \Gamma$.

\begin{figure}[htb]
\begin{center}
\subfloat[{The eigenvalues $\lambda^{(2)}(\xi)$ and $\lambda^{(3)}(\xi)$ as a function of $\xi \in \Gamma$ (left) and as a function of $\xi_2 \in [-1,1]$ for $\xi_1=0$ (right).}]{\includegraphics[width=0.48\textwidth]{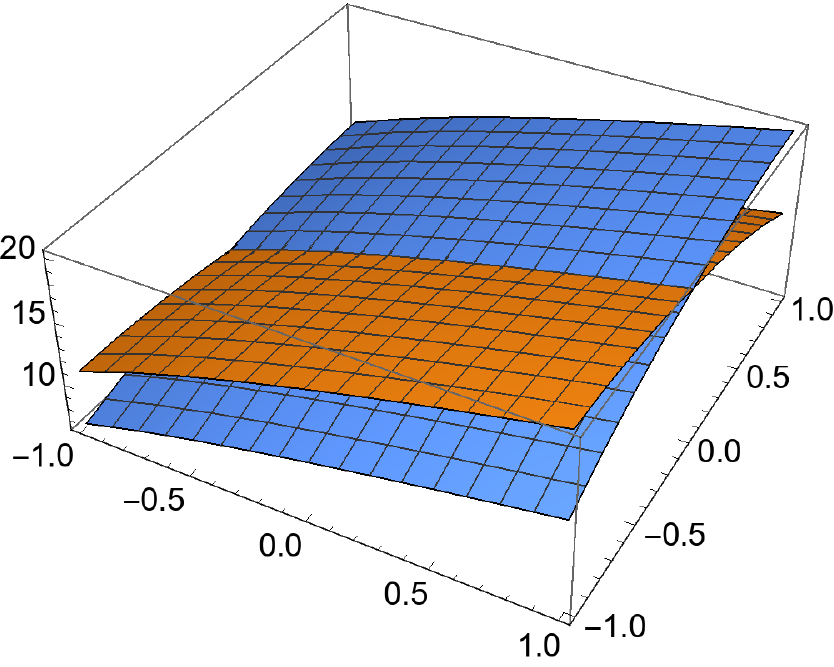} \quad \includegraphics[width=0.48\textwidth]{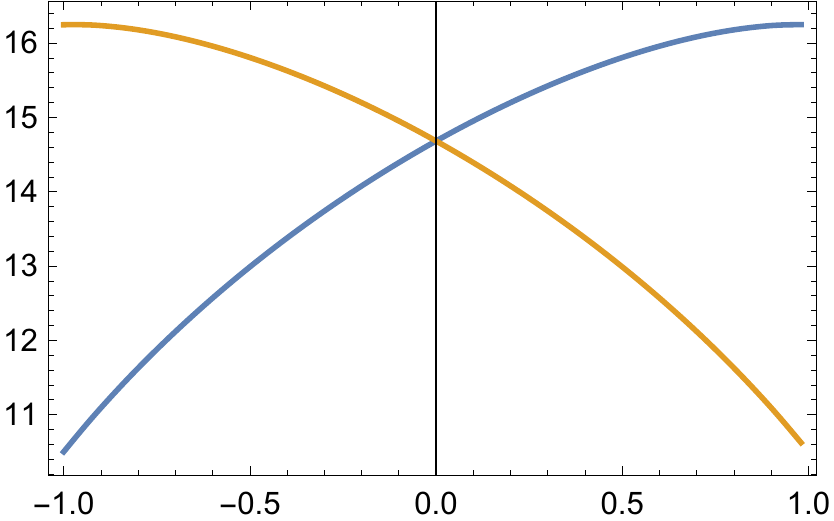}}\quad
\subfloat[{Mode $u^{(2)}(\xi)$ at $\xi=(0,-1/2)$.}]{\includegraphics[width=0.22\textwidth]{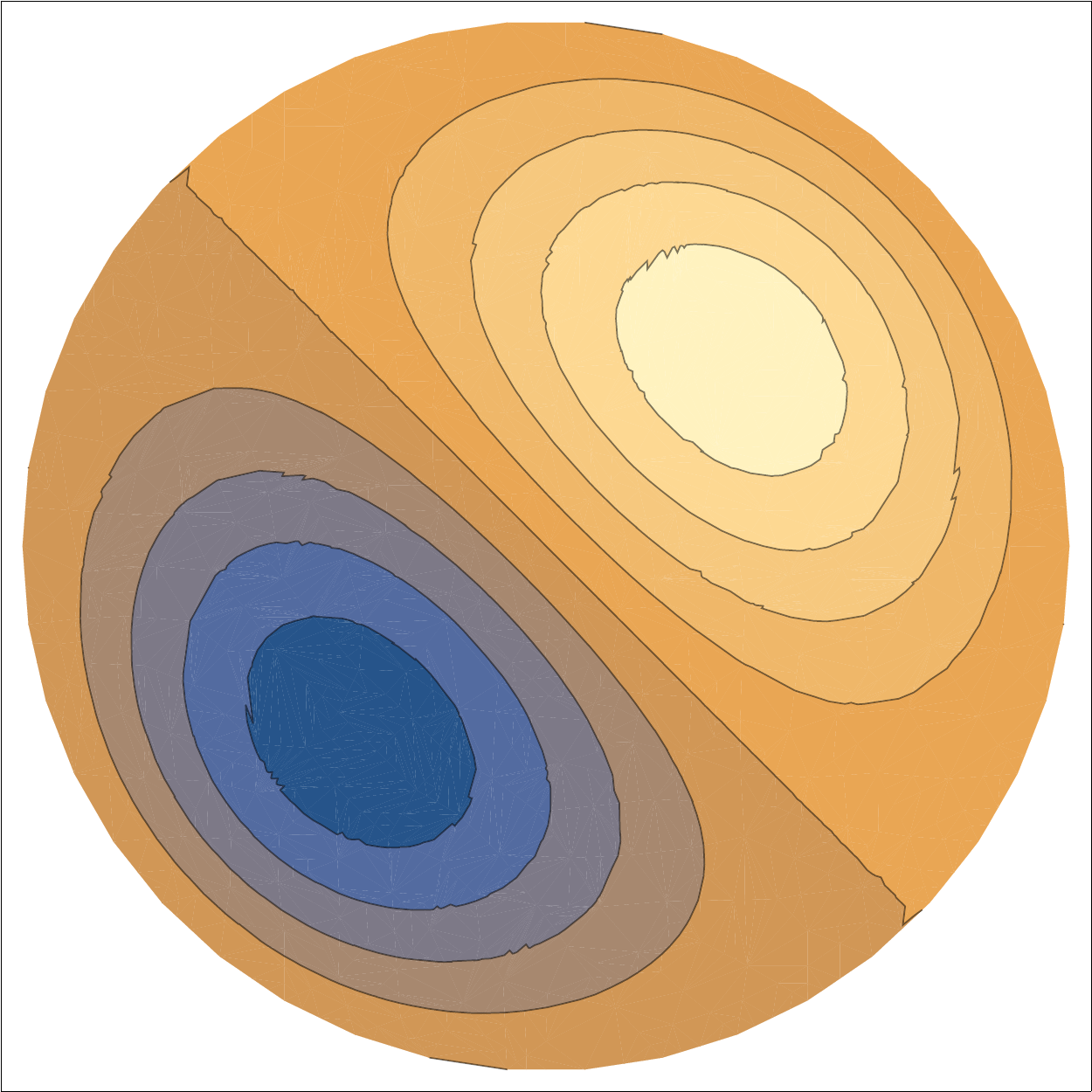}}\quad
\subfloat[{Mode $u^{(2)}(\xi)$ at $\xi=(0,1/2)$.}]{\includegraphics[width=0.22\textwidth]{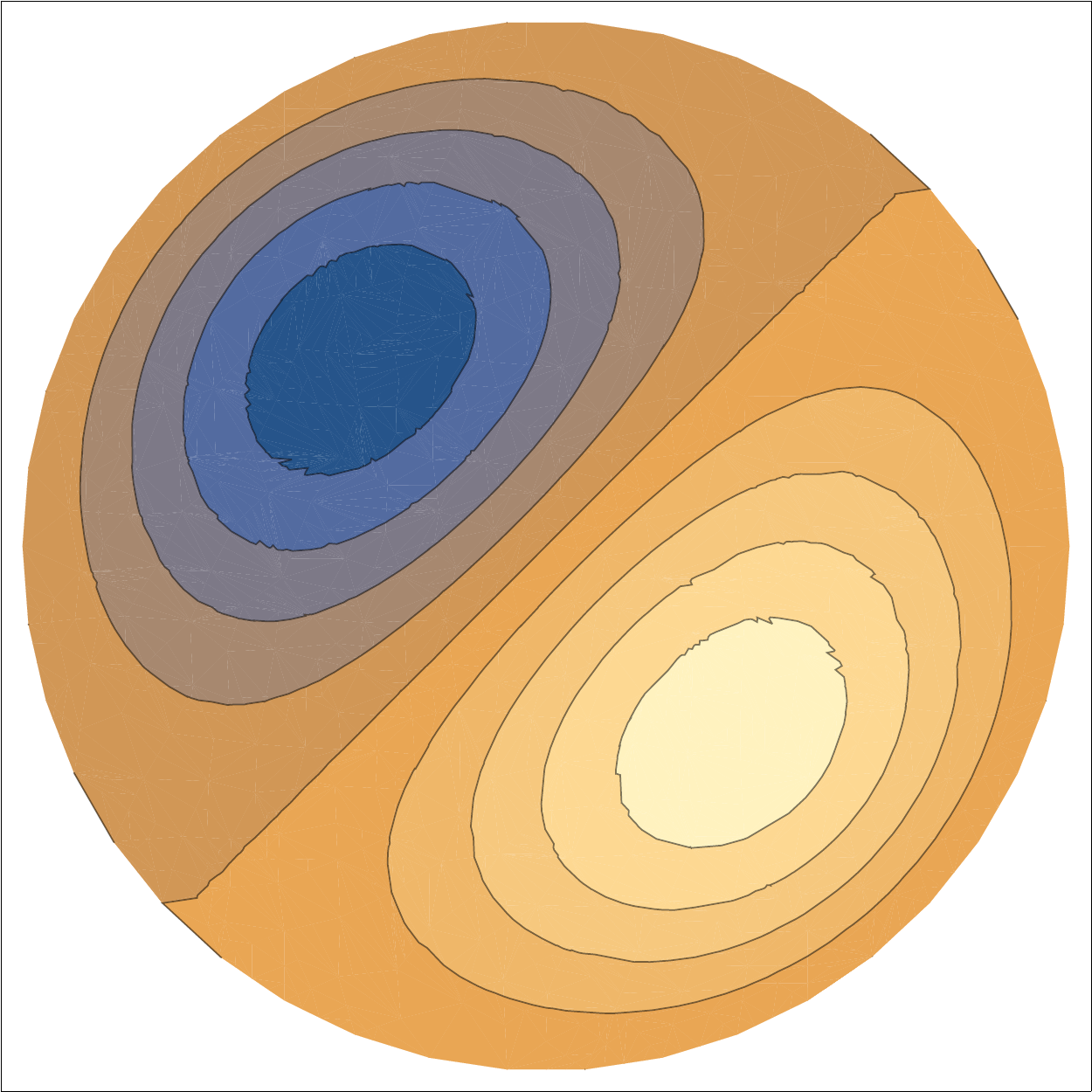}}\quad
\subfloat[{Mode $u^{(3)}(\xi)$ at $\xi=(0,-1/2)$.}]{ \includegraphics[width=0.22\textwidth]{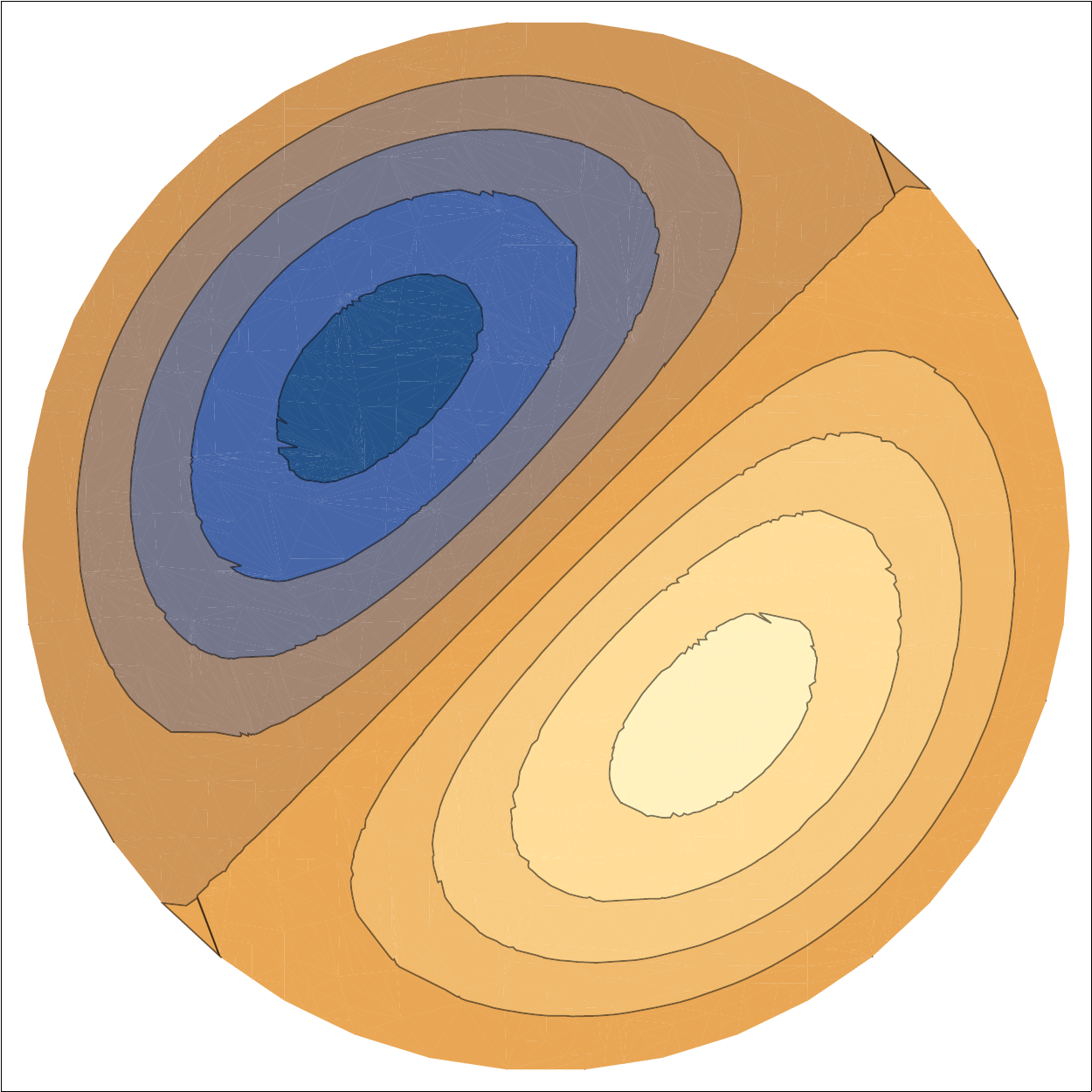}}\quad
\subfloat[{Mode $u^{(3)}(\xi)$ at $\xi=(0,1/2)$.}]{\includegraphics[width=0.22\textwidth]{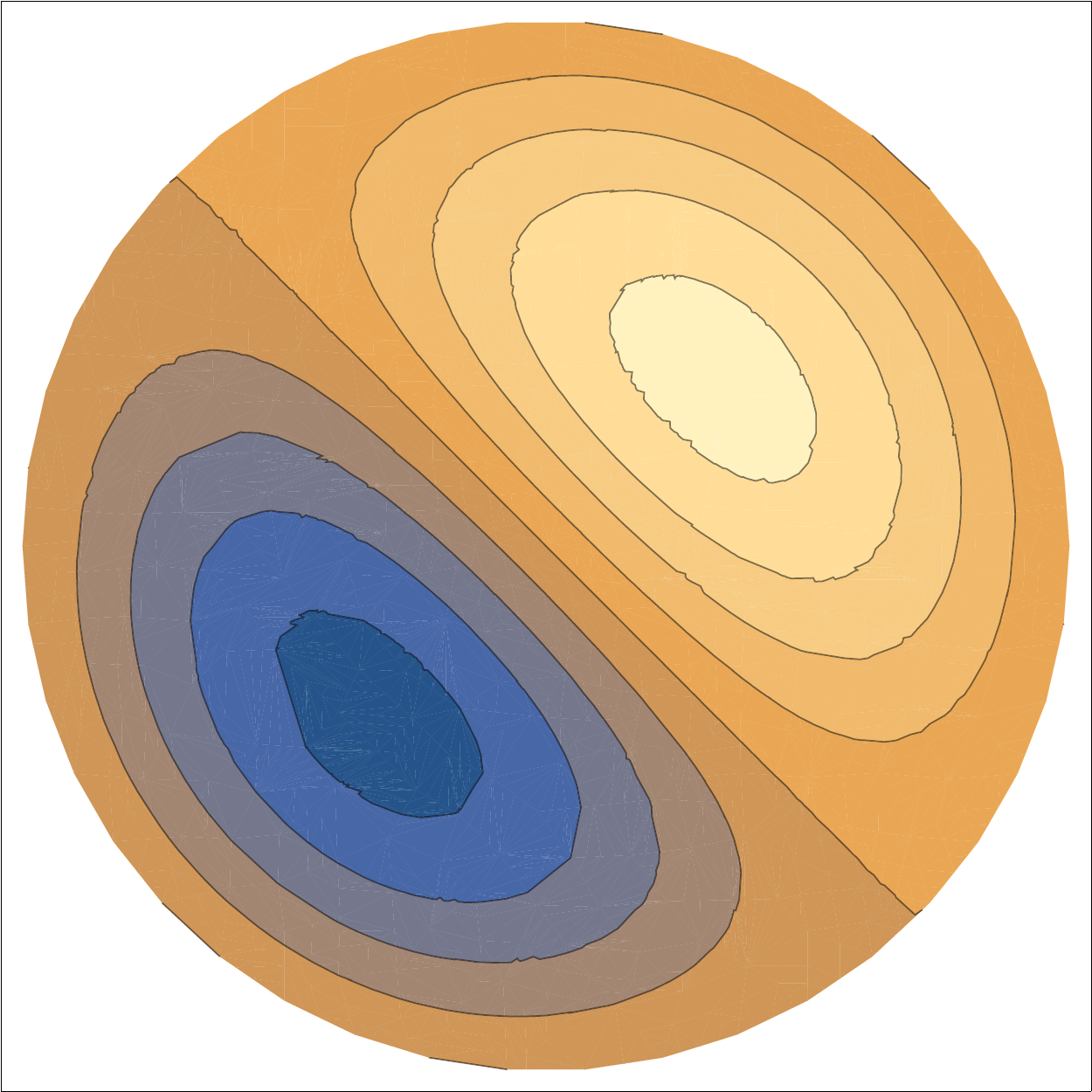}}\quad
\caption{The eigenpairs $(\lambda^{(2)}, u^{(2)})$ and $(\lambda^{(3)},u^{(3)})$ of the example problem.}
\label{fig:subspace}
\end{center}
\end{figure}

In Figure \ref{fig:subspace} we have visualised the second and third smallest eigenvalues $\lambda^{(2)}(\xi)$ and $\lambda^{(3)}(\xi)$. We see that the two eigenvalues cross at $\xi_2 = 0$ as is expected due to symmetry. Due to this eigenvalue crossing the associated eigenmodes switch places so that the eigenfunctions $u^{(2)}(\xi)$ and $u^{(3)}(\xi)$ as defined by the enumeration \eqref{exenum} are in fact discontinuous at $\xi_2 = 0$. This is problematic if we want to construct a polynomial approximation of the eigenmodes on $\Gamma$. In theory it might be possible to give up \eqref{exenum} and enumerate the eigenpairs in such a way that the eigenmodes are smooth on $\Gamma$, but in practise this might be difficult to achieve: we would have to know a priori which eigenmode to choose as our solution candidate at each point $\xi \in \Gamma$.

\section{Shell Eigenproblem}
\label{sec:shell_eigenproblem}

\begin{figure}[htb]
  \centering
  \subfloat[{Component: $w$; Does not change sign along the axial direction.}]{ \includegraphics[width=0.45\textwidth]{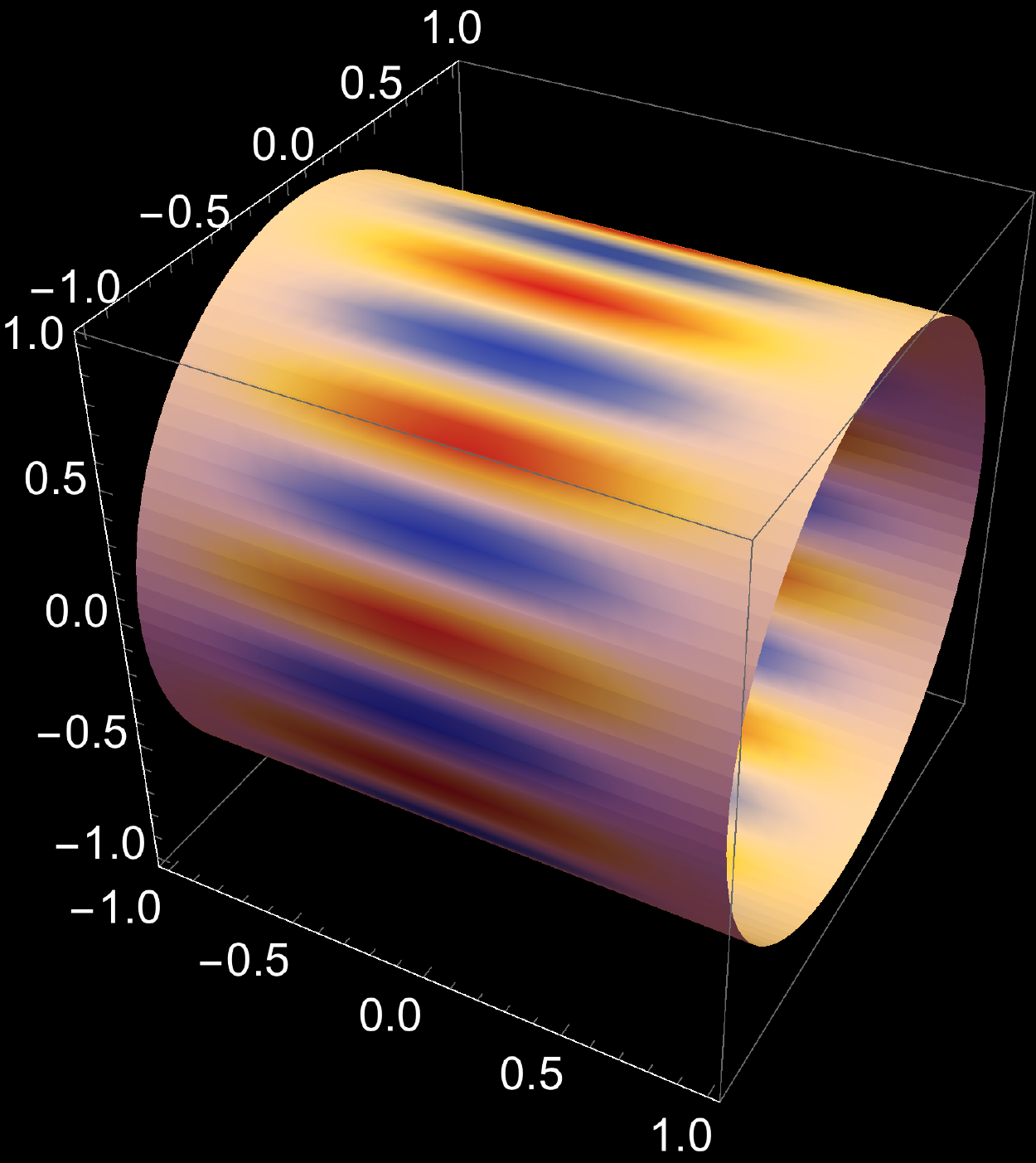}}\quad
  \subfloat[{Component: $\theta$.}]{\includegraphics[width=0.45\textwidth]{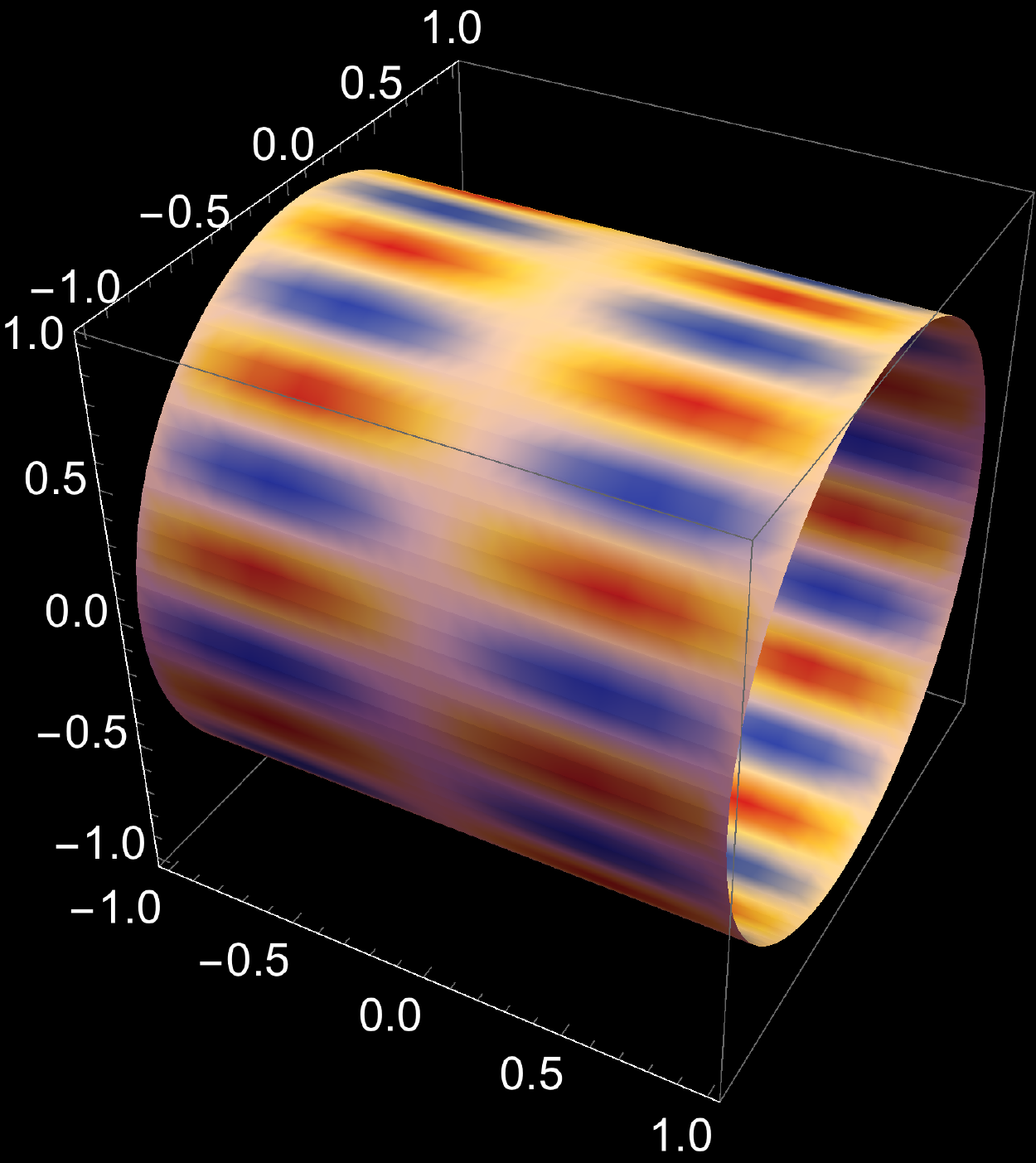}}
  \caption{Deterministic smallest eigenmode: Cylinder with radius = 1, $x\in[-1,1]$, $t=1/1000$, $k=11$. 
  Displacement indicated with temperature colour scheme.}
  \label{fig:d3d}
\end{figure}

Assuming a time harmonic displacement field, the free vibration problem for a general
shell of thickness \( t\) leads to the following 
{\em eigenvalue} problem: Find $\ub(t)$ and $\omega^2(t) \in \mathbb{R}$ such that
\begin{equation}
\label{pbl-strong}
\begin{cases}
t A_m \ub(t) + t A_s \ub(t) + t^3 A_b \ub(t) =\omega^2(t) M(t) \ub(t) \\
+ \:\: \mbox{boundary conditions.}
\end{cases}
\end{equation}

Above, $\ub(t)$ represents the shell displacement field, 
while $\omega^2(t)$ represents the square of the eigenfrequency.
The differential operators $A_m$, $A_s$ and $A_b$ account for 
membrane, shear, and bending potential energies, respectively and
are {\em independent of $t$}.
Finally, $M(t)$ is the inertia operator, which in this case can be split
into the sum $M(t) = t M^l + t^3 M^r$, with $M^l$ (displacements)
and $M^r$ (rotations) {\em independent of $t$}.
Many well-known shell models fall into this framework.

Let us next consider the variational formulation of 
problem~\eqref{pbl-strong}. Accordingly, we introduce the space
$V$ of admissible displacements, 
and consider the problem: Find:  $(\ub(t), \omega^2(t)) \in V\times\mathbb{R}$ such that
\begin{equation}
\label{pbl}
t a_m(\ub(t),\vb)+ t a_s(\ub(t),\vb)+ t^3
a_b(\ub(t),\vb)= \omega^2(t)  m(t; \ub(t),\vb) \quad \forall \vb \in V,
\end{equation}
where
\(a_m(\cdot,\cdot)\), \(a_s(\cdot,\cdot)\), \(a_b(\cdot,\cdot)\) and $m(t; \cdot,\cdot)$
are the bilinear forms associated with the operators $A_m$, $A_s$, $A_b$ and
$M(t)$, respectively.
Obviously, the space $V$ and the three bilinear forms depend on the chosen shell
model (see for instance~\cite{CB:book}). 
Here we consider two shell models, the Reissner-Naghdi and the so-called mathematical
shell model and restrict to cylindrical (parabolic) shells. 
These shell models are formally defined in the Appendix.

For $t > 0$ we denote by $H_t(D)$ the Hilbert space equipped with the inner product $m(t; \cdot,\cdot)$ and consisting of functions for which the associated norm is finite. The bilinear forms \(a_m(\cdot,\cdot)\), \(a_s(\cdot,\cdot)\), \(a_b(\cdot,\cdot)\) are continuous and elliptic. From the compact embedding of $V$ into $H_t(D)$ it follows that the problem \eqref{pbl} admits a countable number of real eigenvalues and corresponding eigenfunctions that form an orthonormal basis of $H_t(D)$ \cite{babuskaosborn91, boffi10}. We shall use the notation $(\lambda^{(i)}(t),\ub^{(i)}(t))$ for $i \in \N$ to refer to the different eigenpairs.

\subsection{Asymptotics}
Parameter-dependent asymptotics of the smallest eigenpairs 
are known rigorously for parabolic and elliptic shells, and
well-understood for hyperbolic ones. However, the inherent dynamics of the associated Fourier
modes have not been addressed in the literature, although it is clear that the practitioners
have been well aware of many of the issues. Here we give a brief outline of the central features
on cylindrical shells.

First, let us for simplicity consider just one component of the displacement field.
The transverse profile $w$ of the smallest eigenmode is parabolic, see \cite{bhp}. This means that
we can omit torsional modes from our discussion. Also, it is natural to choose as the Fourier
basis, $x \in [-\pi,\pi]$, modes, i.e., the functions of type
\begin{equation} \label{eq:mode}
  w_{mn}(x,y) = \sin (\frac m 2 x) \cos (n y) \pmod{\text{boundary layers } \hat{u}(t)},
\end{equation}
where the wavenumbers $m,n \in \mathbb{Z}$, leading to general eigenbasis 
\[
  W = \left\{ w_{mn}(x,y) \right\}_{m,n=0}^{\infty}.
\]
Similar expansions for different components can be derived.
Notice, that in (\ref{eq:mode}) the boundary layers are functions of $t$. In the case of
free vibration, the boundary layers do not carry significant amount of energy and thus,
it is meaningful to track modes in the parameter space using wavenumbers $m$ and $n$ only,
since for every component of the eigenmode they are the same even if the trigonometric basis function may differ. 
Hence, we can use notation
$(\lambda_{mn},\u_{mn})$ for a corresponding eigenpair.

For a fixed value of $t$ the eigenvalues can be ordered
\[
  \lambda_{\min}(t)=\lambda_{m_i n_i}^{(1)}(t) <\lambda_{m_j n_j}^{(2)}(t) < \ldots <\lambda_{m_q n_q}^{(k)}(t) < \ldots
\]
and naturally this is an ordering of the eigenmodes as well. For parabolic shells, it is known
that asymptotically
\begin{theorem}[{\cite{bhp},\cite{abhl1}}]
\label{asymptotics}
For cylindrical shells of revolution, the smallest eigenvalue $\lambda$ and the integer valued
wavenumber $k$ scale as functions of dimensionless thickness $t$:
\[
  \lambda \sim t, \quad k \sim t^{-1/4}.
\]
\end{theorem}

The theoretical results above imply that for some fixed $t_1$ the smallest eigenpair is of type
\[
  (\lambda_{1n_i}(t_1), \u_{1n_i}(t_1))
\]
with
\[
\lambda_{\min}(t_1) = \lambda_{1n_i}(t_1) \leq \lambda_{1n_j}(t_1),  \qquad n_j \geq n_i,
\]
but for some $t_2 < t_1$
\[
  \lambda_{\min}(t_2) = \lambda_{1n_j}(t_2) \leq \lambda_{1n_i}(t_2), \qquad n_j \geq n_i.
\]
Because of continuity, the crossing of modes can only occur at double eigenvalues.
Therefore, as the parameter $t$ tends to zero, the orderings of the modes can change continuously.
\begin{remark}
  The crossings are always sharp, since for instance energy ratios of the modes are bounded below
  by the membrane eigenmodes ($t=0$).
\end{remark}
Due to this characteristic of the shell eigenproblems, they are interesting model problems
for stochastic eigenproblems. As we allow for random material properties, 
for instance the Young's modulus $E$, it is to be expected that the orderings can be perturbed
in the parameter space and crossings will occur. In other words, for a given thickness $t$
the observed smallest eigenmodes may vary over different realizations of $E$.

\section{Stochastic Shell Eigenproblem}
\label{sec:stochastic_eigenproblem}

The stochastic reformulation of the shell eigenvalue problem arises from introducing uncertainties in the physical coefficients. The underlying assumption is that these uncertainties may be parametrized using a countable number of independent random variables. We note that special care must now be taken in defining the eigenmodes of the stochastic problem so that different realizations are in fact comparable.

\subsection{Parametrization of the random input}

We assume that the Young's modulus $E$ is a random field on the physical domain $D$. We take the conventional approach in stochastic finite elements and assume that $E$ can be written in an expansion of the form
\begin{equation}
\label{kl}
E(\x, \xi) = E_0(\x) + \sum_{m=1}^{\infty} E_m(\x) \xi_m, \quad \x \in D, \quad \xi \in \Gamma,
\end{equation}
where $\xi = (\xi_1, \xi_2, \ldots)$ represents a vector of mutually independent random variables that take values in a suitable domain $\Gamma \subset \R^{\infty}$. Typically the parametrization \eqref{kl} is assumed to result from a Karhunen-Lo\'eve expansion of the input random field, i.e., $E$ is written as a linear combination of the eigenfunctions of the associated covariance operator. 

For the sake of simplicity we assume here that the random variables $\{\xi_m\}_{m=1}^{\infty}$ are uniformly distributed. Hence, after possible rescaling, we have $\Gamma := [-1,1]^{\infty}$. We let $\mu$ denote the underlying uniform product probability measure and $L^2_{\mu}(\Gamma)$ the corresponding weighted $L^2$-space. For functions in $L^2_{\mu}(\Gamma)$ we define the expected value 
\[
\E[v] = \int_{\Gamma} v(\xi) \ d \mu(\xi)
\]
and variance $\Var[v] = \E[(v - \E[v])^2]$.

In order to guarantee positivity and boundedness of the coefficient $E$ we assume that $E_0 \in L^{\infty}(D)$ and
\[
\essinf_{\x \in D} E_0(\x) > \sum_{m=1}^{\infty} \norm{E_m}{L^{\infty}(D)}.
\]
In Section \ref{sec:numerical_experiments} we consider examples with algebraic decay of the series \eqref{kl}:
\[
\norm{E_m}{L^{\infty}(D)} \le Cm^{-\sigma}, \quad \sigma > 1, \quad m = 1, 2, \ldots
\]
For bounds on the decay of the Karhunen-Lo\'eve eigenpairs we refer to \cite{bieri09}.

\subsection{The stochastic eigenproblem}

The random nature of the Young's modulus is inherited by the underlying shell model itself. Hence, the operators $A_m$, $A_s$ and $A_b$ in the eigenproblem \eqref{pbl-strong} now depend on $\xi \in \Gamma$. The stochastic extension of the free vibration problem for a shell of thickness \( t\) reads: find $\ub(t,\cdot)$ and $\omega^2(t,\cdot)$ such that for every $\xi \in \Gamma$
\begin{equation}
\label{sproblem}
\begin{cases}
t A_m(\xi) \ub(t,\xi) + t A_s(\xi) \ub(t,\xi) + t^3 A_b(\xi) \ub(t,\xi) =\omega^2(t,\xi) M(t) \ub(t,\xi) \\
+ \:\: \mbox{boundary conditions.}
\end{cases}
\end{equation}
We assume the eigenvector $\ub(t,\xi)$ to be normalized in $H_t(D)$ for all $\xi \in \Gamma$.

The variational form the problem \eqref{sproblem} is: find functions $\ub(t, \cdot) \! : \Gamma \to V$ and $\omega^2(t, \cdot) \! : \Gamma \to \R$ such that for all $\xi \in \Gamma$ we have
\begin{align}
\label{svarproblem}
t a_m(\xi; \ub(t,\xi),\vb) + t a_s(\xi; \ub(t,\xi),\vb)+ & t^3
a_b(\xi; \ub(t,\xi),\vb) = \nonumber \\
& \omega^2(t,\xi) m(t;\ub(t,\xi),\vb) \quad \forall \vb \in V.
\end{align}
Here \(a_m(\xi; \cdot,\cdot)\), \(a_s(\xi; \cdot,\cdot)\), and \(a_b(\xi; \cdot,\cdot)\) are bilinear forms associated with the stochastic operators $A_m(\xi)$, $A_s(\xi)$ and $A_b(\xi)$ respectively.

\subsection{Stochastic subspaces}

Stochastic eigenvalue problems have proven difficult to assess mathematically and to solve numerically. The algorithms typically suggested in the literature (e.g. \cite{andreevschwab12, hakulakaarniojalaaksonen15}) restrict to simple eigenmodes only. As noted, for shells of revolution it does not make sense to assume that the smallest eigenmode is simple. In fact the ansatz given in \ref{sec:1D_models} reveals that each eigenvalue is associated to an eigenspace of dimension two, when material properties of the shell are constant in the angular direction. Moreover, eigenvalue crossings may appear, and in certain cases we observe two eigenspaces inextricably intertwining.

Consider the problem \eqref{svarproblem} for a fixed value of $t$. Assume an ordering of the eigenvalues such that
\[
\lambda^{(1)} (t,\xi) \le \lambda^{(2)} (t,\xi) \le \ldots \quad \forall \xi \in \Gamma
\]
and denote by $\ub^{(1)}(t,\xi), \ub^{(2)}(t,\xi), \ldots$ associated eigenfunctions that are orthonormal in $H_t(D)$. We refer to the eigenpair $(\lambda^{(1)}(t,\cdot), \ub^{(1)}(t,\cdot))$ as the effective smallest mode of the problem. Due to possible eigenvalue crossings the eigenpairs $\{ (\lambda^{(i)}(t,\cdot), \ub^{(i)}(t,\cdot)) \}_{i=1}^S$ are not necessarily analytic nor even continuous on $\Gamma$. Therefore, instead of considering individual eigenmodes of the problem, it often makes more sense to consider the subspace associated to a cluster of eigenvalues.

Suppose that the $S \in \N$ smallest eigenvalues are strictly separated from the rest of the spectrum, i.e.,
\[
\inf_{\xi \in \Gamma} | \lambda^{(S)}(t,\xi) - \lambda^{(S+1)}(t,\xi) | > \varepsilon, \quad \varepsilon > 0. 
\]
We aim to find a basis $\{ \v^{(i)}(t,\cdot) \}_{i=1}^S \subset V$ such that
\[
\spa \{ \ub^{(i)}(t,\xi) \}_{i=1}^S = \spa \{ \v^{(i)}(t,\xi) \}_{i=1}^S \quad \forall \xi \in \Gamma.
\]
In Section \ref{sec:strategies} we suggest algorithms for approximately computing the basis vectors $\{ \v^{(i)}(t,\cdot) \}_{i=1}^S$.

Needless to say, the basis for an eigenspace is not uniquely defined, since the vectors $\{ \v^{(i)}(t,\cdot) \}_{i=1}^S$ may be chosen in a variety of ways for each value of $\xi \in \Gamma$. However, as we illustrate in Section \ref{sec:strategies}, the basis vectors may always be anchored to a certain reference basis computed for the deterministic problem.

\section{Solution Strategies}
\label{sec:strategies}

We propose two different approaches for the spatial discretization of the variational shell eigenvalue problem \eqref{svarproblem}:
\begin{itemize}
\item [A)] Apply the ansatz given in \ref{sec:1D_models} and discretize the reduced 1D problem;
\item [B)] Discretize the full 2D problem.
\end{itemize}
In this paper we employ the $p$-version of the finite element method to form the appropriate approximation spaces.

For solving the resulting spatially discretized problems we again propose two strategies:
\begin{itemize}
\item [1)] Sample solution statistics from an ensemble of deterministic solutions computed at predefined points in $\Gamma$ (stochastic collocation);
\item [2)] Employ spectral inverse and subspace iteration algorithms to seek an approximate polynomial representation for the solution (stochastic Galerkin).
\end{itemize}
The smallest eigenvalue of the dimensionally reduced problem in case A is typically simple and the corresponding eigenmode may often be computed separately. For the full 2-dimensional problem in case B, on the other hand, the eigenvalues are clustered and therefore we need to consider a higher dimensional subspace.

\subsection{Galerkin discretization in space}

We assume standard finite element discretization spaces $V_p \subset H_0^1(D)$ of varying polynomial order $p \in \N$. This results in the approximation estimates
\[
\inf_{v_p \in V_p} \norm{v - v_p}{L^2(D)} \le Ch^{p+1} \norm{v}{H^{p+1}(D)}
\]
and
\[
\inf_{v_p \in V_p} \norm{v - v_p}{H_0^1(D)} \le Ch^p \norm{v}{H^{p+1}(D)},
\]
where $h$ is the mesh discretization parameter.

Consider the discretized variational equation
\begin{align}
\label{dvarform}
t a_m(\xi; \ub_p(t,\xi),\vb_p) + & t a_s(\xi; \ub_p(t,\xi),\vb_p) + t^3
a_b(\xi; \ub_p(t,\xi),\vb_p) = \nonumber \\
& \omega_p^2(t,\xi) m(t; \ub_p(t,\xi),\vb_p) \quad \forall \vb_p \in V_p.
\end{align}
From the theory of Galerkin approximation for variational eigenvalue problems we obtain the following bounds for the discretization error \cite{babuskaosborn91, boffi10}.

\begin{theorem}
\label{thm:herror}
For $t > 0$ and $\xi \in \Gamma$ let $\mbf{u}(t, \xi)$ be an eigenfunction of \eqref{svarproblem} associated with an eigenvalue $\lambda(t,\xi)$ of multiplicity $m$. Let $\lambda_p^{(n_1)}(t,\xi), \ldots, \lambda_p^{(n_m)}(t,\xi)$ be eigenvalues of \eqref{dvarform} that converge to $\lambda(t,\xi)$ and let $\mbf{u}_p^{(n_1)}(t, \xi), \ldots, \mbf{u}_p^{(n_m)}(t, \xi)$ denote associated eigenfunctions. If the components of $\mbf{u}(t, \xi)$ are in $H^{1+p}(D)$, then there exists $C > 0$ and $\mbf{u}_p(t, \xi) \in \spa \{ \mbf{u}_p^{(n_i)}(t,\xi) \}_{i=1}^m$ such that
\[
| \lambda(t,\xi) - \lambda_p^{(n_i)}(t,\xi) | \le C h^{2p}, \quad i = 1,\ldots, m
\]
and
\[
\norm{\mbf{u}(t, \xi) - \mbf{u}_p(t, \xi)}{H_t(D)} \le C h^{1+p}
\]
as $p \to \infty$.
\end{theorem}

For a fixed $t > 0$ the spatially discretized variational form \eqref{dvarform} can be written as a parametric matrix eigenvalue problem: find $\lambda_p \!: \Gamma \to \mbb{R}$ and $\y_p \!: \Gamma \to \mbb{R}^N$ such that
\begin{equation}
\label{matrixevp}
\mbf{K}(\xi) \y_p (\xi) = \lambda_p(\xi) \mbf{M} \y_p (\xi) \quad \forall \xi \in \Gamma,
\end{equation}
where $N = \dim{V_p}$. Here $\mbf{M}$ is the classical mass matrix whereas
\[
\mbf{K}(\xi) = \mbf{K}^{(0)} + \sum_{m = 1}^{\infty} \mbf{K}^{(m)} \xi_m
\]
is a stochastic stiffness matrix whose each term corresponds to a term in the series \eqref{kl}. For notational simplicity we have now omitted the inherent dependence on $t$. For any fixed $\xi \in \Gamma$ the problem \eqref{matrixevp} reduces to a positive-definite generalized matrix eigenvalue problem. 

\subsection{The eigenspace of interest}

Denote by $\langle \cdot, \cdot \rangle_{\R^N_\mbf{M}}$ the inner product induced by $\mbf{M}$ and by $\norm{ \ \cdot \ }{\R^N_\mbf{M}}$ the associated norm. We define an ordering
\begin{equation}
\label{denumeration}
\lambda^{(1)}_p(\xi) \le \lambda^{(2)}_p(\xi) \le \ldots \le \lambda^{(N)}_p(\xi) \quad \forall \xi \in \Gamma
\end{equation}
for the eigenvalues of \eqref{matrixevp} and let $\y_p^{(1)}(\xi), \y_p^{(2)}(\xi), \ldots, \y_p^{(N)}(\xi)$ denote associated eigenvectors that are orthonormal in $\langle \cdot, \cdot \rangle_{\R^N_\mbf{M}}$. In the context of this paper we assume that the eigenspace of interest is the subspace corresponding to the eigenvalues $\{ \lambda^{(i)}_p \}_{i=1}^S$ for some $S \in \N$. The underlying assumption is that these eigenvalues are sufficiently well separated from the rest of the spectrum, i.e.,
\begin{equation}
\label{separated}
\inf_{\xi \in \Gamma} | \lambda^{(S)}_p(\xi) - \lambda^{(S+1)}_p(\xi) | > \varepsilon
\end{equation}
with some $\varepsilon > 0$.

\begin{remark}
Note that in the 1D ansatz we fix the angular components of the solution before discretization. Therefore the eigenvalues of the discrete problem \eqref{dvarform} only represent a subset of the eigenvalues of the continuous problem \eqref{svarproblem}.
\end{remark}

Assume a subspace
\[
U(\xi) := \spa \{ \mbf{b}^{(i)}(\xi) \}_{i=1}^S \quad \forall \xi \in \Gamma
\]
spanned by a set of vectors $\mbf{b}^{(i)} \! : \Gamma \to \R^N$. In order to obtain a unique representation of this subspace, we choose a reference basis $\{ \bar{\y}^{(i)} \}_{i=1}^S \subset \R^N$ and define the projections
\begin{equation}
\label{projections}
P_U(\bar{\y}^{(i)})(\xi) := \sum_{j=1}^S \Pi_{ij}(\xi) \mathbf{b}^{(j)}(\xi) \quad i = 1, \ldots, S
\end{equation}
with $\Pi_{ij}(\xi) = \langle \bar{\y}^{(i)}, \mathbf{b}^{(j)}(\xi) \rangle_{\R^N_{\mbf{M}}}$. Note that for each $\xi \in \Gamma$ these projected vectors only depend on the subspace $U(\xi)$ and not on the chosen basis $\{ \mathbf{b}^{(i)}(\xi) \}_{i=1}^S$. Moreover, if the matrix $\Pi(\xi) \in \R^{S \times S}$ is nonsingular, then the vectors form a basis for the space $U(\xi)$. Gram-Schmidt algorithm could be applied pointwise for $\xi \in \Gamma$ in order to make this basis orthonormal.

\subsection{Stochastic collocation}

We introduce an anisotropic sparse grid collocation operator defined with respect to a finite multi-index set, see e.g. \cite{andreevschwab12}. Let $(\mbb{N}_0^{\infty})_c$ denote the set of all multi-indices with finite support, i.e.,
\[
(\mbb{N}_0^{\infty})_c := \{ \a \in \mbb{N}_0^{\infty} \ | \ \# \supp(\a) < \infty \},
\]
where $\supp(\a) = \{ m \in \mbb{N} \ | \ \a_m \not= 0 \}$. Given a finite $\mc{A} \subset (\mbb{N}_0^{\infty})_c$ we define the greatest active dimension $M_{\mc{A}} := \max \{ m \in \N \ | \ \exists \a \in \mc{A} \textrm{ s.t. } \a_m \not= 0 \}$. For $\a, \b \in \mc{A}$ we write $\a \le \b$ if $\a_m \le \b_m$ for all $m \ge 1$. We call the multi-index set $\mc{A}$ monotone if whenever $\b \in (\mbb{N}_0^{\infty})_c$ is such that $\b \le \a$ for some $\a \in \mc{A}$, then $\b \in \mc{A}$.

Let $L_p$ be the univariate Legendre polynomial of degree $p$. Denote by $\{ \chi_k^{(p)} \}_{k=0}^p$ the zeros of $L_{p+1}$. We define the one-dimensional Lagrange interpolation operators $\mc{I}^{(m)}_p$ via
\[
(\mc{I}^{(m)}_p v)(\xi_m) = \sum_{k=0}^p v(\chi_k^{(p)}) \ell_k^{(p)} (\xi_m),
\]
where $\{ \ell_k^{(p)} \}_{k=0}^p$ are the related Lagrange basis polynomials of degree $p$. Given a finite set $\mc{A} \subset (\mbb{N}_0^{\infty})_c$ we may now define the sparse collocation operator as
\begin{equation}
\label{colloperator}
\mc{I}_{\mc{A}} := \sum_{\alpha \in \mc{A}} \bigotimes_{m \in \supp{\alpha}} \left( \mc{I}_{\alpha_m}^{(m)} - \mc{I}_{\alpha_m - 1}^{(m)} \right).
\end{equation}

The operator \eqref{colloperator} may be rewritten in a computationally more convenient form
\begin{equation}
\label{interpolationreformulation}
\mc{I}_{\mc{A}} = \sum_{\alpha \in \mc{A}} \sum_{\g \in \mc{G}_{\a}} (-1)^{\norm{\alpha - \gamma}{1}} \bigotimes_{m \in \supp(\gamma)} \mc{I}_{\gamma_m}^{(m)},
\end{equation}
where $\mc{G}_{\a} := \{ \gamma \in (\mbb{N}^{\infty}_0)_c \ | \ \a - 1 \le \g \le \a \}$. We see that the complete grid of collocation points is now given by
\[
X_{\mc{A}} := \bigcup_{\a \in \mc{A}} \bigcup_{\g \in \mc{G}_{\a}} \prod_{m \ge 1} \{ \chi_k^{(\g_m)} \}_{k=0}^{\g_m}.
\]
Observe that for every $\chi \in X_{\mc{A}}$ we have $\chi_m = 0$ when $m > M_{\mc{A}}$. For monotone multi-index sets we have
\[
X_{\mc{A}} = \bigcup_{\a \in \mc{A}} \prod_{m \ge 1} \{ \chi_k^{(\a_m)} \}_{k=0}^{\a_m}
\]
and the number of collocation points admits the bound
\[
\# X_{\mc{A}} = \sum_{\a \in \mc{A}} \prod_{m \ge 1} (\a_m + 1) \le (\# \mc{A})^2
\]
as shown in \cite{andreevschwab12}.

Assume first that \eqref{separated} holds with $S=1$. Fix a finite set $\mc{A} \subset (\mbb{N}_0^{\infty})_c$. We propose Algorithm \ref{alg:ssc} for computing the first eigenmode of the problem \eqref{matrixevp}.
\begin{algorithm}[Stochastic collocation for simple eigenvalues]
\label{alg:ssc}
Consider the multiparametric eigenvalue problem \eqref{matrixevp} with the enumeration \eqref{denumeration}. Compute a reference vector $\bar{\y}^{(1)} = \y^{(1)}_p(0)$. For each $\chi \in X_\mc{A}$ do
\begin{enumerate}
\item[(1)] Solve the problem at $\xi = \chi$ for the eigenpair $( \lambda^{(1)}_p(\chi), \y^{(1)}_p(\chi) )$.
\item[(2)] Set $\mbf{p}^{(1)}(\chi) = \sgn \left( \langle \bar{\y}^{(1)}, \y^{(1)}_p(\chi) \rangle_{\R^N_{\mbf{M}}} \right) \y^{(1)}_p(\chi)$.
\end{enumerate}
Return $(\mc{I}_{\mc{A}}(\lambda^{(1)}_p),\mc{I}_{\mc{A}}(\mbf{p}^{(1)}))$ as the approximate eigenpair.
\end{algorithm}

If the first eigenmode is not simple or strictly separated, then we may try to aim for a higher dimensional subspace. In general this means that we have to give up information on the individual eigenmodes. Assume now that \eqref{separated} holds with $S \ge 2$. We propose Algorithm \ref{alg:msc} for computing the eigenspace of interest for the problem \eqref{matrixevp}.

\begin{algorithm}[Stochastic collocation for subspaces]
\label{alg:msc}
Consider the multiparametric eigenvalue problem \eqref{matrixevp} with the enumeration \eqref{denumeration}. Compute a reference basis $\bar{\y}^{(i)} = \y^{(i)}_p(0)$ for $i = 1, \ldots, S$. For each $\chi \in X_\mc{A}$ do
\begin{enumerate}
\item[(1)] Solve the eigensystem \eqref{matrixevp} at $\xi = \chi$ for the eigenvectors $\{ \y^{(i)}_p(\chi) \}_{i=1}^S$.
\item[(2)] Let $U(\chi) = \spa \{ \y^{(i)}_p(\chi) \}_{i=1}^S$ and compute the projected vectors $\mbf{p}^{(i)}(\chi) = P_U(\bar{\y}^{(i)})(\chi)$ for $i = 1,\ldots, S$.
\end{enumerate}
Return $\{ \mc{I}_{\mc{A}}(\mbf{p}^{(i)}) \}_{i=1}^S$ as a basis for the approximate eigenspace.
\end{algorithm}

\begin{remark}
Typically we might be interested in the statistics of the solution rather than its explicit form. Applying the one-dimensional Gauss-Legendre quadrature rules on the components of \eqref{interpolationreformulation} yields formulas for the expected value as well as higher moments of the solution.
\end{remark}

\subsection{Spectral inverse iterations}

In view of stochastic Galerkin methods we construct a suitable basis of orthogonal polynomials. We then consider the spectral inverse iteration and its subspace iteration variant for computing the coefficients of the solution in the constructed basis. We refer to \cite{hakulakaarniojalaaksonen15, hakulalaaksonen17} for a detailed description of the algorithms.

We define multivariate Legendre polynomials
\[
\Lambda_{\a}(\xi) := \prod_{m \in \supp{\a}} L_{\a_m}(\xi_m), \quad \a \in (\mbb{N}_0^{\infty})_c
\]
with the normalization $\E[\Lambda_{\a}^2] = 1$. The system $\{ \Lambda_{\a} \ | \ \a \in (\mbb{N}_0^{\infty})_c \}$ forms an orthonormal basis of $L^2_{\mu}(\Gamma)$. Therefore, we may write any $v \in L^2_{\mu}(\Gamma)$ in a series
\[
v(\xi) = \sum_{\a \in (\mbb{N}_0^{\infty})_c} v_{\a} \Lambda_{\a} (\xi),
\]
where the expansion coefficients are given by $v_{\a} = \E[v \Lambda_{\a}]$.

Due to orthogonality of the Legendre polynomials we have $\E[\Lambda_{\a}] = \delta_{\a 0}$ and $\E[\Lambda_{\a} \Lambda_{\beta}] = \delta_{\a \beta}$ for all $\a, \beta \in (\mbb{N}_0^{\infty})_c$. We set
\begin{align*}
c_{\a \beta \gamma} & := \E[\Lambda_{\a} \Lambda_{\beta} \Lambda_{\gamma}], \quad \a, \beta, \gamma \in (\mbb{N}_0^{\infty})_c \\
c_{m \a \beta} & := \E[\xi_m \Lambda_{\a} \Lambda_{\beta}], \quad m \in \mbb{N}, \quad \a, \beta \in (\mbb{N}_0^{\infty})_c \\
c_{0 \a \beta} & := \delta_{\a \b}, \quad \a, \beta \in (\mbb{N}_0^{\infty})_c
\end{align*}
and define moment matrices $G^{(m)} \in \R^{P \times P}$ for $m \in \N_0$ and $G^{(\a)} \in \R^{P \times P}$ for $\a \in \mc{A}$ via $[G^{(m)}]_{\a \b} = c_{m \a \b}$ and $[G^{(\a)}]_{\b \g} = c_{\a \b \g}$.

\begin{remark}
Functions expressed in a truncated series
\[
v_{\mc{A}}(\xi) = \sum_{\a \in \mc{A}} v_{\a} \Lambda_{\a} (\xi)
\]
admit easy formulas for the mean and variance:
\[
\E[v_{\mc{A}}] = v_0, \quad \Var[v_{\mc{A}}] = \left( \sum_{\a \in \mc{A}} v_{\a}^2 \right) - v_0^2.
\]
\end{remark}

Fix a finite set of multi-indices $\mc{A} \subset (\mbb{N}_0^{\infty})_c$ and set $P= \#\mc{A}$. Denote $J := \{1,2,\ldots,N\}$. Given $\hat{s} = \{ s_{\a} \}_{\a \in \mc{A}} \in \R^P$ and $\hv = \{ v_{\a i} \}_{\a \in \mc{A}, i \in J} \in \R^{PN}$ we set
\[
\mc{P}_{\mc{A}}(\hat{s}) := \sum_{\a \in \mc{A}} s_{\a} \Lambda_{\a}
\]
and
\[
\mc{P}_{\mc{A}}(\hv) := \left\{ \sum_{\a \in \mc{A}} v_{\a i} \Lambda_{\a} \right\}_{i \in J}.
\]
We define the matrices
\begin{align*}
\Delta (\hat{s}) & := \sum_{\a \in \mc{A}} G^{(\a)} s_{\a}, \\
\widehat{\mbf{K}} & := \sum_{m=0}^{M_{\mc{A}}} G^{(m)} \otimes \mbf{K}^{(m)}, \\
\widehat{\mbf{M}} & := I_P \otimes \mbf{M}, \\
\mbf{T}(\hat{s}) & := \Delta(\hat{s}) \otimes \mbf{I}_N,
\end{align*}
where and $I_P \in \R^{P \times P}$ and $\mbf{I}_N \in \R^{N \times N}$ are identity matrices. Moreover, we define the nonlinear function $F \! : \R^{P} \times \R^{PN} \to \R^P$ via
\[
F_{\a}(\hat{s}, \hv) := \hat{s} \cdot G^{(\a)} \hat{s} - \hv \cdot (G^{(\a)} \otimes \mbf{M}) \hv, \quad \a \in \mc{A}
\]
and let $F^v \! : \R^{PN} \times \R^{PN} \to \R^P$ denote the associated bilinear form given by $F^v_{\a}(\hv,\hw) := \hv \cdot (G^{(\a)} \otimes \mbf{M}) \hw$.

The spectral inverse iteration is now defined in Algorithm \ref{alg:sii}. Assuming that \eqref{separated} holds for $S = 1$, we expect the algorithm to converge to an approximation of the first eigenpair of the system.

\begin{algorithm}[Spectral inverse iteration]
\label{alg:sii}
Consider the multiparametric eigenvalue problem \eqref{matrixevp} with the enumeration \eqref{denumeration}. Fix $tol > 0$ and let $\hby^{(0)} = \{ y^{(0)}_{\a i} \}_{\a \in \mc{A}, i \in J} \in \R^{PN}$ be an initial guess for the eigenvector. For $k = 1,2,\ldots$ do
\begin{enumerate}
\item[(1)] Solve $\hz = \{ z_{\a i} \}_{\a \in \mc{A}, i \in J} \in \R^{PN}$ from the linear system
\begin{equation}
\label{eq:linearsys}
\widehat{\mbf{K}} \hz = \widehat{\mbf{M}} \hby^{(k-1)}.
\end{equation}
\item[(2)] Solve $\hat{s} = \{ s_{\a} \}_{\a \in \mc{A}} \in \R^P$ from the nonlinear system
\begin{equation}
\label{eq:nonlinearsys}
F(\hat{s},\hz) = 0
\end{equation}
with the initial guess $s_{\a} = \norm{ \hz }{\R^P \otimes \R^N_{\mbf{M}}} \delta_{\a 0}$ for $\a \in \mc{A}$.
\item[(3)] Solve $\hby^{(k)} = \{ y_{\a i}^{(k)} \}_{\a \in \mc{A}, i \in J} \in \R^{PN}$ from the linear system
\[
\mbf{T}(\hat{s}) \hby^{(k)} = \hz.
\]
\item[(4)] Stop if $\delta^{(k)} := \norm{ \hby^{(k)} - \hby^{(k-1)} }{\R^P \otimes \R^N_{\mbf{M}}} < tol$.
\end{enumerate}
Solve $\hat{\lambda}^{(k)} \in \R^P$ from the equation
\begin{equation}
\label{eq:evalue}
\Delta(\hat{s}) \hat{\lambda}^{(k)} = \hat{e}_1,
\end{equation}
where $\hat{e}_1 = \{ \delta_{\a 0} \}_{\a \in \mc{A}} \in \R^P$. Return $(\mc{P}_{\mc{A}} (\hat{\lambda}^{(k)}), \mc{P}_{\mc{A}}(\hby^{(k)}))$ as the approximate eigenpair.
\end{algorithm}

The asymptotic convergence of Algorithm \ref{alg:sii} has been studied in \cite{hakulalaaksonen17}. In particular we have the following result.

\begin{theorem}
\label{thm:iterconv}
Let $\hby \in \R^{PN}$ be a fixed point of the Algorithm \ref{alg:sii}. Let $(\hat{s}, \hz) \in \R^P \times \R^{PN}$ be such that the equations \eqref{eq:linearsys} and \eqref{eq:nonlinearsys} hold. Assume that $\Delta(\hat{s})$ is invertible and let $\hat{\lambda}$ denote the corresponding solution to \eqref{eq:evalue}. Then there exists $\Theta = \Theta(\hat{s}, \hz) \in \R$ such that the iterates of Algorithm \ref{alg:sii} satisfy
\[
\norm{ \hby^{(k)} - \hby }{\R^P \otimes \R^N_{\mbf{M}}} \le \Theta \norm{ \hby^{(k-1)} - \hby }{\R^P \otimes \R^N_{\mbf{M}}}, \quad k \in \N
\]
for $\hby^{(k)}$ sufficiently close to $\hby$. Furthermore, there exists $C > 0$ such that
\[
\norm{ \hat{\lambda}^{(k)} - \hat{\lambda} }{\R^P} \le C \norm{ \hby^{(k)} - \hby }{\R^P \otimes \R^N_{\mbf{M}}}, \quad k \in \N.
\]
\end{theorem}

In \cite{hakulalaaksonen17} it has been argued that the value of $\Theta$ in Theorem \ref{thm:iterconv} is bounded by some $\Theta^* \in \R$ such that
\[
\Theta^* \approx \frac{\sup_{\xi \in \Gamma} \lambda_p^{(1)}(\xi)}{\inf_{\xi \in \Gamma} \lambda_p^{(2)}(\xi)}.
\]
Hence, the speed of convergence of Algorithm \ref{alg:sii} is characterized by what is essentially the largest ratio of the two smallest eigenvalues of the problem \eqref{matrixevp}.

Again, if the first eigenmode is not simple or strictly separated, then we may try to compute a higher dimensional subspace. To this end we extend the previous algorithm to a spectral subspace iteration. As before, information on the individual eigenmodes is in general lost in the process.

A Galerkin projection of the equation \eqref{projections} to the basis $\{ \Lambda_{\a} \}_{\a \in \mc{A}}$ leads to approximate projections, which can be used to measure the convergence of our spectral subspace iteration algorithm. Let $\{ \bar{\y}^{(i)} \}_{i=1}^S \subset \R^N$ be a reference basis and
\[
U(\xi) = \spa \left\{ \mc{P}_{\mc{A}}(\hat{\mbf{b}}^{(j)})(\xi) \right\}_{j=1}^S,
\]
a subspace spanned by a set of vectors $\hat{\mbf{b}}^{(j)} = \{ b_{\a i}^{(j)} \}_{\a \in \mc{A}, i \in J} \in \R^{PN}$. We define the approximate projections
\[
\hat{P}_U(\bar{\y}^{(i)}) = \sum_{j=1}^S \hat{\Pi}_{ij} \hat{\mbf{b}}^{(j)} \quad i = 1, \ldots, S,
\]
where
\[
\hat{\Pi}_{ij} = \sum_{\a \in \mc{A}} \langle \bar{\y}^{(i)}, \mbf{b}^{(j)}_{\a} \rangle_{\R^N_{\mbf{M}}} G^{(\a)} \otimes \mbf{I}_N
\]
and $\mbf{b}_{\a}^{(j)} := \{ b^{(j)}_{\a i }\}_{i \in J}$.

Assume that \eqref{separated} holds for some $S \ge 2$. We propose Algorithm \ref{alg:ssi} for computing the eigenspace of interest for the problem \eqref{matrixevp}.

\begin{algorithm}[Spectral subspace iteration]
\label{alg:ssi}
Consider the multiparametric eigenvalue problem \eqref{matrixevp} with the enumeration \eqref{denumeration}. Compute a reference basis $\bar{\y}^{(i)} = \y^{(i)}_p(0)$ for $i = 1, \ldots, S$. Fix $tol > 0$ and let $\hby^{(0,i)} = \{ y^{(0,i)}_{\a i} \}_{\a \in \mc{A}, i \in J} \subset \R^{PN}$ for $i = 1, \ldots, S$ be an initial guess for the basis of the subspace. For $k = 1,2,\ldots$ do
\begin{enumerate}
\item[(1)] For each $i = 1, \ldots, S$ solve $\hz^{(i)} = \{ z^{(i)}_{\a i} \}_{\a \in \mc{A}, i \in J}\in \R^{PN}$ from the linear system
\[
\widehat{\mbf{K}} \hz^{(i)} = \widehat{\mbf{M}} \hby^{(k-1,i)}.
\]
\item[(2)] For $i = 1, \ldots, S$ do
\begin{enumerate}
\item[(2.0)] Set
\[
\hz^{(1)} = \sum_{i=1}^S \hz^{(i)}.
\]
\item[(2.1)] Set
\[
\hw^{(i)} = \hz^{(i)} - \sum_{j=1}^{i-1} \mbf{T}\left( F^v(\hz^{(i)},\hby^{(k,j)}) \right)\hby^{(k,j)}.
\]
\item[(2.2)] Solve $\hat{s}^{(i)} = \{ s_{\a}^{(i)} \}_{\a \in \mc{A}} \in \R^P$ from the nonlinear system
\[
F(\hat{s}^{(i)},\hw^{(i)}) = 0
\]
with the initial guess $s^{(i)}_{\a} = \norm{ \hw^{(i)} }{\R^P \otimes \R^N_{\mbf{M}}} \delta_{\a 0}$ for $\a \in \mc{A}$.
\item[(2.3)] Solve $\hby^{(k,i)} = \{ y^{(k,i)}_{\a i} \}_{\a \in \mc{A}, i \in J} \in \R^{PN}$ from the linear system
\[
\mbf{T}(\hat{s}^{(i)}) \hby^{(k,i)} = \hw^{(i)}.
\]
\end{enumerate}
\item[(3)] Let $U(\xi) = \spa \{ \mc{P}_{\mc{A}} (\hby^{(i)})(\xi) \}_{i=1}^S$. Compute the projected vectors $\hat{\mbf{p}}^{(k,i)} = \hat{P}_U(\bar{\y}^{(i)})$ for $i = 1,\ldots, S$.
\item[(4)] Stop if
\[
\delta^{(k)} := \left( \sum_{i=1}^S \norm{\hat{\mbf{p}}^{(k,i)}-\hat{\mbf{p}}^{(k-1,i)}}{\R^P \otimes \R^N_{\mbf{M}}}^2 \right)^{1/2} < tol.
\]
\end{enumerate}
Return $\{ \mc{P}_{\mc{A}}(\hat{\mbf{p}}^{(k,i)}) \}_{i=1}^S$ as a basis for the approximate eigenspace.
\end{algorithm}

\begin{remark}
The motivation behind step (2.0) in Algorithm \ref{alg:ssi} is the following: In the event of an eigenvalue crossing we expect the sum of the eigenvectors to be smooth on $\Gamma$ even if the individual eigenvectors as defined by \eqref{denumeration} are not.  
\end{remark}

\subsection{Convergence of the polynomial approximations}

The stochastic collocation and stochastic Galerkin strategies proposed in this section result in polynomial approximations of the eigenspace of interest. The accuracy of these approximations is ultimately determined by the smoothness of the solution as well as the choice of the multi-index set $\mc{A} \subset (\mbb{N}_0^{\infty})_c$.

It has been shown in \cite{andreevschwab12} that, under certain conditions, simple eigenpairs of a parameter dependent elliptic operator are in fact complex analytic functions of the input parameters. In this case there exists a sequence of multi-index sets $\mc{A}_{\epsilon}\subset (\mbb{N}_0^{\infty})_c$ so that the approximation error is bounded by $(\# \mc{A}_{\epsilon})^{-r}$ for some $r > 0$ as $\epsilon \to 0$. In other words the collocated eigenpairs and the fixed points of the spectral inverse iteration converge to the exact eigenpair at an algebraic rate with respect to the number of multi-indices. We refer to \cite{andreevschwab12, hakulalaaksonen17} for details. If we assume that the subspace associated to the eigenvalues $\{ \lambda^{(i)}_p(\xi) \}_{i=1}^S$ is complex analytic as a function of $\xi$, then we expect to see a similar rate of convergence for the basis vectors generated by Algorithms \ref{alg:msc} and \ref{alg:ssi}.
\section{Numerical Experiments}
\label{sec:numerical_experiments}

In this section we apply the different solution strategies to the shell eigenvalue problem with random material properties. The goal of our numerical experiments is to test the functionality of the proposed strategies, to illustrate the issue of eigenvalue crossings, and to analyze the asymptotic behaviour of the solutions as the shell thickness tends to zero.

The midsurface of the shell is $D = [-1,1] \times [0,2\pi]$. For the 1D solver we use a mesh with 16 elements. For the 2D solver we use a mesh with $16 \times 8$ quadrilateral elements in the axial and angular dimensions respectively. The Naghdi shell model is used in all examples except for the last one 
(featuring a general 2D material uncertainty), where we compare the Naghdi and mathematical shell models.

We set the Poisson ratio equal to $\nu = 1/3$ and assume the random material model \eqref{kl} with $E_0(\x) = \bar{E}$ and $E_m(\x) = \bar{E}(m+1)^{-2}\phi_m(\x)$. We set
\begin{equation}
\label{axial}
\phi_m(\x) = \sin(\pi m x) \quad m = 1,2, \ldots \quad \x = (x,y) \in D
\end{equation}
in the case of axial uncertainty and
\begin{equation}
\label{general}
\phi_m(\x) = \left\{ \begin{array}{l} \sin(\pi m x) \quad m = 1,3, \ldots \\ \sin(m y) \quad m = 2,4, \ldots \end{array} \right. \quad \x = (x,y) \in D
\end{equation}
in the case of general uncertainty. Without loss of generality, we assume a scaling of the results so that the mean value of the Young's modulus is $\bar{E} = 1$.

In each example we are interested in the smallest eigenvalue and the corresponding eigenfunction or alternatively the subspace associated to the cluster of few smallest eigenvalues. In the 1D solver we choose the wavenumber of the angular component accordingly.

As in \cite{andreevschwab12, hakulakaarniojalaaksonen15, hakulalaaksonen17, bierischwab09,bieriandreevschwab09} we use monotone multi-index sets of the type
\begin{equation}
\label{miset}
\mc{A}_{\epsilon} := \left\{ \a \in (\N_0^{\infty})_c \ \left\vert \ \ \prod_{m \in \supp \a} \eta_m^{\a_m} > \epsilon  \right\} \right.,
\end{equation}
which can be generated using the algorithm given in \cite{bieriandreevschwab09}. Here $\{ \eta_m \}_{m=1}^{\infty} \subset \R$ is a suitably chosen decreasing sequence such that $0 < \eta_m < 1$ for all $m \ge 1$.

In the following examples we have focused on the fourth component $\theta$ of the eigenmode, since it is usually considered the most challenging one to approximate. Similar results would in general be observed for the other components.

\subsection{Calibration}
\label{sec:numex_calibration}

\begin{figure}[htb]
\begin{center}
\subfloat[{The values $| \E[\lambda_p] - \E[\lambda_*] |$.}]{\includegraphics[width=0.48\textwidth]{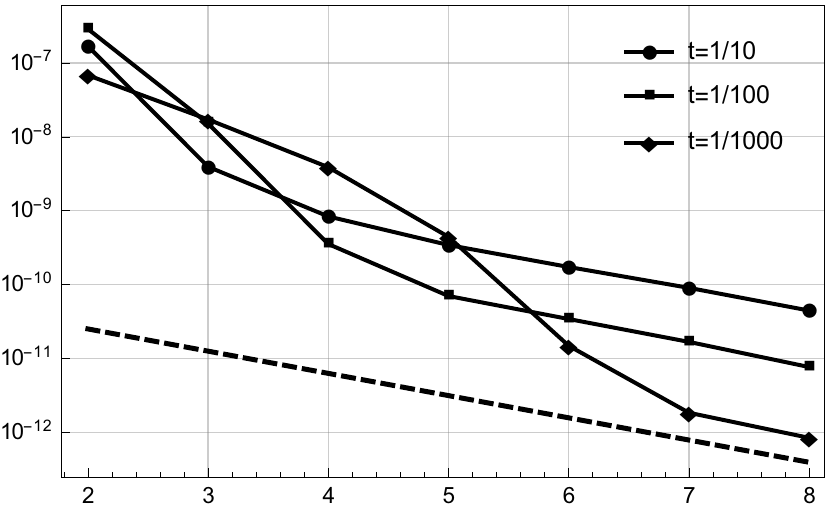}}\quad
\subfloat[{The values $\norm{\E[\theta_p] - \E[\theta_*]}{L^2([-1,1])}$.}]{\includegraphics[width=0.48\textwidth]{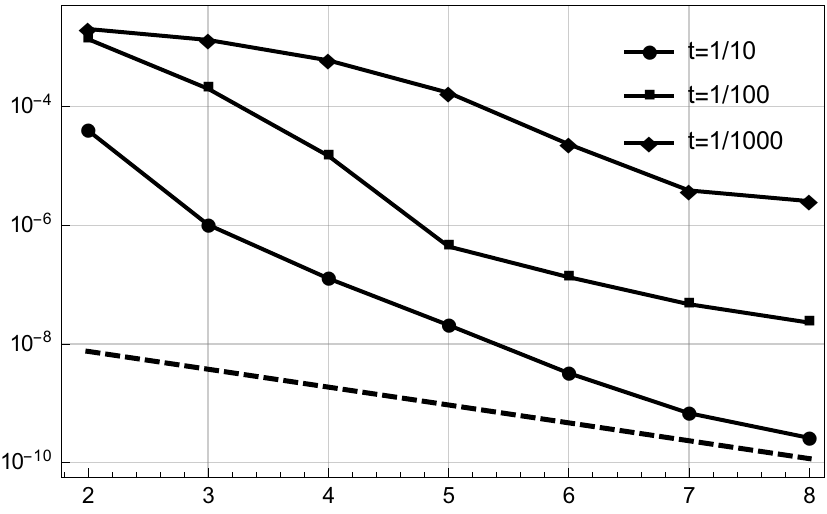}}\quad
\subfloat[{The values $| \Var[\lambda_p] - \Var[\lambda_*] |$.}]{\includegraphics[width=0.48\textwidth]{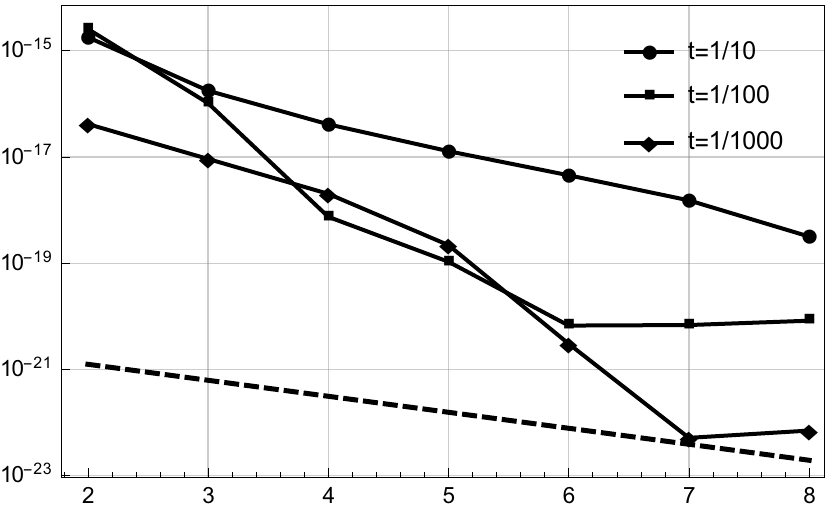}}\quad
\subfloat[{The values $\norm{\Var[\theta_p] - \Var[\theta_*]}{L^2([-1,1])}$.}]{\includegraphics[width=0.48\textwidth]{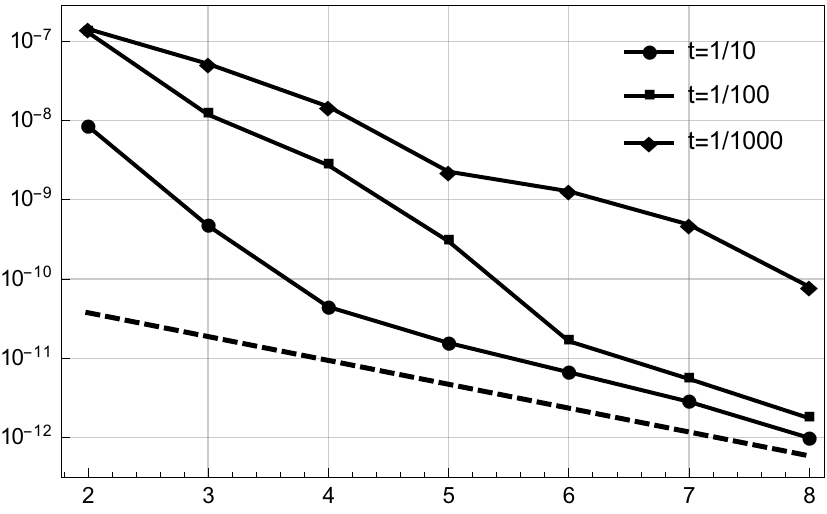}}
\caption{Spatial convergence of the solution $(\lambda_p, \theta_p)$ (the eigenvalue and the fourth axial field component) as computed by Algorithm \ref{alg:sii} to the overkill solution $(\lambda_*, \theta_*)$. A log plot of the solution statistics as a function of the polynomial order $p$. The dashed lines represent the rate $1/2^p$.}
\label{fig:spatialconv}
\end{center}
\end{figure}

We first calibrate the spatial discretization using the 1D Galerkin solver. The model of axial uncertainty \eqref{axial} is assumed. The spatial convergence of the solution as computed by Algorithm \ref{alg:sii} has been presented in Figure \ref{fig:spatialconv}. We have used the values $t = 1/10$, $t = 1/100$ and $t = 1/1000$ for the shell thickness. In each case we have
\[
\bar{\lambda}_{1/2} := \frac{\lambda_p^{(1)}(0)}{\lambda_p^{(2)}(0)} < 0.32
\]
for the ratio of the two smallest discrete eigenvalues, and therefore we expect \eqref{separated} to hold with $S = 1$. We set $\epsilon = 10^{-4}$ in \eqref{miset} so that $M_{\mc{A}_{\epsilon}} = 99$ and $\# \mc{A}_{\epsilon} = 358$. As a reference we have used an overkill solution computed with $p = 10$.

The convergence graphs indicate numerical locking due to (unknown)
error amplification factor $C(t)$ which increases as $t \to 0$.
However, even for $t=1/1000$ with sufficiently high $p$ we get convergence.

\subsection{Validation}
\label{sec:numex_validation}

\begin{figure}[p!]
\begin{center}
\subfloat[{$\E[u]$.}]{\includegraphics[width=0.17\textwidth]{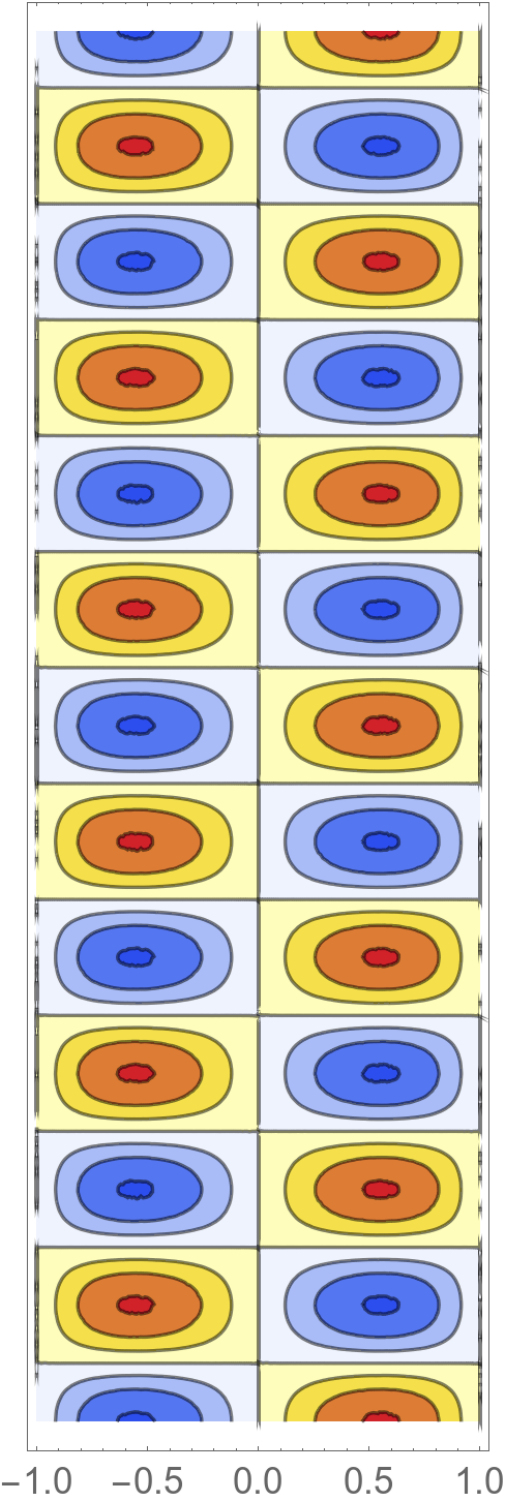}}\quad
\subfloat[{$\E[v]$.}]{\includegraphics[width=0.17\textwidth]{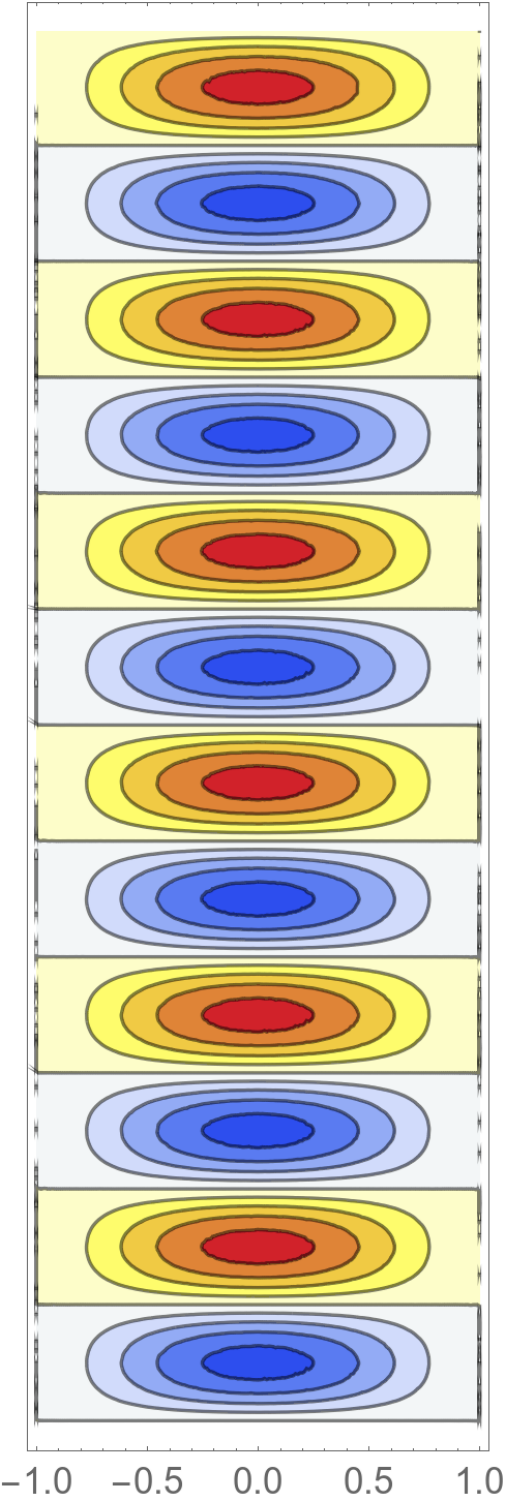}}\quad
\subfloat[{$\E[w]$.}]{\includegraphics[width=0.17\textwidth]{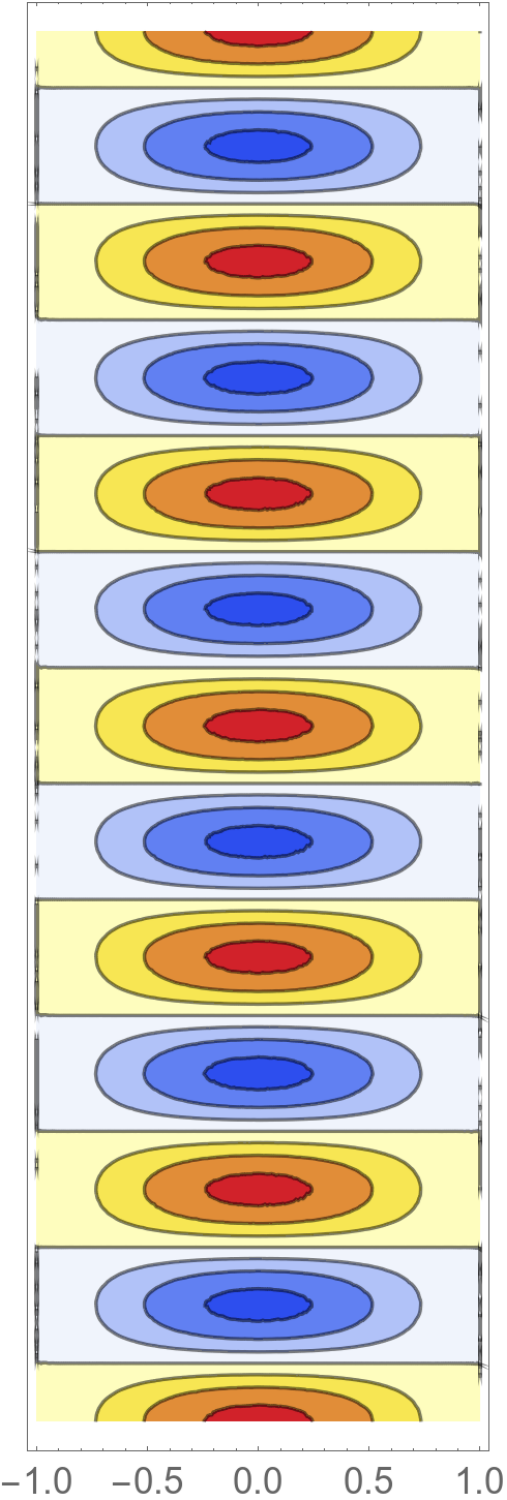}}\quad
\subfloat[{$\E[\theta]$.}]{\includegraphics[width=0.17\textwidth]{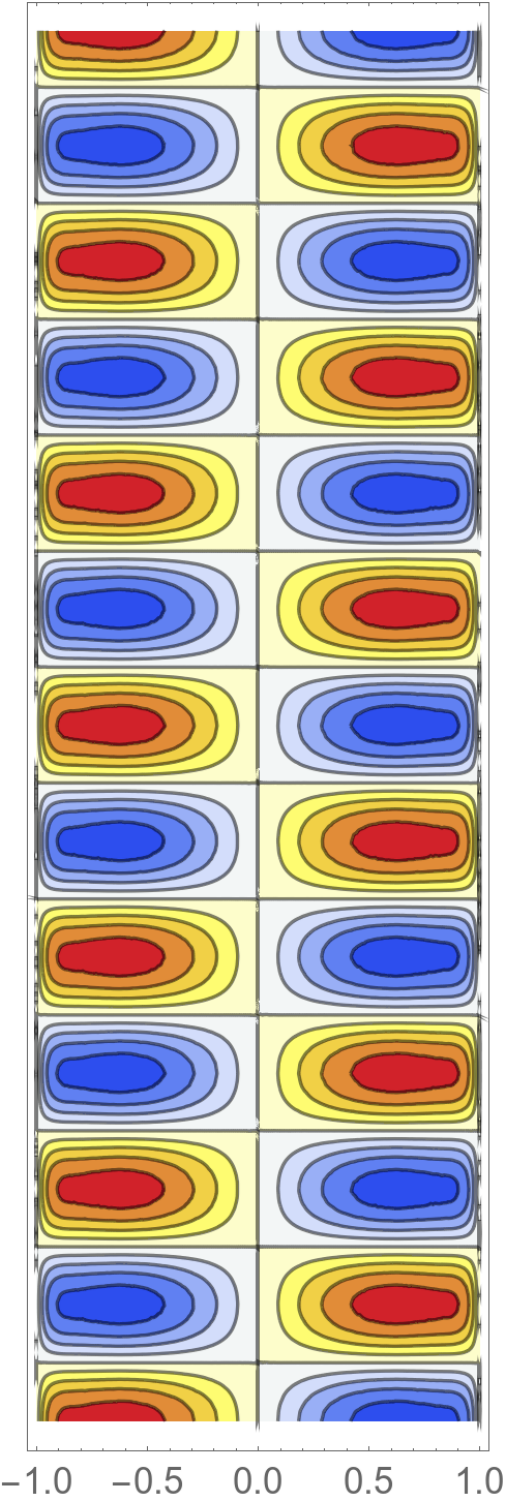}}\quad
\subfloat[{$\E[\psi]$.}]{\includegraphics[width=0.17\textwidth]{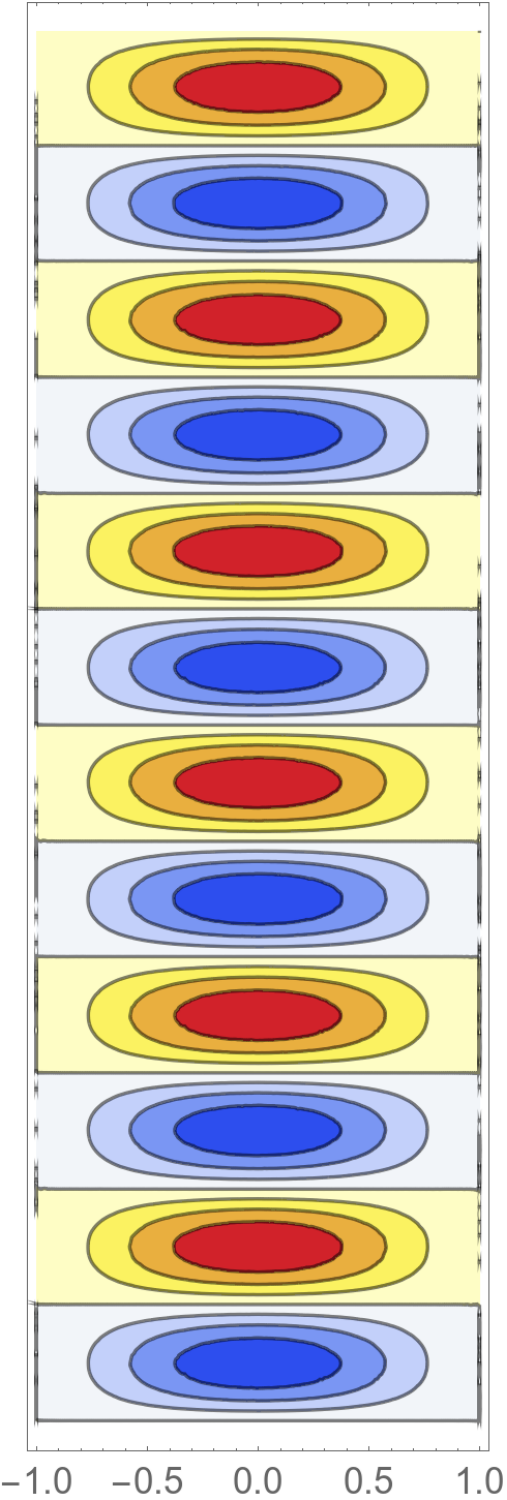}}\quad
\subfloat[{$\Var[u]$.}]{\includegraphics[width=0.17\textwidth]{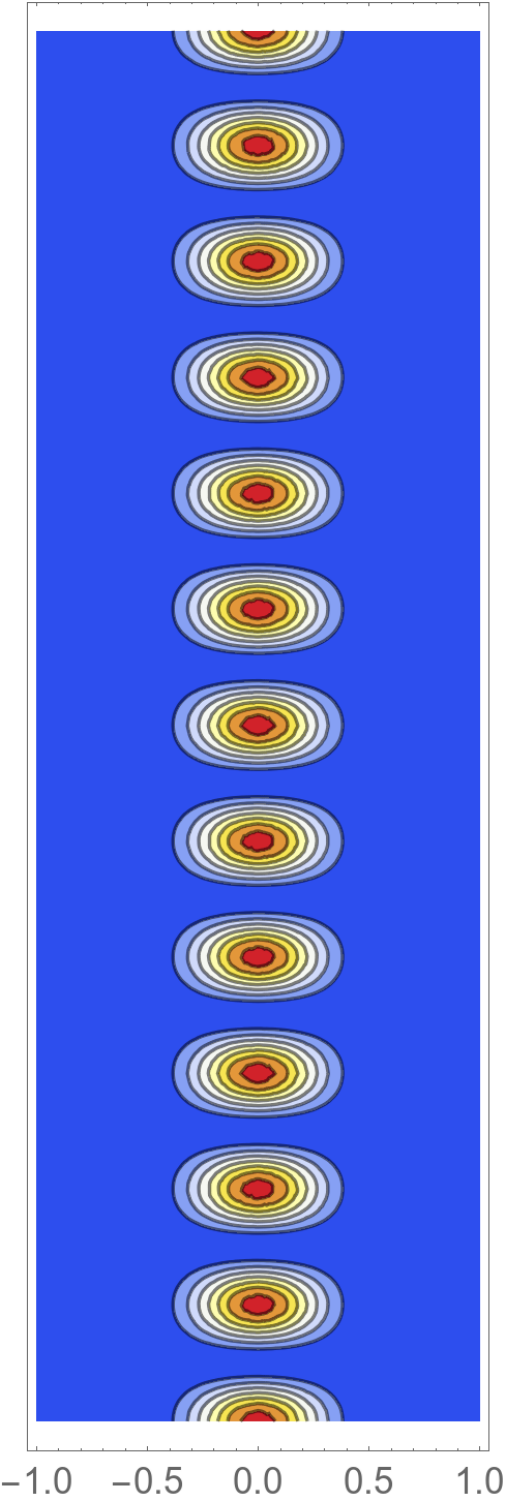}}\quad
\subfloat[{$\Var[v]$.}]{\includegraphics[width=0.17\textwidth]{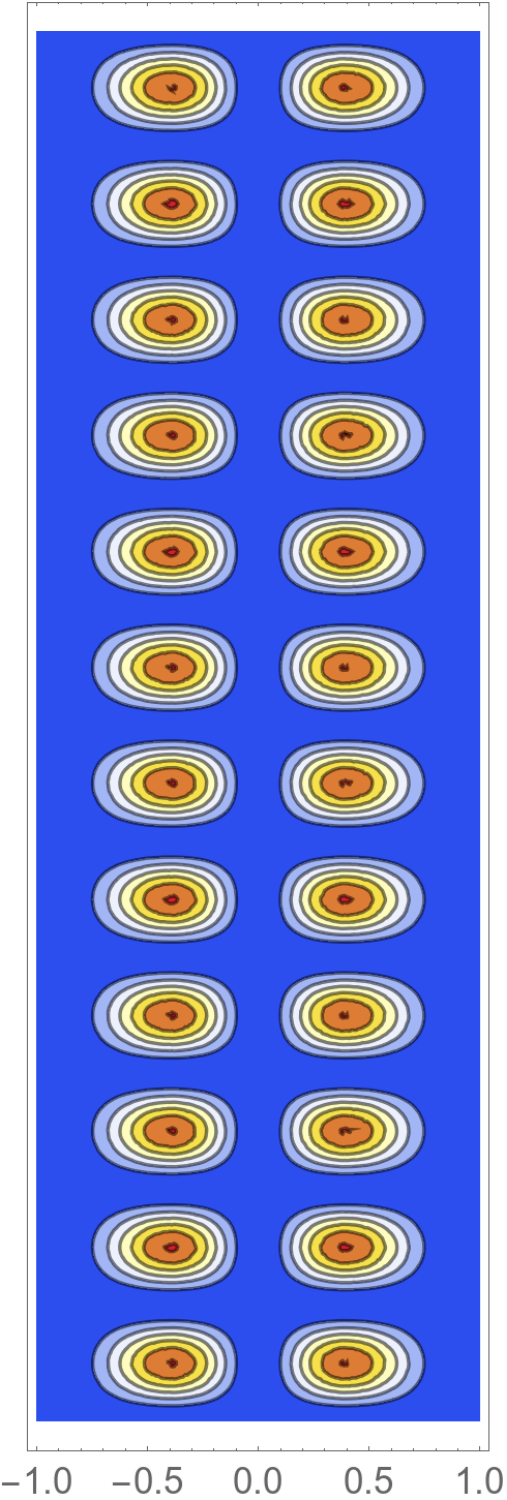}}\quad
\subfloat[{$\Var[w]$.}]{\includegraphics[width=0.17\textwidth]{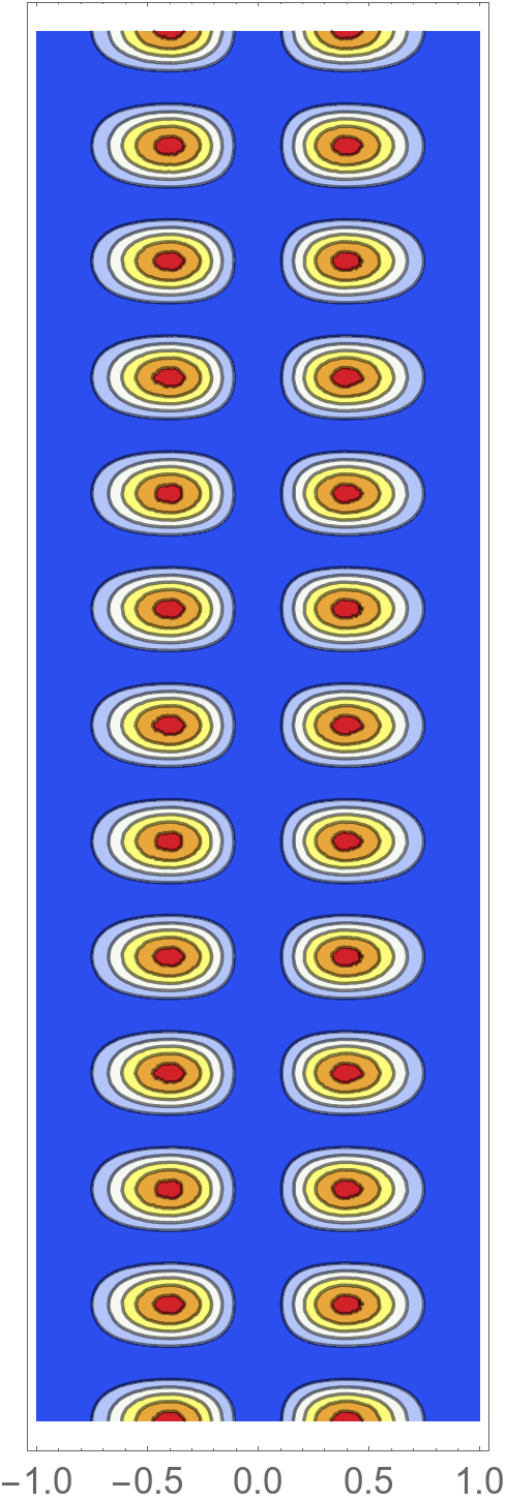}}\quad
\subfloat[{$\Var[\theta]$.}]{\includegraphics[width=0.17\textwidth]{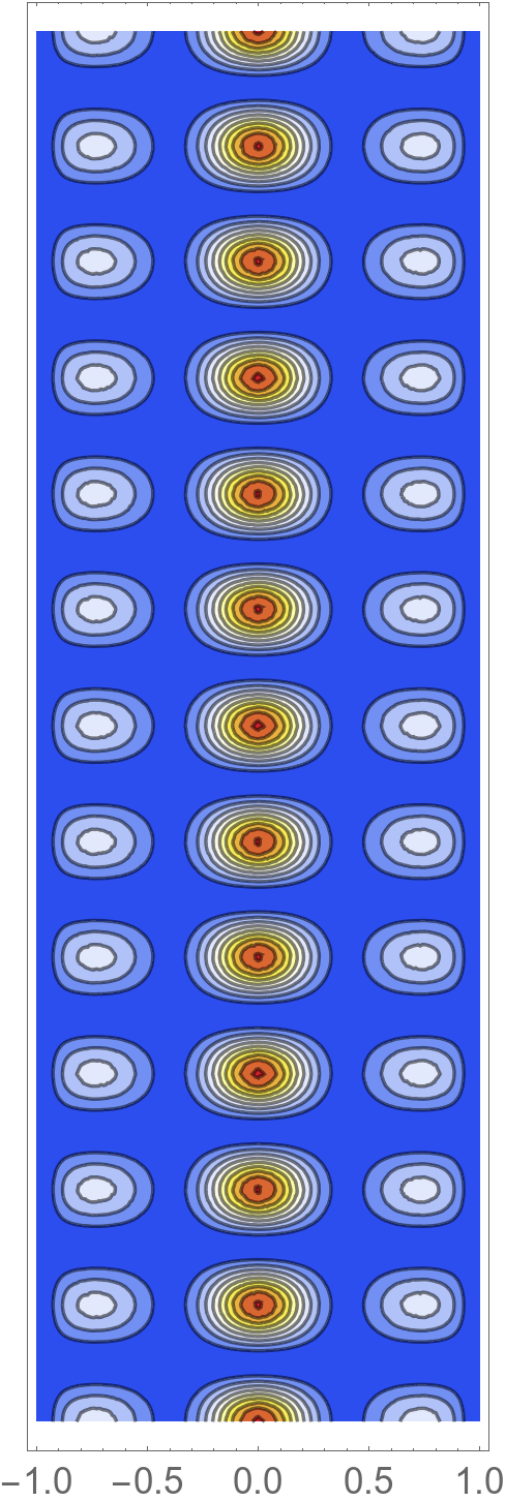}}\quad
\subfloat[{$\Var[\psi]$.}]{\includegraphics[width=0.17\textwidth]{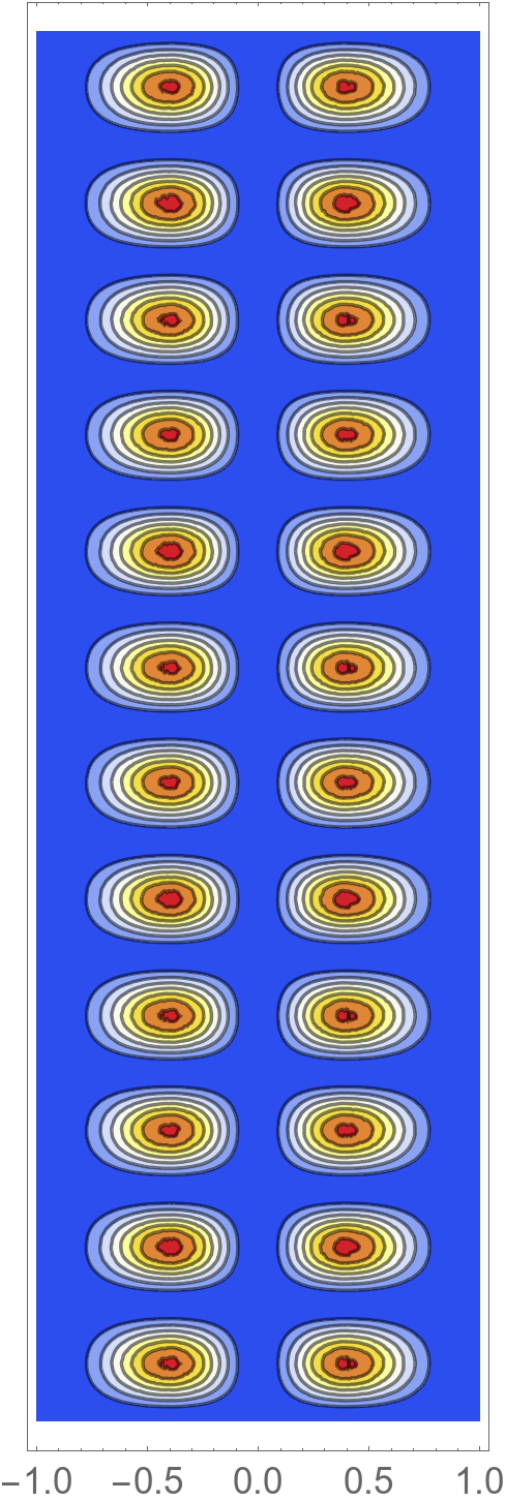}}
\caption{The solution for $t=1/100$ in the case of axial uncertainty. Expected value and variance of the first basis function of the eigenspace. A contour plot of the five field components.}
\label{fig:comps}
\end{center}
\end{figure}

We next validate both the 1D solver and the 2D solver. To this end we assume the model of axial uncertainty \eqref{axial} and set $t = 1/100$. We set $p = 8$ for the 1D solver and $p=6$ for the 2D solver. Statistics of the first basis function of the two dimensional eigenspace have been presented in Figure \ref{fig:comps}. The wavenumber of the axial component is $k = 6$.

\begin{remark}
As noted, the eigenpairs defined by the enumeration \eqref{denumeration} are not necessarily smooth and therefore computation of each eigenmode individually is not justified. However, in the case of axial uncertainty each eigenvalue is actually a double eigenvalue (see \ref{sec:1D_models}). Thus, when computing the associated two-dimensional subspace we may in fact evaluate the two eigenvalues explicitly since they admit the same value for every $\xi \in \Gamma$.
\end{remark}

\begin{figure}[htb]
\begin{center}
\subfloat[{The values $| \E[\lambda_{\mc{A}_{\epsilon}}] - \E[\lambda_*] |$.}]{\includegraphics[width=0.48\textwidth]{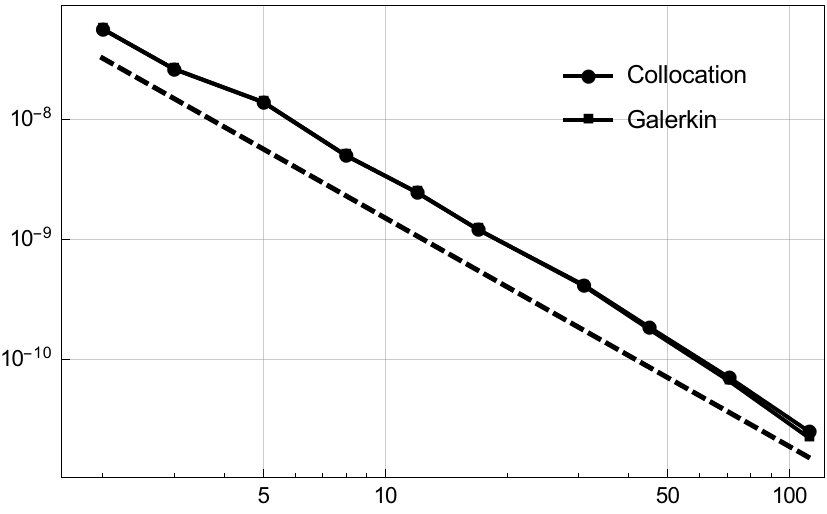}}\quad
\subfloat[{The values $\norm{\E[\theta_{\mc{A}_{\epsilon}}] - \E[\theta_*]}{L^2([-1,1])}$.}]{\includegraphics[width=0.48\textwidth]{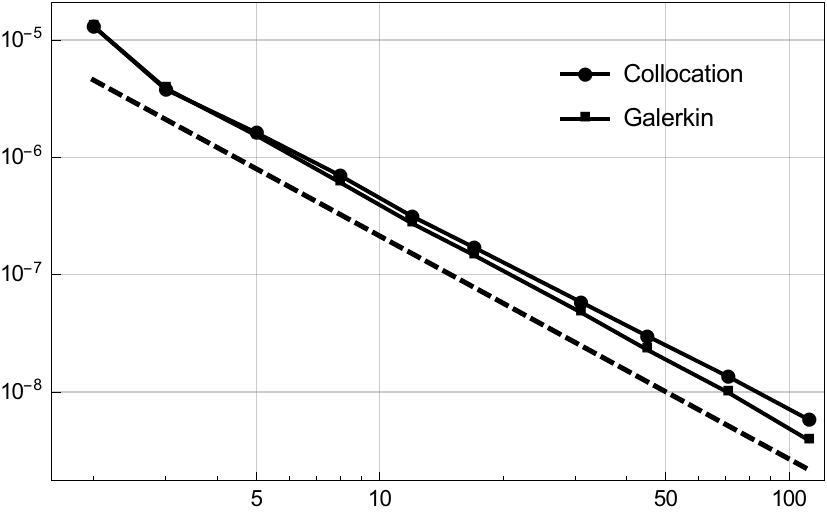}}\quad
\subfloat[{The values $| \Var[\lambda_{\mc{A}_{\epsilon}}] - \Var[\lambda_*] |$.}]{\includegraphics[width=0.48\textwidth]{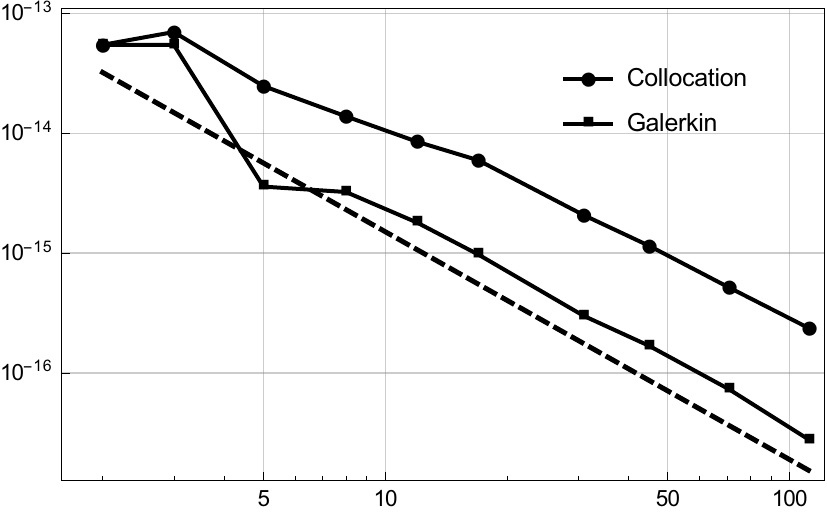}}\quad
\subfloat[{The values $\norm{\Var[\theta_{\mc{A}_{\epsilon}}] - \Var[\theta_*]}{L^2([-1,1])}$}]{\includegraphics[width=0.48\textwidth]{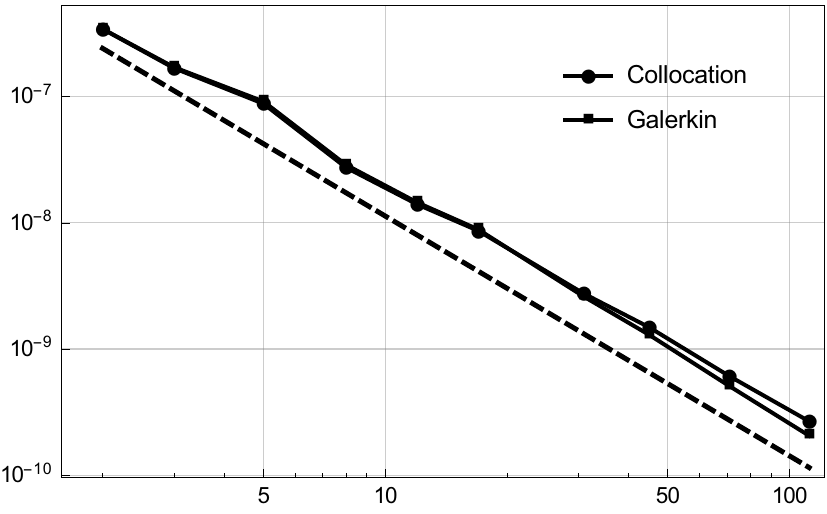}}
\caption{Stochastic convergence of the 1D solution $(\lambda_{\mc{A}_{\epsilon}}, \theta_{\mc{A}_{\epsilon}})$ (the eigenvalue and the fourth axial field component) as computed by Algorithms \ref{alg:ssc} (collocation) and \ref{alg:sii} (Galerkin) to the overkill solution $(\lambda_*, \theta_*)$. A log-log plot of the solution statistics as a function of the basis size $\# \mc{A}_{\epsilon}$. The dashed lines represent the rate $(\# \mc{A}_{\epsilon})^{-1.9}$.}
\label{fig:stochasticconv1d}
\end{center}
\end{figure}

In Figure \ref{fig:stochasticconv1d} we have illustrated the convergence of the 1D solution towards an overkill solution as a function of $\# \mc{A}_{\epsilon}$. For the overkill solution we have set $\epsilon = 10^{-4}$ in \eqref{miset} so that $M_{\mc{A}_{\epsilon}} = 99$ and $\# \mc{A}_{\epsilon} = 358$. As before we have $\bar{\lambda}_{1/2} < 0.32$. We see that the solutions computed using Algorithm \ref{alg:sii} almost perfectly agree with the solutions computed using Algorithm \ref{alg:ssc}.

\begin{figure}[htb]
\begin{center}
\subfloat[{The values $| \E[\lambda_{\mc{A}_{\epsilon}}] - \E[\lambda_*] |$.}]{\includegraphics[width=0.48\textwidth]{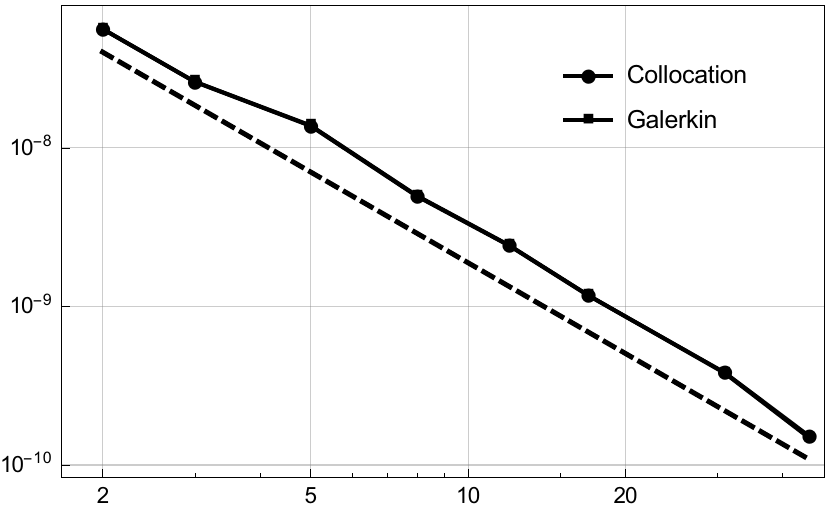}}\quad
\subfloat[{The values $\norm{\E[\theta_{\mc{A}_{\epsilon}}] - \E[\theta_*]}{L^2(D)}$.}]{\includegraphics[width=0.48\textwidth]{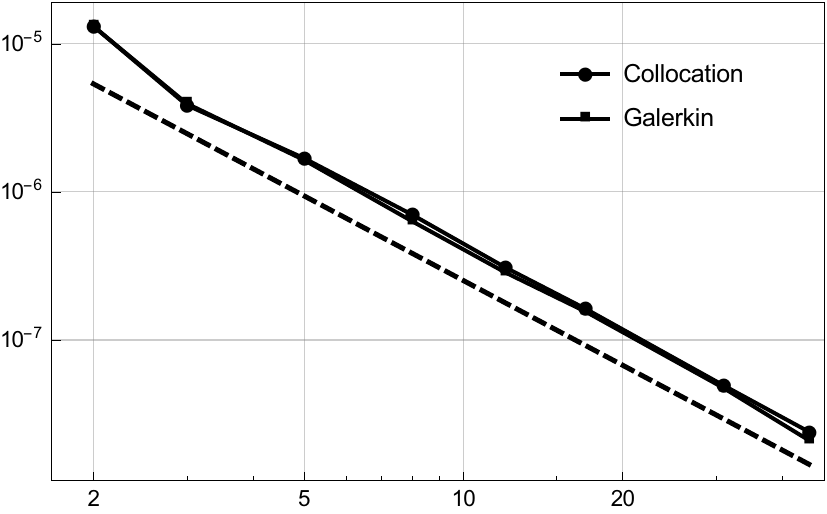}}\quad
\subfloat[{The values $| \Var[\lambda_{\mc{A}_{\epsilon}}] - \Var[\lambda_*] |$.}]{\includegraphics[width=0.48\textwidth]{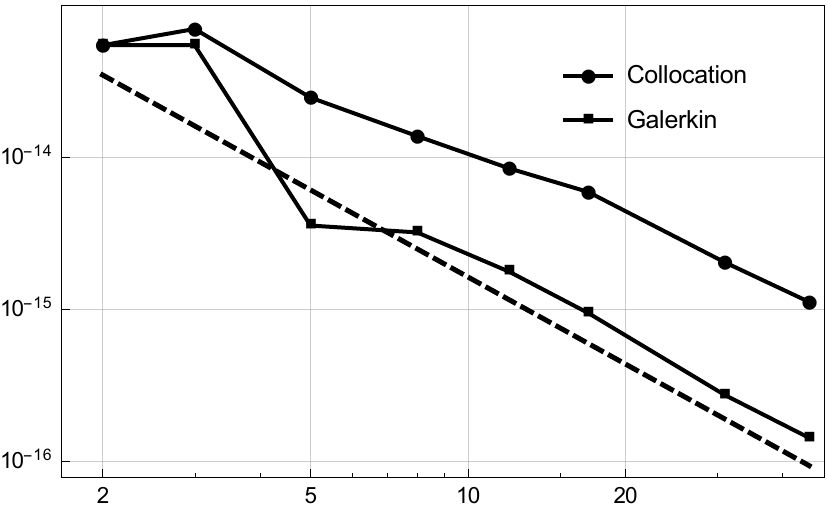}}\quad
\subfloat[{The values $\norm{\Var[\theta_{\mc{A}_{\epsilon}}] - \Var[\theta_*]}{L^2(D)}$.}]{\includegraphics[width=0.48\textwidth]{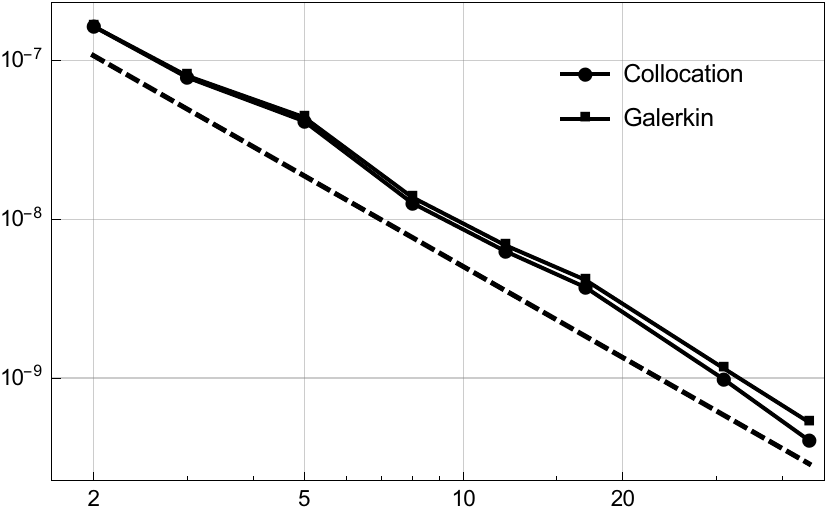}}
\caption{Stochastic convergence of the 2D solution $(\lambda_{\mc{A}_{\epsilon}}, \theta_{\mc{A}_{\epsilon}})$ (the eigenvalue and the fourth field component) as computed by Algorithms \ref{alg:msc} (collocation) and \ref{alg:ssi} (Galerkin) to the overkill solution $(\lambda_*, \theta_*)$. A log-log plot of the solution statistics as a function of the basis size $\# \mc{A}_{\epsilon}$. The dashed lines represent the rate $(\# \mc{A}_{\epsilon})^{-1.9}$.}
\label{fig:stochasticconv2d}
\end{center}
\end{figure}

In Figure \ref{fig:stochasticconv2d} we have illustrated the convergence of the 2D solution towards an overkill solution as a function of $\# \mc{A}_{\epsilon}$. For the overkill solution we have set $\epsilon = 5 \cdot 10^{-4}$ in \eqref{miset} so that $M_{\mc{A}_{\epsilon}} = 43$ and $\# \mc{A}_{\epsilon} = 116$. In the 2D model we have
\[
\bar{\lambda}_{2/3} := \frac{\lambda_p^{(2)}(0)}{\lambda_p^{(3)}(0)} < 0.93
\]
and the results have been computed using Algorithms \ref{alg:msc} and \ref{alg:ssi} with $S=2$. We have shown results for only one of the two basis vectors but practically identical results would be observed for the other basis vector as well. We see that the solutions computed using Algorithm \ref{alg:ssi} are again in excellent agreement with the solutions computed using Algorithm \ref{alg:msc}. Moreover, the convergence rates agree with the ones obtained for the 1D model.

\begin{figure}[htb]
\begin{center}
\subfloat[{The values $\delta^{(k)}$ in Algorithm \ref{alg:sii}.}]{\includegraphics[width=0.48\textwidth]{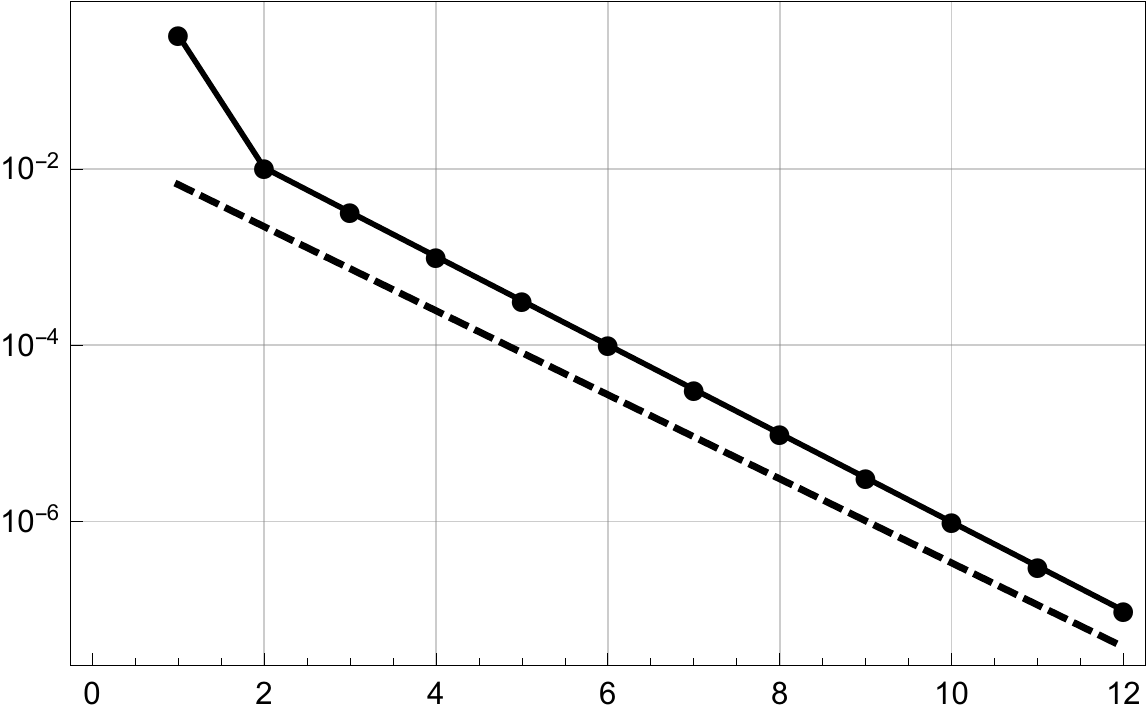}}\quad
\subfloat[{The values $\delta^{(k)}$ in Algorithm \ref{alg:ssi}.}]{\includegraphics[width=0.48\textwidth]{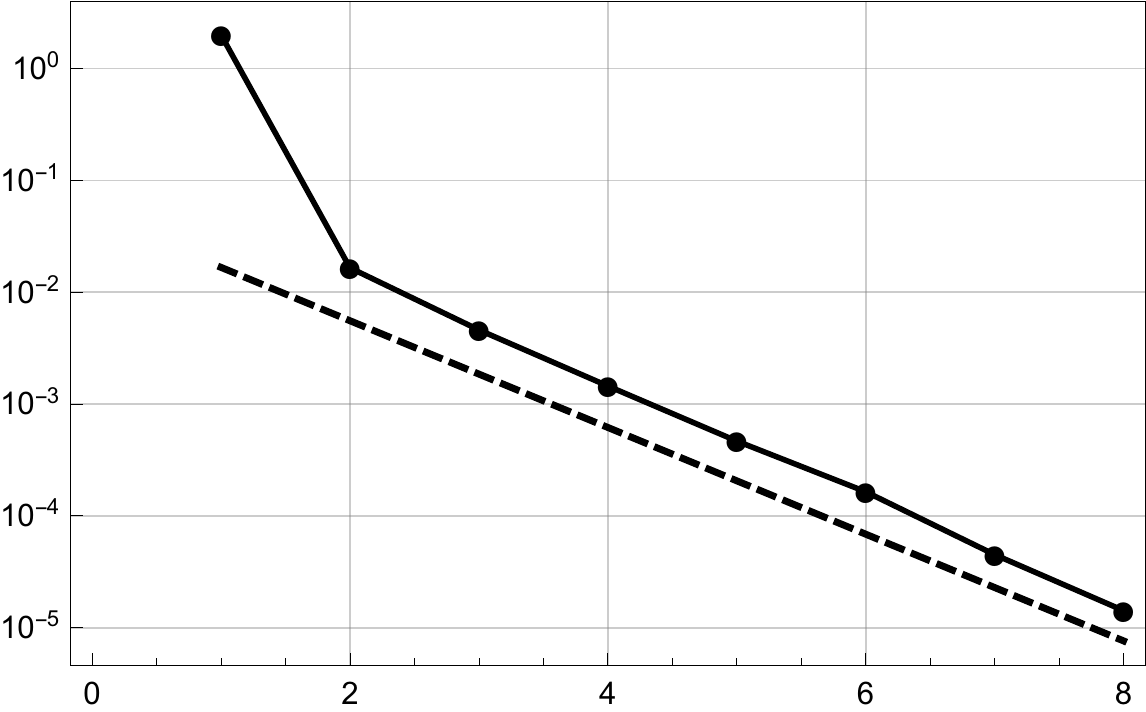}}\quad
\caption{Convergence of the spectral inverse iteration and spectral subspace iteration algorithms. A log plot of the values $\delta^{(k)}$ as a function of $k$. The dashed lines represent the rate $1/3^k$.}
\label{fig:iterationconv}
\end{center}
\end{figure}

The convergence of the inverse iteration (Algorithm \ref{alg:sii}) and the spectral subspace iteration (Algorithm \ref{alg:ssi}) for $p = 6$ and $\epsilon = 5 \cdot 10^{-4}$ so that $\# \mc{A}_{\epsilon} = 116$ has been presented in Figure \ref{fig:iterationconv}. The 1D and 2D systems have approximately $N = 4.9 \cdot 10^2$ and $N = 2.4 \cdot 10^4$ degrees of freedom in space and $5.6 \cdot 10^4$ and $2.8 \cdot 10^6$ degrees of freedom in total respectively. The computation times are around 15 seconds in the 1D case and around 50 minutes in the 2D one on a single core 3.40GHz CPU (Intel Xeon E3-1230 v5) with 32 GiB memory. The respective computation times for the collocation method are an order of magnitude longer with sequential computations.

\subsection{Eigenvalue crossings}
\label{sec:numex_subspace}

\begin{figure}[htb]
\begin{center}
\subfloat[{Eigenvalues at $\xi_3 = \xi_4 = \ldots = 1/2$ as a function of $(\xi_1, \xi_2) \in [-1,1]^2$.}]{\includegraphics[width=0.48\textwidth]{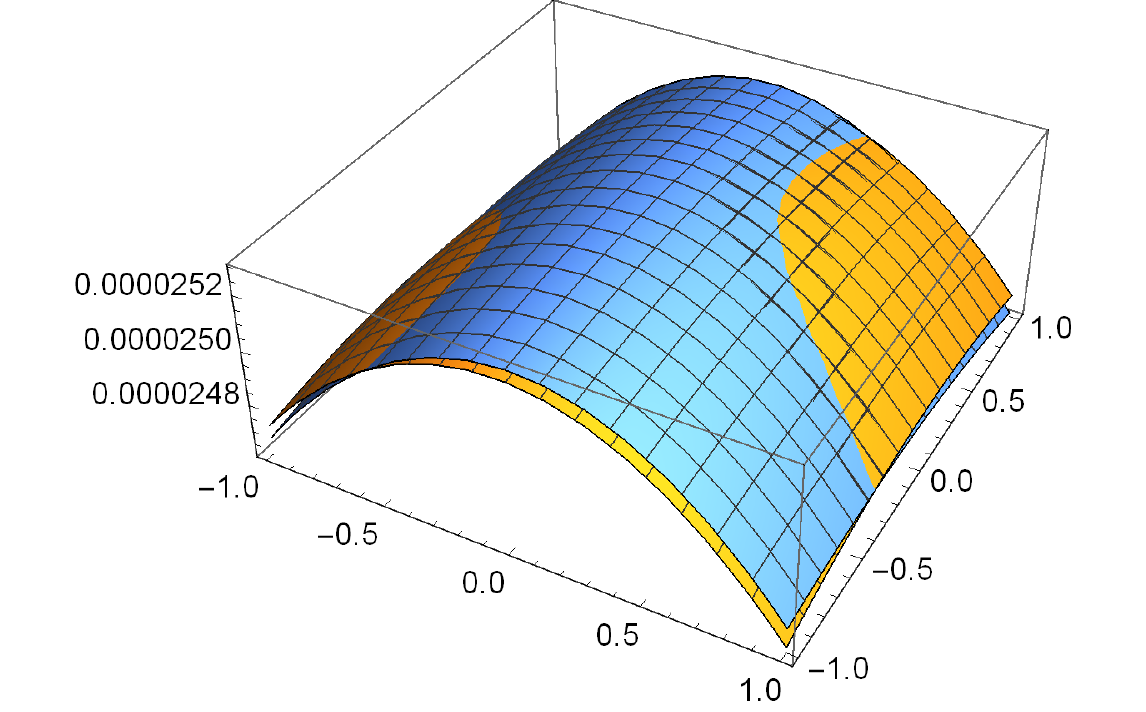}}\quad
\subfloat[{Eigenvalues at $\xi_1 = \xi_3 = \xi_4 = \ldots = 1/2$ as a function of $\xi_2 \in [-1,1]$.}]{\includegraphics[width=0.48\textwidth]{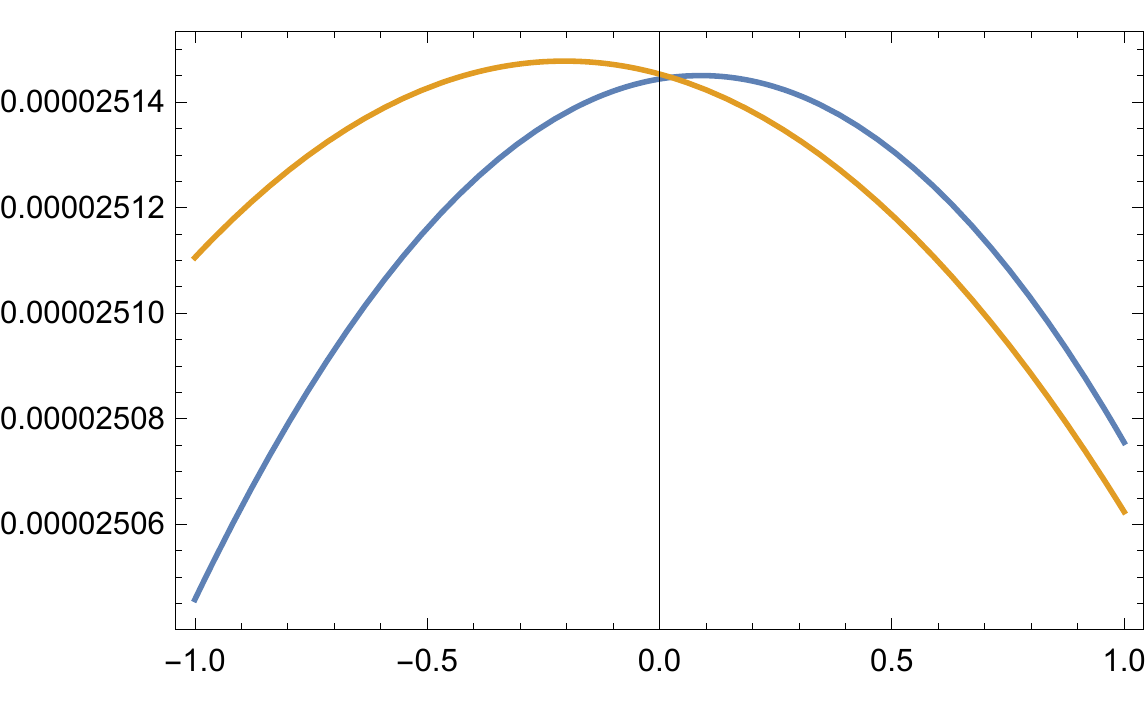}}\quad
\caption{The smallest eigenvalues of the 1D model and axial uncertainty with wavenumbers $k=6$ and $k=7$. Solution computed by Algorithm \ref{alg:sii} for $t = 0.0067$.}
\label{fig:crossing}
\end{center}
\end{figure}

The aim of our next experiment is to illustrate that for a suitably chosen value of $t$ we observe an eigenvalue crossing and as a result the effective smallest eigenmode is in fact discontinuous on $\Gamma$. Consider the axial uncertainty model \eqref{axial} and a value $t = 0.0067$ for the shell thickness. We set $p = 8$ and $\epsilon = 10^{-4}$ in \eqref{miset} so that $M_{\mc{A}_{\epsilon}} = 99$ and $\# \mc{A}_{\epsilon} = 358$. We compute the smallest eigenvalue of the problem by applying Algorithm \ref{alg:sii} on the 1D model with wavenumbers $k=6$ and $k=7$. The results are shown in Figure \ref{fig:crossing} from which a crossing of the eigenvalues is evident. Moreover, since the two modes have different wavenumbers in the angular component, the effective smallest eigenmode must be discontinuous in the proximity of the crossing.

\begin{figure}[htb]
\begin{center}
\subfloat[{The values $\norm{\E[\theta_{\mc{A}_{\epsilon}}] - \E[\theta_*]}{L^2(D)}$.}]{\includegraphics[width=0.48\textwidth]{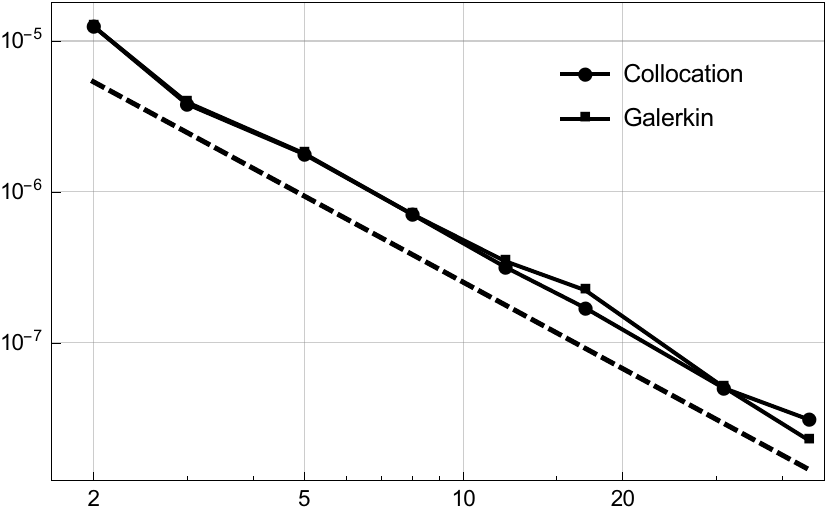}}\quad
\subfloat[{The values $\norm{\Var[\theta_{\mc{A}_{\epsilon}}] - \Var[\theta_*]}{L^2(D)}$.}]{\includegraphics[width=0.48\textwidth]{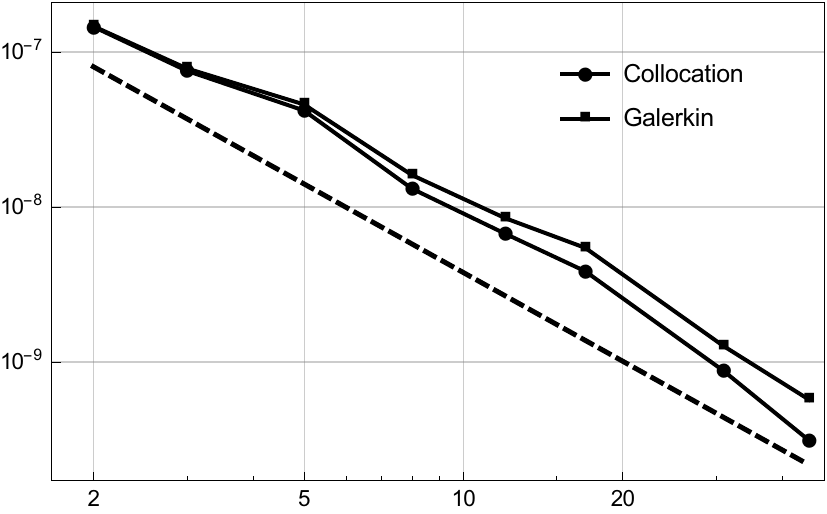}}
\caption{Stochastic convergence of the 2D solution $\theta_{\mc{A}_{\epsilon}}$ (the fourth field component) as computed by Algorithms \ref{alg:msc} (collocation) and \ref{alg:ssi} (Galerkin) to the overkill solution $\theta_*$ in the case of an eigenvalue crossing. A log-log plot of the solution statistics as a function of the basis size $\# \mc{A}_{\epsilon}$. The dashed lines represent the rate $(\# \mc{A}_{\epsilon})^{-1.9}$.}
\label{fig:crossingconv}
\end{center}
\end{figure}

In Figure \ref{fig:crossingconv} we have illustrated the convergence of the 2D solution in the case of the eigenvalue crossing ($t = 0.0067$) as a function of $\# \mc{A}_{\epsilon}$. As a reference we have again used an overkill solution for which $\epsilon = 5 \cdot 10^{-4}$ in \eqref{miset} so that $M_{\mc{A}_{\epsilon}} = 43$ and $\# \mc{A}_{\epsilon} = 116$. We now have
\[
\bar{\lambda}_{4/5} := \frac{\lambda_p^{(4)}(0)}{\lambda_p^{(5)}(0)} < 0.81
\]
and the results have been computed using Algorithms \ref{alg:msc} and \ref{alg:ssi} with $S=4$. We have shown results for only one of the four basis vectors. Similar results would be observed for the other basis vectors as well. We see that the solutions computed using Algorithm \ref{alg:ssi} are once again in excellent agreement with the solutions computed using Algorithm \ref{alg:msc}.

\subsection{Asymptotic analysis in the stochastic setting}
\label{sec:numex_asymptotic}

\begin{figure}[htb]
\begin{center}
\subfloat[{The eigenvalue.}]{\includegraphics[width=0.48\textwidth]{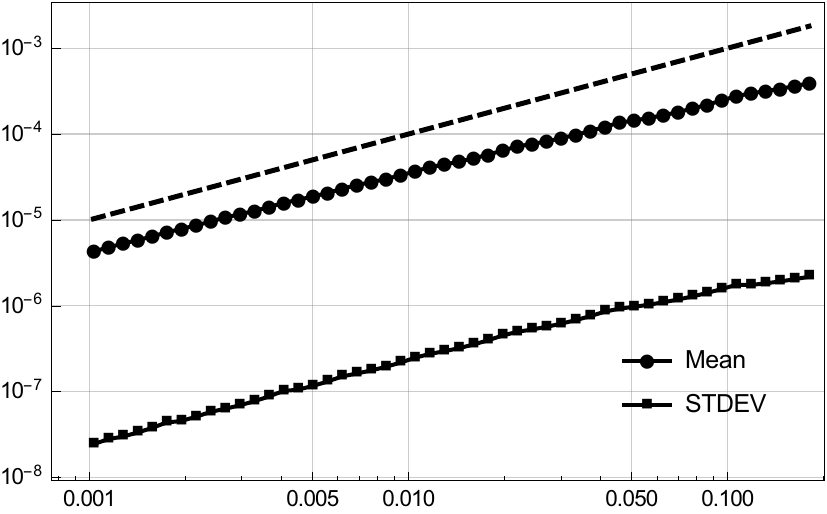}}\quad
\subfloat[{The deterministic wavenumber.}]{\includegraphics[width=0.48\textwidth]{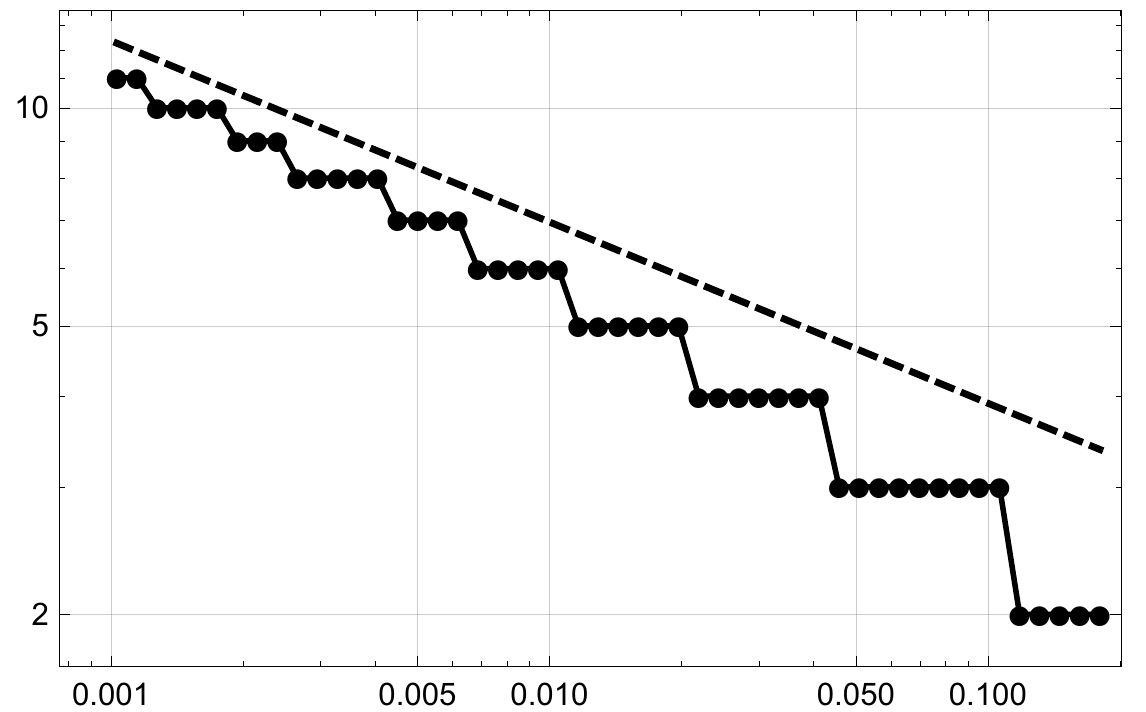}}
\caption{Asymptotics of the smallest eigenvalue and the associated deterministic wavenumber. A log-log plot of the solution statistics as a function of $t$. The dashed lines represent the rates $t$ and $t^{-1/4}$ respectively.}
\label{fig:tconv}
\end{center}
\end{figure}

We examine the asymptotics of the solution as $t \to 0$. To this end we assume the axial uncertainty model \eqref{axial} and employ the 1D Galerkin solver. Again we set $p = 8$ and $\epsilon = 10^{-4}$ in \eqref{miset} so that $M_{\mc{A}_{\epsilon}} = 99$ and $\# \mc{A}_{\epsilon} = 358$. In each case we select the wavenumber that gives the smallest eigenvalue at $\xi = 0$ and compute the corresponding eigenmode using Algorithm \ref{alg:sii}. In Figure \ref{fig:tconv} we have illustrated the behaviour of the eigenvalue and the associated wavenumber as a function of $t$. We see that the mean and variance of the eigenvalue converge to zero whereas the wavenumber grows as predicted by Theorem \ref{asymptotics}.

\subsection{General model of uncertainty}
\label{sec:numex_revolution}

\begin{figure}[p!]
\begin{center}
\subfloat[{$\E[u]$.}]{\includegraphics[width=0.17\textwidth]{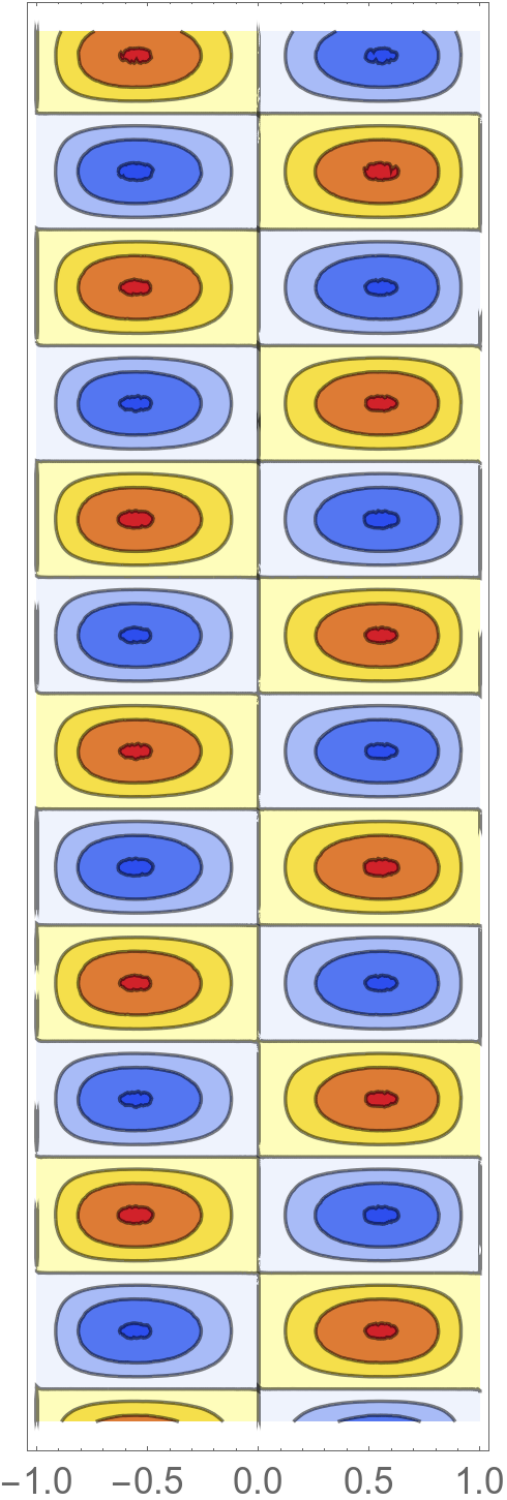}}\quad
\subfloat[{$\E[v]$.}]{\includegraphics[width=0.17\textwidth]{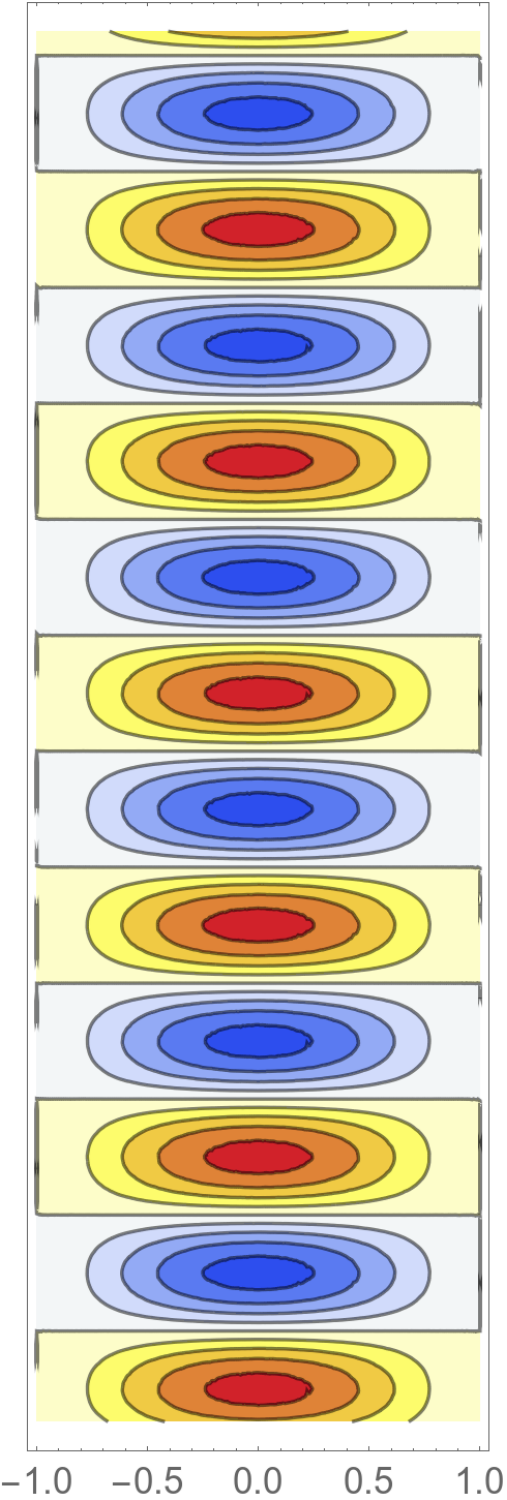}}\quad
\subfloat[{$\E[w]$.}]{\includegraphics[width=0.17\textwidth]{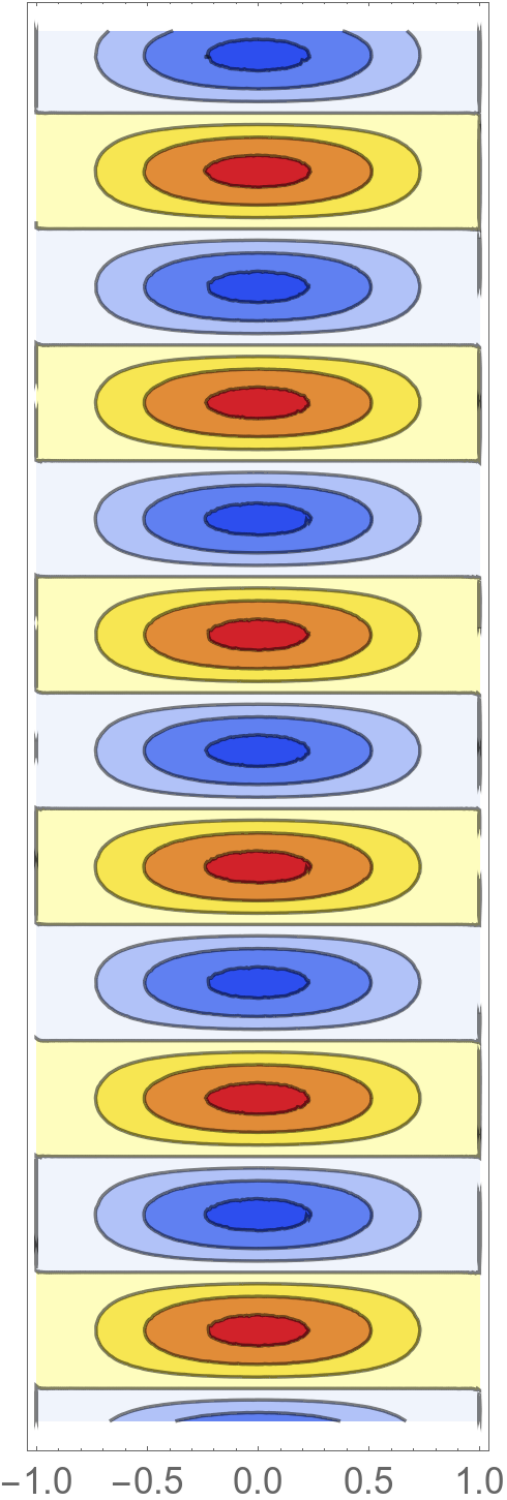}}\quad
\subfloat[{$\E[\theta]$.}]{\includegraphics[width=0.17\textwidth]{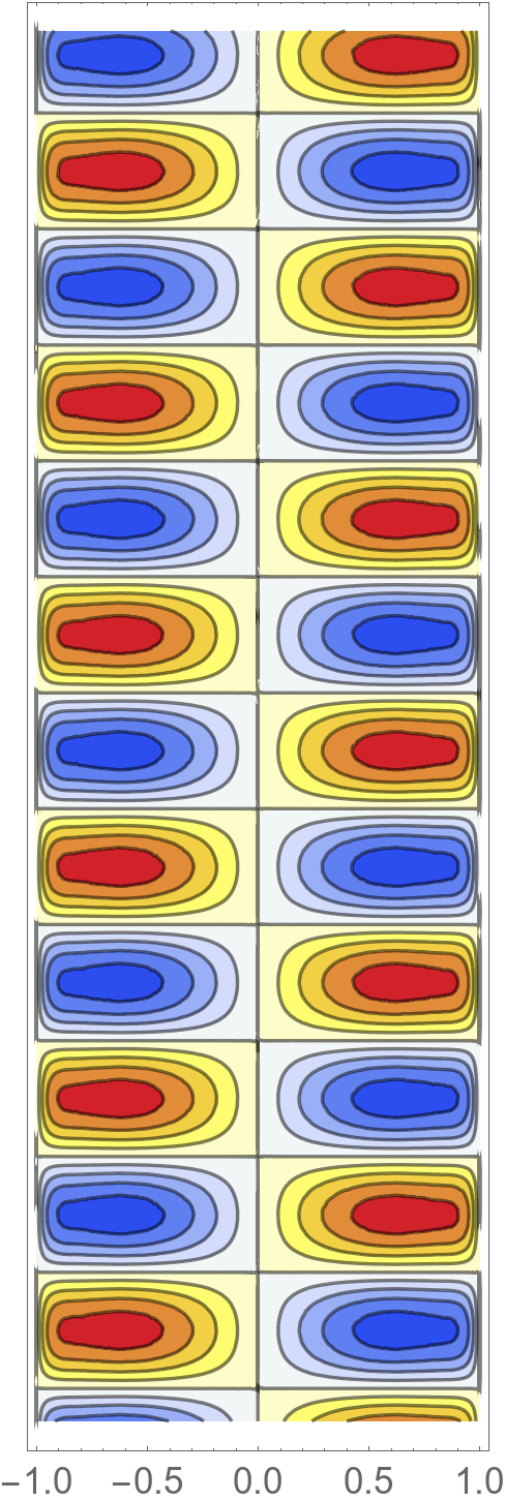}}\quad
\subfloat[{$\E[\psi]$.}]{\includegraphics[width=0.17\textwidth]{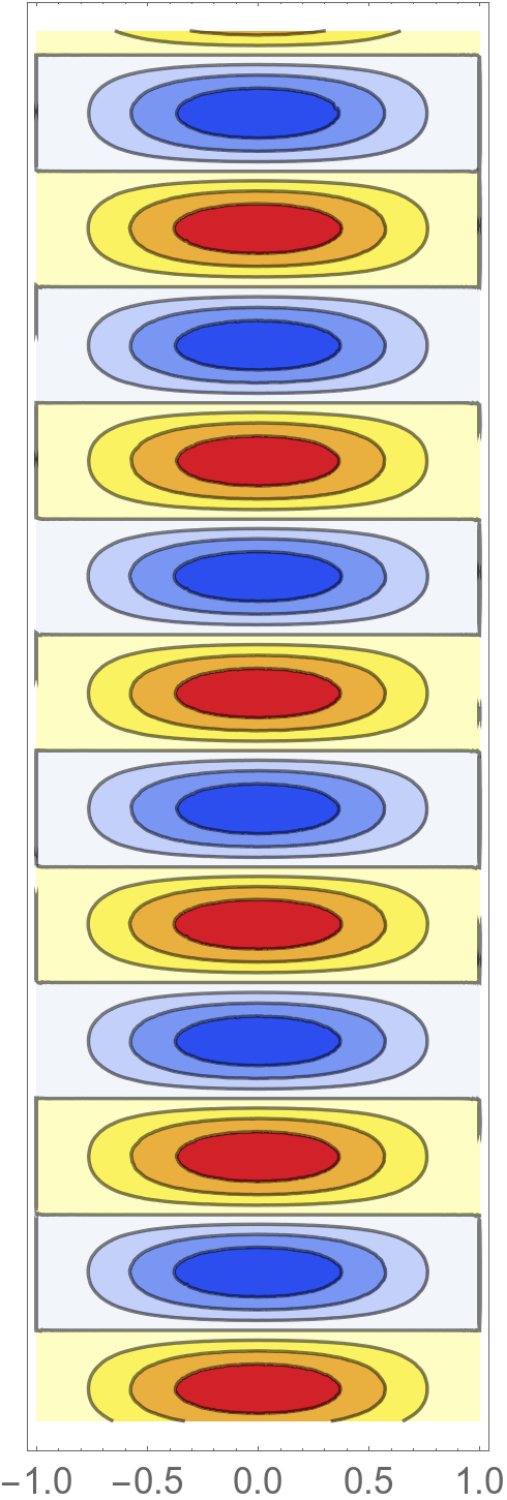}}\quad
\subfloat[{$\Var[u]$.}]{\includegraphics[width=0.17\textwidth]{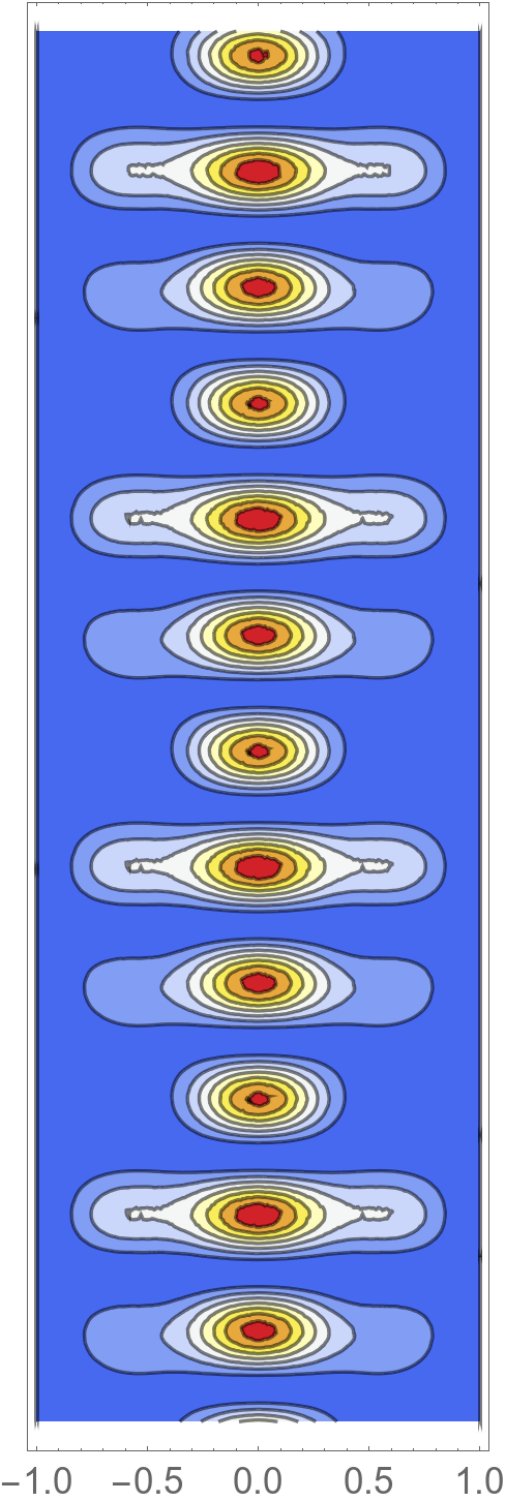}}\quad
\subfloat[{$\Var[v]$.}]{\includegraphics[width=0.17\textwidth]{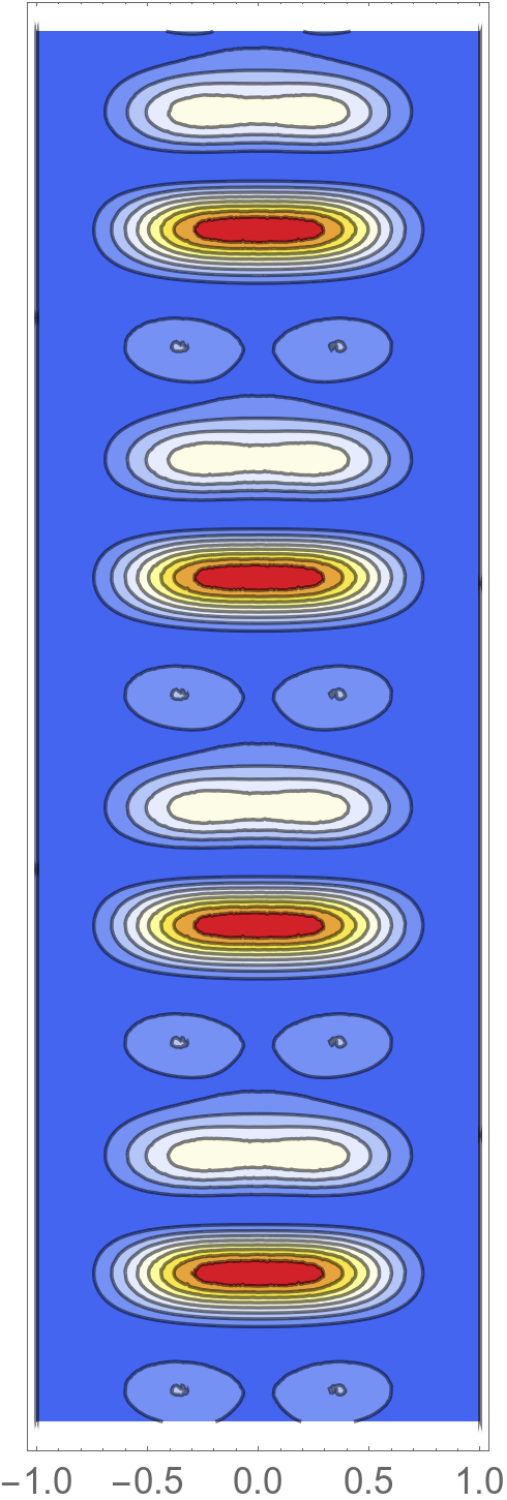}}\quad
\subfloat[{$\Var[w]$.}]{\includegraphics[width=0.17\textwidth]{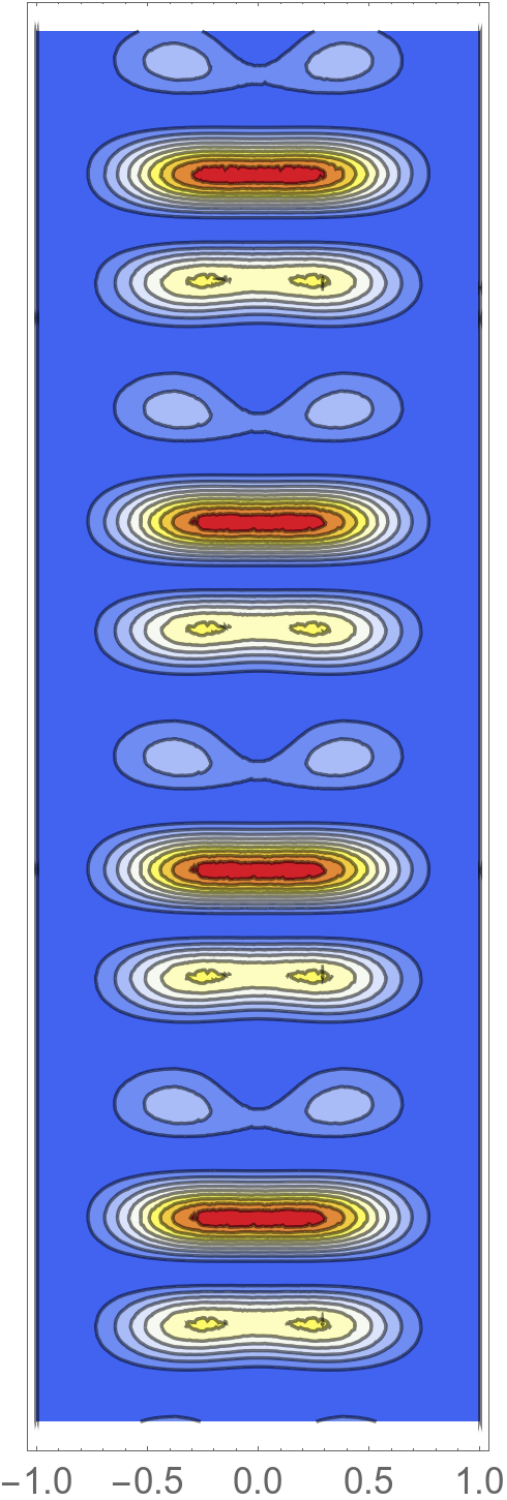}}\quad
\subfloat[{$\Var[\theta]$.}]{\includegraphics[width=0.17\textwidth]{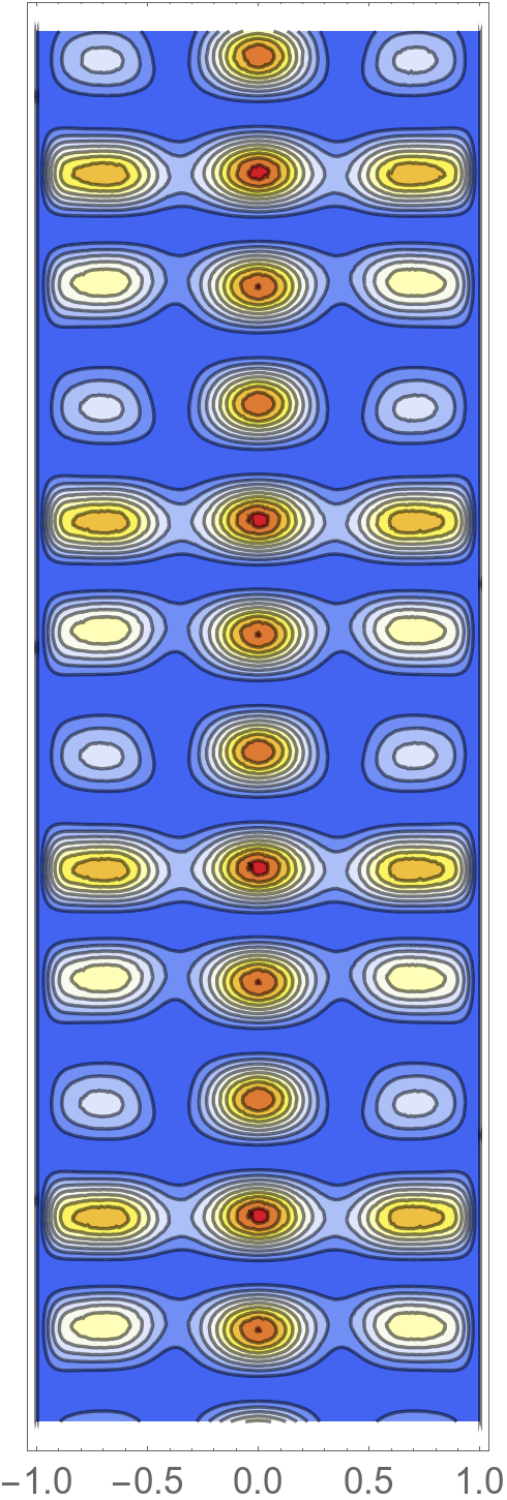}}\quad
\subfloat[{$\Var[\psi]$.}]{\includegraphics[width=0.17\textwidth]{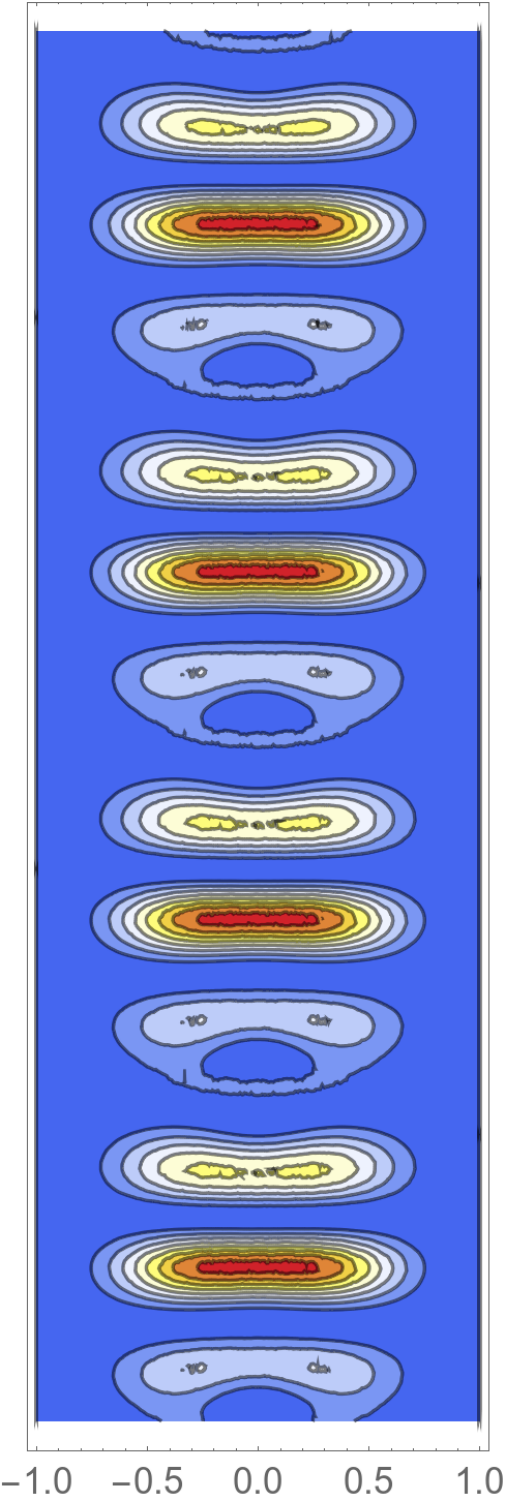}}
\caption{The solution for $t=1/100$ in the case of general uncertainty and the Naghdi shell model. Expected value and variance of the first basis function of the eigenspace. A contour plot of the five field components.}
\label{fig:generalnaghdicomps}
\end{center}
\end{figure}

\begin{figure}[p!]
\begin{center}
\subfloat[{$\E[u]$.}]{\includegraphics[width=0.17\textwidth]{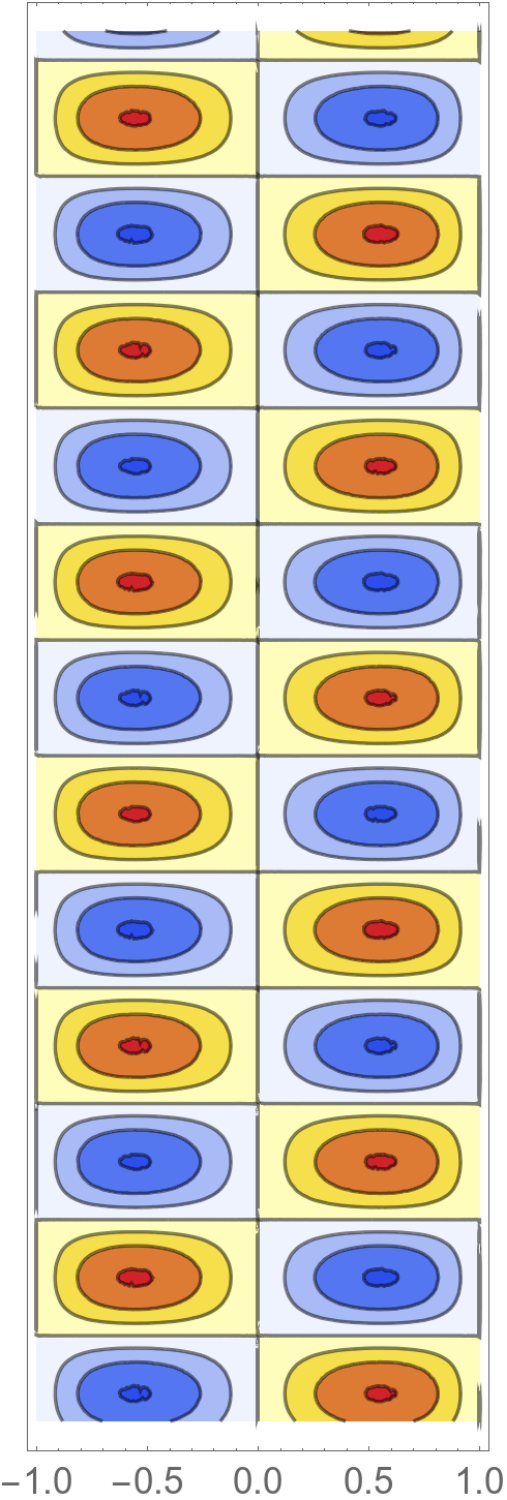}}\quad
\subfloat[{$\E[v]$.}]{\includegraphics[width=0.17\textwidth]{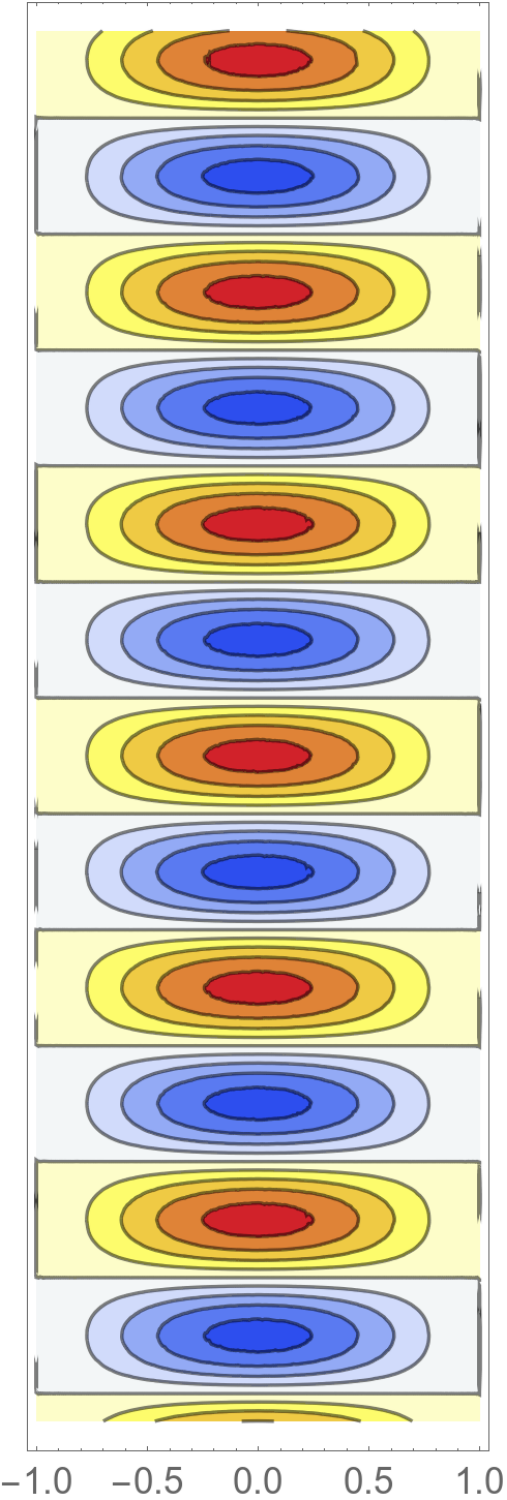}}\quad
\subfloat[{$\E[w]$.}]{\includegraphics[width=0.17\textwidth]{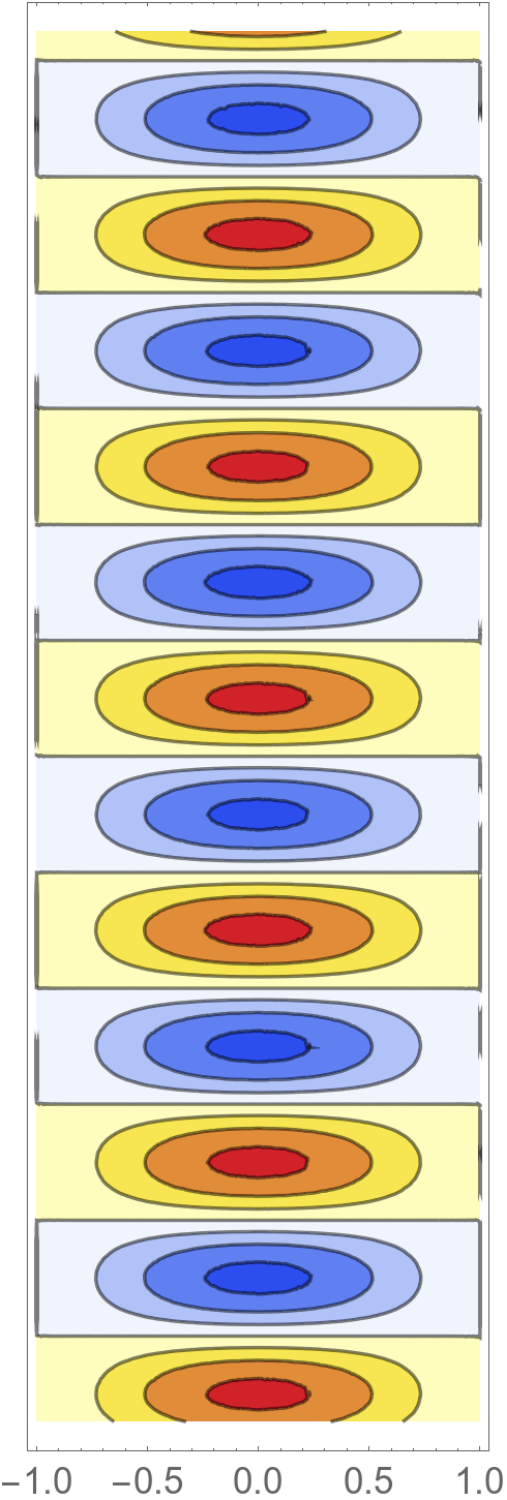}}\quad
\subfloat[{$\E[\theta]$.}]{\includegraphics[width=0.17\textwidth]{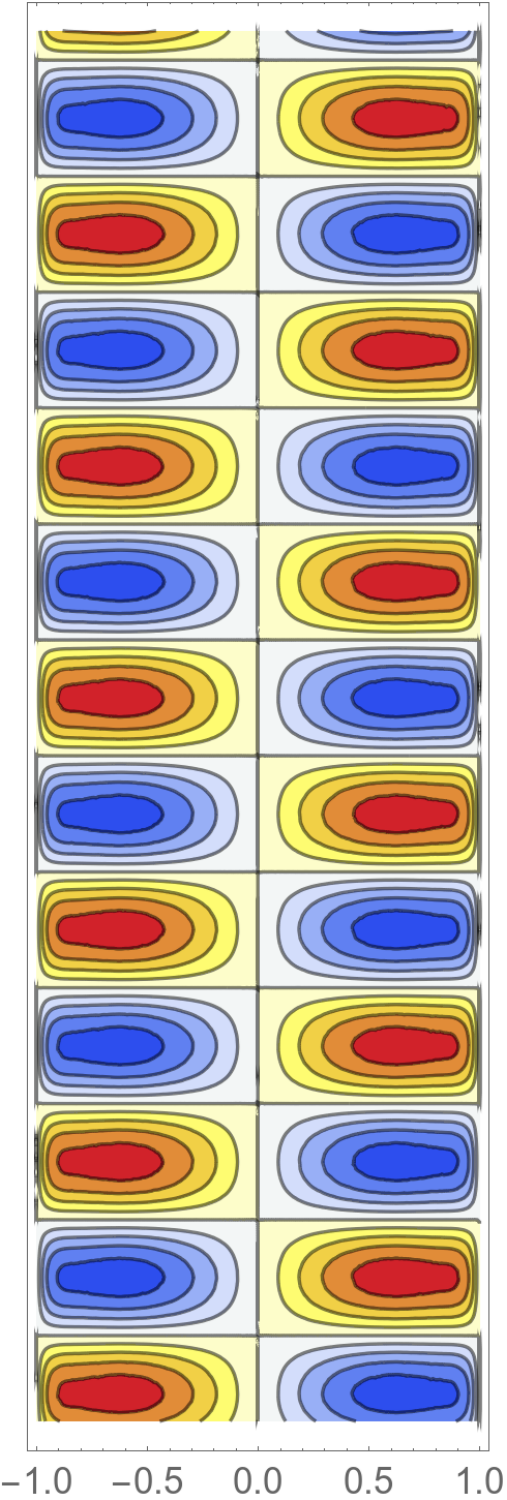}}\quad
\subfloat[{$\E[\psi]$.}]{\includegraphics[width=0.17\textwidth]{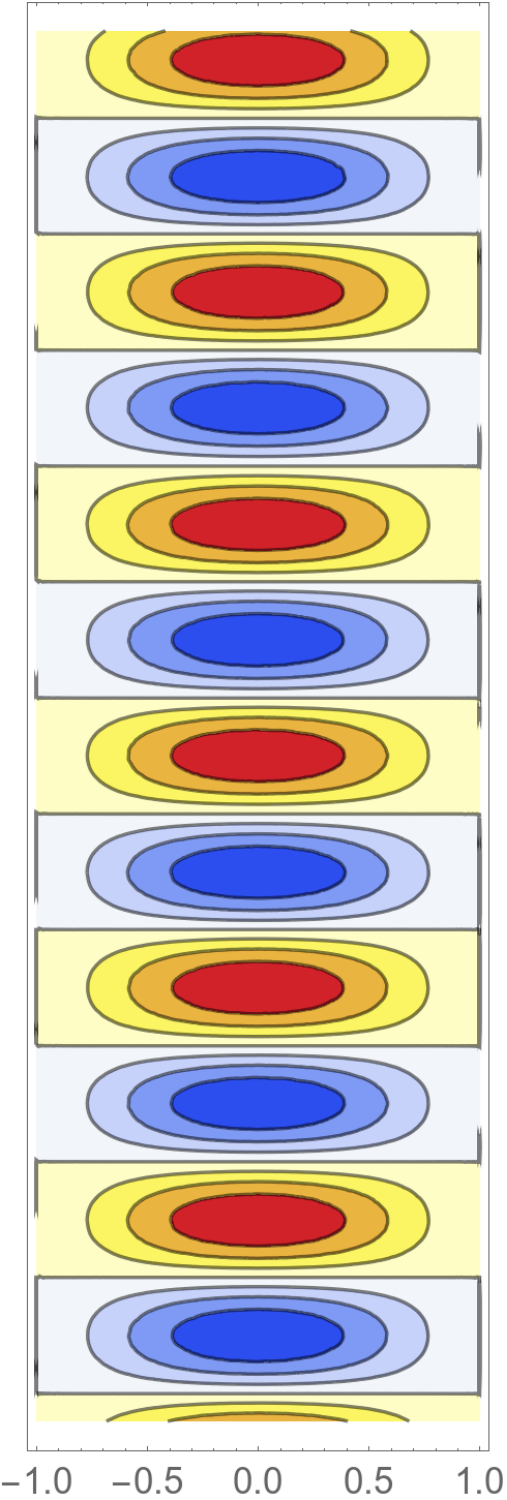}}\quad
\subfloat[{$\Var[u]$.}]{\includegraphics[width=0.17\textwidth]{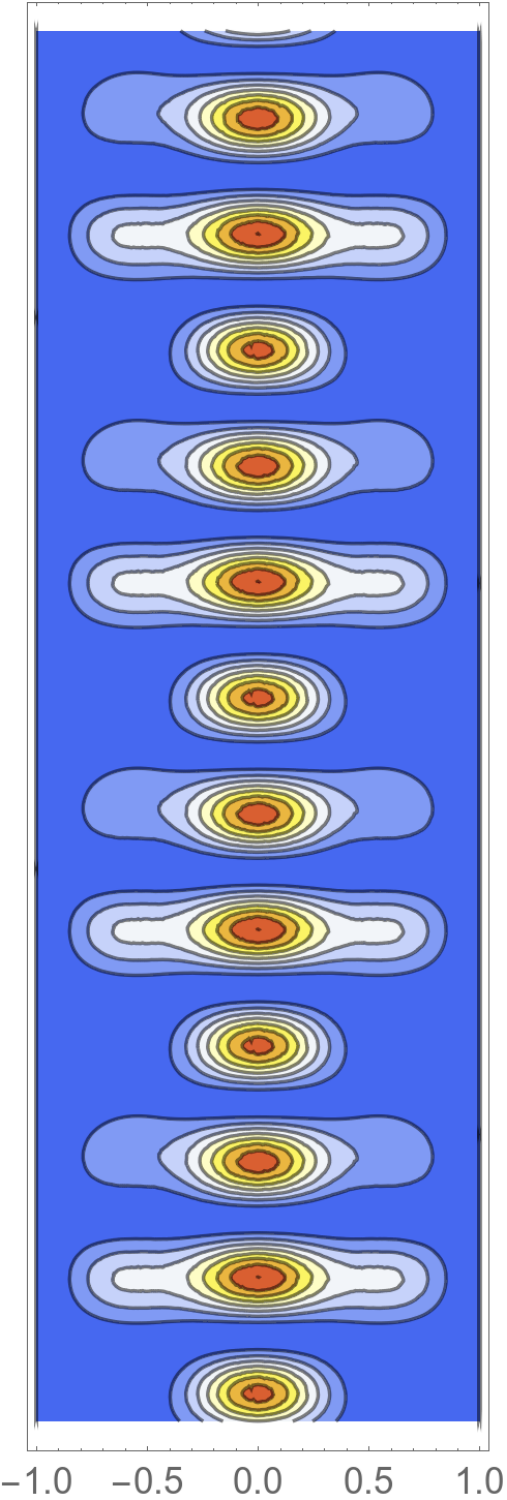}}\quad
\subfloat[{$\Var[v]$.}]{\includegraphics[width=0.17\textwidth]{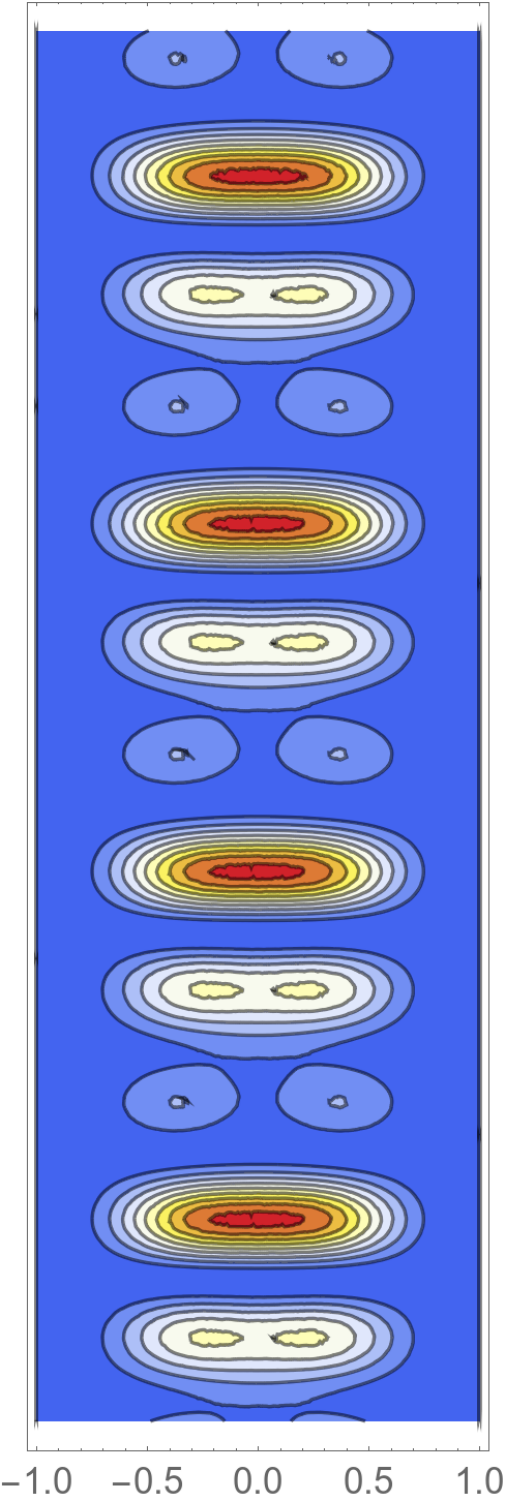}}\quad
\subfloat[{$\Var[w]$.}]{\includegraphics[width=0.17\textwidth]{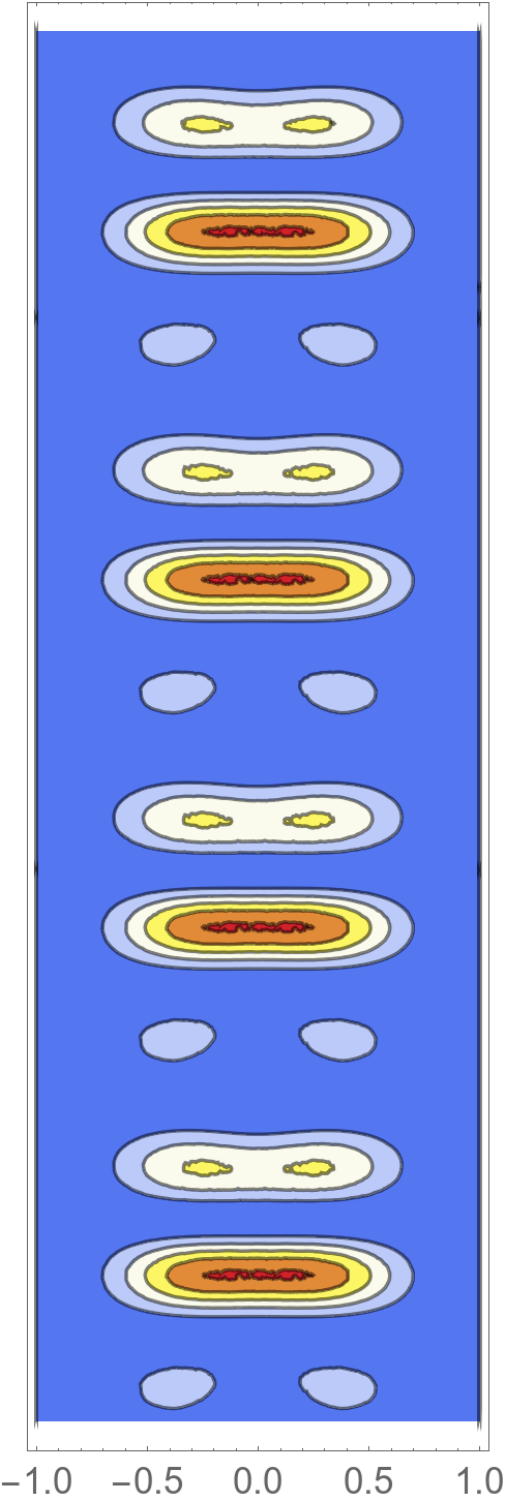}}\quad
\subfloat[{$\Var[\theta]$.}]{\includegraphics[width=0.17\textwidth]{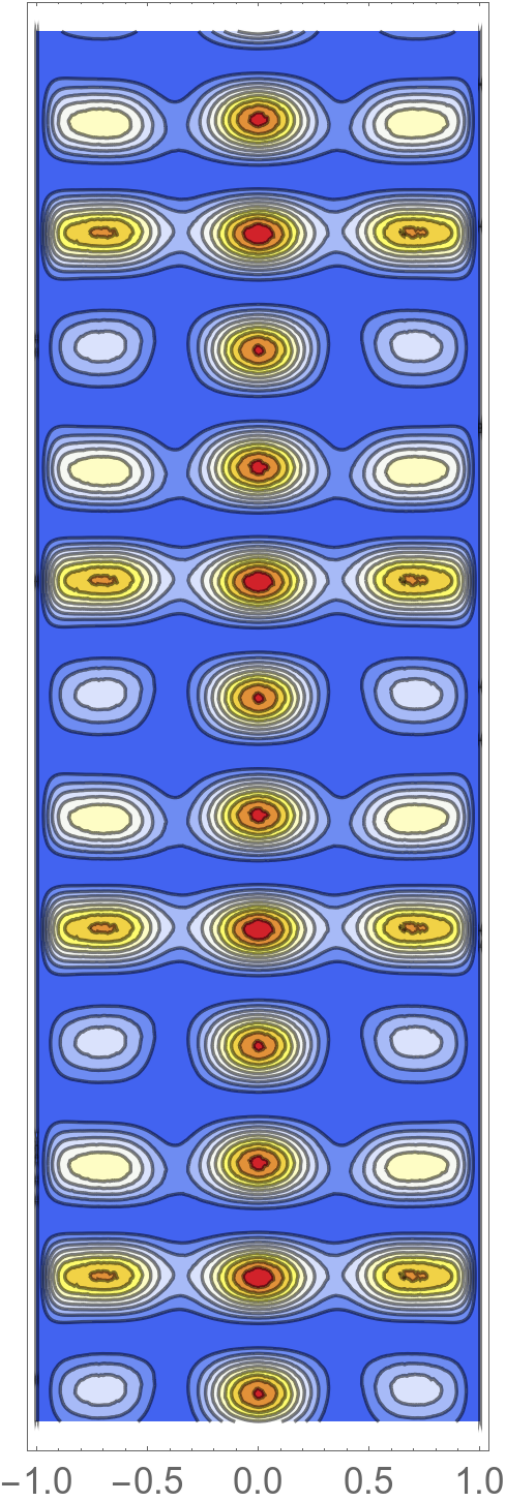}}\quad
\subfloat[{$\Var[\psi]$.}]{\includegraphics[width=0.17\textwidth]{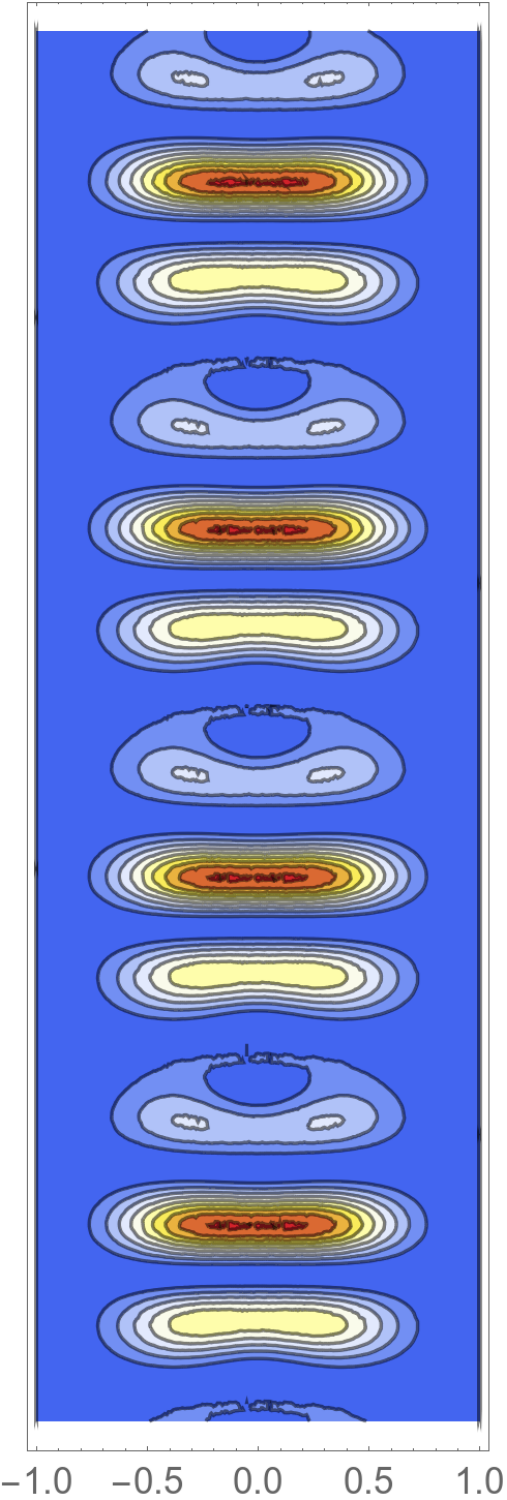}}
\caption{The solution for $t=1/100$ in the case of general uncertainty and the mathematical shell model. Expected value and variance of the first basis function of the eigenspace. A contour plot of the five field components.}
\label{fig:generalshallowcomps}
\end{center}
\end{figure}

Finally we assume the general model of uncertainty \eqref{general}. We set $t=1/100$ and compute the two-dimensional eigenspace associated to the two smallest eigenvalues using the 2D solver. Statistics of the first basis function have been presented in Figures \ref{fig:generalnaghdicomps} and \ref{fig:generalshallowcomps} for the Naghdi and mathematical shell models respectively. The axial and angular components of the solution no longer separate and therefore the 1D ansatz in \ref{sec:1D_models} is not valid in the case of general uncertainty.

\begin{remark}
As expected the results for the Naghdi and mathematical shell models differ slightly. There are at least two reasons for this: not only are the underlying mathematical models different but also the choice of the basis for the two-dimensional eigenspace may differ between the two examples. The Figures \ref{fig:generalnaghdicomps} and \ref{fig:generalshallowcomps} only illustrate the first basis function of the eigenspace.
\end{remark}

\section{Conclusions}
\label{sec:conclusions}

Shells of revolution have natural eigenclusters due to symmetries,
moreover, the eigenpairs depend on a deterministic parameter, the dimensionless thickness.
The stochastic subspace iteration algorithms presented here are capable 
of resolving the eigenclusters.

It is interesting that the crossing of eigenpairs is not only a theoretical
concept but something that can be demonstrated in a standard engineering problem.
We emphasise that the effect is not an artefact of our particular numerical experiment.

Not surprisingly, the effect of the chosen material model on the
asymptotics in relation to the deterministic parameter is shown to be negligible.
Both the smallest eigenvalue and its standard deviation have the same dependence on the  
dimensionless thickness. The logical next step would be to consider problems with 
varying thickness. This is not straightforward, however, since the dimension reduction 
assumes constant thickness. 

\appendix

\section{Shell Models}
\label{sec:shell_models}

In this appendix the two shell models used in the experiments are discussed.
After the introduction of the two-dimensional models, the corresponding one-dimensional
ones for the shells of revolution are derived by applying suitable ansatz.

\subsection{Shell Geometry}\label{sec:shell_geometry}
In this work we study thin shells of revolution.
They can formally be characterized as domains in $\mathbb{R}^3$ of type
\begin{equation}
	\Omega = \left\{ \mathbf{x} + z \mathbf{n(x)} \ | \ \mathbf{x} \in D, -d/2 < z < d/2\right\},
\end{equation}
where $d$ is the (constant) thickness of the shell, 
$D$ is a (mid)surface of revolution, 
and $\mathbf{n(x)}$ is the unit normal to $D$.
For realistic geometries we assume principal curvature coordinates, where only four parameters,
the radii of principal curvature $R_1$, $R_2$, and the so-called Lam\'e parameters,
$A_1$, $A_2$, which relate coordinates changes to arc lengths, are needed to specify
the curvature and the metric on $D$. The displacement vector field of
the midsurface $\mathbf{u} = \{u,v,w\}$ can be interpreted as projections to
directions
\begin{equation}
\mathbf{e}_1 = \frac{1}{A_1}\frac{\partial \Psi}{\partial x_1},\quad \mathbf{e}_2 = \frac{1}{A_2}\frac{\partial \Psi}{\partial x_2},\quad \mathbf{e}_3 = \mathbf{e}_1\times \mathbf{e}_2,
\end{equation}
where $\Psi(x_1,x_2)$ is a suitable parametrization of the surface of revolution,
$\mathbf{e}_1, \mathbf{e}_2$ are the unit tangent vectors along the principal curvature lines, and
$\mathbf{e}_3$ is the unit normal.
In other words 
\[
\mathbf{u} = u\,\mathbf{e}_1+ v\,\mathbf{e}_2+ w\,\mathbf{e}_3.
\]
There is another option, however, which makes it possible to consider \textit{non-realisable}
shell geometries. We can simplify the model above by assuming that $D$ is a rectangular 
domain expressed in the coordinates $x_1$ and $x_2$. Furthermore, we assume that the curvature 
tensor $\{b_{ij}\}$ of the midsurface is constant and write $a = b_{11}, b = b_{22}$, and 
$c = b_{12} = b_{21}$.

Let us consider a cylindrical shell generated by a function $f_1(x_1) = 1$, $x_1 \in [-x_0,x_0]$, $x_0 > 0$.
In this case the product of the Lam\'e parameters (metric), $A_1(x_1) A_2(x_1) = 1$, and the reciprocal
curvature radii are $1/R_1(x_1) = 1$ and $1/R_2(x_1) = 0$, since 
\begin{equation}\label{eq:lame}
A_1(x_1) = \sqrt{1+[f_1 '(x_1)]^2}, \quad A_2(x_1) = f_1 (x_1),
\end{equation}
and
\begin{equation}\label{eq:curvature}
R_1(x_1) = - \frac{A_1(x_1)^3}{f ''(x_1)}, \quad R_2(x_1) = A_1(x_1) A_2(x_1).
\end{equation}
Thus, in the simplified model we can choose
$a=0$, $b=1$, and $c=0$, and arrive at a very good approximation of the exact geometry.
On the other hand, if we consider $f_2(x_1) = 1 + x_1^2$, it is clear that the approximation by
constant curvature tensor degenerates as $x_0$ increases.

\subsection{Two-Dimensional Models}\label{sec:2D_models}
Our two-dimensional shell models are the so-called Reissner-Naghdi model \cite{malinen},
and the mathematical shell model \cite{pms},
where the transverse deflections are approximated with low-order
polynomials. The resulting vector field has five components
$\u=(u,v,w,\theta,\psi)$, where the first three are the
displacements and the latter two are the rotations
in the axial and angular directions, respectively.
Here we adopt the convention that the computational domain $D$ 
is given by the surface parametrisation and the axial/angular coordinates are denoted by $x$ and $y$. 



Deformation energy $\A(\u,\u)$ is divided into bending, membrane, and shear energies,
denoted by subscripts $B$, $M$, and $S$, respectively.
\begin{equation}
\A(\u,\u) = d^2 \A_B(\u,\u)+ \A_M(\u,\u)+ \A_S(\u,\u).
\label{eqn:muodonmuutosenergia}
\end{equation}

Bending, membrane, and shear energies are given as
\begin{eqnarray}
d^2 \A_{B}(\u,\u)
&=& d^2 \int_{\w} E(x,y) \Bigl[\nu (\k_{11}(\u)+\k_{22}(\u))^2 \nonumber \\
& & +(1-\nu) \sum_{i,j=1}^2 \k_{ij}(\u)^2 \Bigr] \ A_1(x,y) A_2(x,y) \ dx\, dy, \label{eqn:taipumaenergia} \\
\A_{M}(\u,\u) &=& 12 \int_{\w}E(x,y) 
    \Bigl[ \nu (\b_{11}(\u)+\b_{22}(\u))^2 \nonumber\\
& & + (1-\nu) \sum_{i,j=1}^2 \b_{ij}(\u)^2\Bigr]
    \ A_1(x,y) A_2(x,y) \ dx\, dy, \label{eqn:kalvoenergia}\\
\A_{S}(\u,\u) &=&
    6(1-\nu)\int_{\w}E(x,y) 
    \Bigl[(\r_1(\u)^2 + \r_2(\u))^2\Bigr] \times \nonumber\\
 & &   \ A_1(x,y) A_2(x,y) \ dx\, dy, \label{eqn:leikkausenergia}
\end{eqnarray}
where $\nu$ is the Poisson ratio (constant), and $E(x,y)$ is the Young's modulus
with scaling $1/(12(1-\nu^2))$.

\subsection{One-Dimensional Models}\label{sec:1D_models}
The shell models above can be further reduced to one-dimensional ones.
For shells of revolution the eigenmodes $\u(x,y)$ have either one the forms
\[
\u_1(x,y) = \left(
  \begin{array}{c}
    u(x) \cos(k \, y) \\
    v(x) \sin(k \, y) \\
    w(x) \cos(k \, y) \\
    \theta(x) \cos(k \, y) \\
    \psi(x) \sin(k \, y) \\
  \end{array}
\right), \quad
    \u_2(x,y) = \left(
      \begin{array}{c}
        u(x) \sin(k \, y) \\
        v(x) \cos(k \, y) \\
        w(x) \sin(k \, y) \\
        \theta(x) \sin(k \, y) \\
        \psi(x) \cos(k \, y) \\
      \end{array}
    \right).
\]

Using the ansatz above the energies can be written in terms of
the harmonic number $k$:
\begin{eqnarray}
d^2 \A_{B}^{1D}(\u,\u)
&=& d^2 \int_{0}^1 E(x) \Bigl[\nu (\k_{11}(\u)+\k_{22}(\u))^2 \nonumber \\
& & +(1-\nu) \sum_{i,j=1}^2 \k_{ij}(\u)^2 \Bigr] \ A_1(x) A_2(x) \ dx , \label{eqn:taipumaenergia1D} \\
\A_{M}^{1D}(\u,\u) &=& 12 \int_{0}^1 E(x) 
    \Bigl[ \nu (\b_{11}(\u)+\b_{22}(\u))^2 \nonumber\\
& & + (1-\nu) \sum_{i,j=1}^2 \b_{ij}(\u)^2\Bigr]
    \ A_1(x) A_2(x) \ dx, \label{eqn:kalvoenergia1D}\\
\A_{S}^{1D}(\u,\u) &=&
    6(1-\nu)\int_{0}^1 E(x)
    \Bigl[(\r_1(\u)^2 + \r_2(\u))^2\Bigr] \times \nonumber\\
   & &  \ A_1(x) A_2(x) \ dx \label{eqn:leikkausenergia1D},
\end{eqnarray}

\subsection{2D Strains}\label{sec:2dstrains}
In this section the 2D strains are presented. Here the choice
$f_1(x_1) = 1$, $x_1 \in [-x_0,x_0]$, $x_0 > 0$, is explicit in the Naghdi
strains by (\ref{eq:lame}) and {\ref{eq:curvature}). Indeed, with the choice $a=0$, $b=1$, and $c=0$ the models
are identical except for the shear strain $\rho_2$.

\subsubsection{Naghdi Shell Model (2D)}
\begin{equation}
    \begin{aligned}
        \k_{11} &= \frac{\D \theta}{\D x},\quad
        \k_{22} = \frac{\D \psi}{\D y}, \quad
        \k_{12}
                 =  \frac{1}{2}\left(
                    \frac{\D \psi}{\D x}
                 +  \frac{\D \theta}{\D y}
        -  \frac{\D v}{\D x}
               \right), \\
        \b_{11} &=  \frac{\D u}{\D x},\quad
    \b_{22} =  \frac{\D v}{\D y}+w,  \quad
    \b_{12} 
             =  \frac{1}{2}\left (
                 \frac{\D v}{\D x}
              +  \frac{\D u}{\D y}
             \right), \\
    \r_{1}  &=  \frac{\D w}{\D x}
              -\theta, \quad
    \r_{2}  =  \frac{\D w}{\D y}
              - v -\psi.
    \end{aligned}
\end{equation}

    
\subsubsection{Mathematical Shell Model (2D)}
\begin{equation}
\begin{aligned}
    \kappa_{11} &= \frac{\partial \theta}{\partial x}, \quad \kappa_{22} = \frac{\partial \psi}{\partial y},
\quad \kappa_{12} = \frac{1}{2}\left(\frac{\partial \theta}{\partial y} + \frac{\partial \psi}{\partial x}\right), \\
\beta_{11} &= \frac{\partial u}{\partial x} + a w,\quad \beta_{22} = \frac{\partial v}{\partial y} + b w,
\quad \beta_{12} = \frac{1}{2}\left(\frac{\partial u}{\partial y} + \frac{\partial v}{\partial x}\right) + c w,\\
\rho_1 &=  \frac{\partial w}{\partial x}- \theta, \quad \rho_2 = \frac{\partial w}{\partial y} - \psi.
\end{aligned}
\end{equation}

\subsection{1D Strains}\label{sec:1dstrains}
From the 2D strains one can derive the 1D strains using either of the ansatz given above.
Here we have used $\u_1(x,y)$.

\subsubsection{Naghdi Shell Model (1D)}
\begin{equation}
    \begin{aligned}
        \k_{11} &= \frac{\D \theta}{\D x},  \quad
        \k_{22} = k \, \psi,  \quad
        \k_{12} = \frac{1}{2}\left(
                    \frac{\D \psi}{\D x}
                 -  k \, \theta
                 -  \frac{\D v}{\D x}
                \right),  \\
        \b_{11} &=  \frac{\D u}{\D x},\quad 
        \b_{22} = k \, v + w,  \quad
        \b_{12} = \frac{1}{2}\left(
                     \frac{\D v}{\D x}
                  - k \, u
               \right), \\
        \r_{1}  &=  \frac{\D w}{\D x} -\theta, \quad
        \r_{2}  = - k \, w - v -\psi.
    \end{aligned}
\end{equation}
    
    
\subsubsection{Mathematical Shell Model (1D)}
\begin{equation}
\begin{aligned}
    \kappa_{11} &= \frac{\partial \theta}{\partial x}, \quad \kappa_{22} = k \psi,
\quad \kappa_{12} = \frac{1}{2}\left(-k \theta + \frac{\partial \psi}{\partial x}\right), \\
\beta_{11} &= \frac{\partial u}{\partial x} + a w,\quad \beta_{22} = k v + b w,
\quad \beta_{12} = \frac{1}{2}\left(-k u + \frac{\partial v}{\partial x}\right) + c w,\\
\rho_1 &= \frac{\partial w}{\partial x} - \theta, \quad \rho_2 = -k w - \psi.
\end{aligned}
\end{equation}

\bibliographystyle{model1-num-names}
\let\v\relax
\bibliography{ShellSFEM}
\end{document}